\documentclass[a4paper,11pt]{article}

\usepackage{geometry}
\usepackage{amsfonts}
\usepackage{amsthm}
\usepackage{amssymb}
\usepackage{amsmath}
\usepackage{upref}
\usepackage{mathrsfs}
\usepackage{color}
\usepackage[citecolor=blue,colorlinks=true]{hyperref}
\usepackage{enumitem}
\usepackage{graphicx}
\usepackage{tikz}

\theoremstyle{plain}
\newtheorem{theorem}{Theorem}[section]
\newtheorem{proposition}[theorem]{Proposition}
\newtheorem{lemma}[theorem]{Lemma}
\newtheorem{corollary}[theorem]{Corollary}

\theoremstyle{definition}
\newtheorem{definition}{Definition}[section]

\theoremstyle{remark}
\newtheorem{remark}{Remark}[section]

\numberwithin{equation}{section}
\def\N{\mathbb{N}}

\def\R{\mathbb{R}}

\def\ds{\displaystyle} 
\def\div{{\rm div}}
\def\ocirc#1{\ifmmode\setbox0=\hbox{$#1$}\dimen0=\ht0
    \advance\dimen0 by1pt\rlap{\hbox to\wd0{\hss\raise\dimen0
    \hbox{\hskip.2em$\scriptscriptstyle\circ$}\hss}}#1\else
    {\accent"17 #1}\fi} 
\def\sgn{\mathrm{sgn}}
\def\eps{\varepsilon}
\def\<{\langle}
\def\>{\rangle}
\def\F{\mathcal{F}}
\def\GG{\mathbf{G}}
\def\P{\mathbb{P}}
\def\E{\mathbb{E}}
\def\T{\mathbb{T}}

\def\order{\mathrm{d}}

\begin{document}

\title{Convergence of approximations to stochastic scalar conservation laws}
\author{Sylvain Dotti\thanks{Aix-Marseille Universit\'e, CNRS, Centrale Marseille, I2M, UMR 7373, 13453 Marseille France }{ } and Julien Vovelle\thanks{Univ Lyon, Universit\'e Claude Bernard Lyon 1, CNRS UMR 5208, Institut Camille Jordan, 43 blvd. du 11 novembre 1918, F-69622 Villeurbanne cedex, France. Julien Vovelle was supported by ANR projects STOSYMAP and STAB.}}
\maketitle

\begin{abstract} We develop a general framework for the analysis of approximations to stochastic scalar conservation laws. Our aim is to prove, under minimal consistency properties and bounds, that such approximations  are converging to the solution to a stochastic scalar conservation law. The weak probabilistic convergence mode is convergence in law, the most natural in this context. We use also a kinetic formulation and martingale methods. Our result is applied to the convergence of the Finite Volume Method in the companion paper \cite{DottiVovelle16b}. 
\end{abstract}

\medskip 


\tableofcontents

\section{Introduction}\label{secIntro}

Let $(\Omega,\F,\P,(\F_t),(\beta_k(t)))$ be a stochastic basis and let $T>0$. Consider the first-order scalar conservation law with stochastic forcing
\begin{equation}\label{stoSCL}
du(x,t)+\div(A(u(x,t)))dt=\Phi(x,u(x,t)) dW(t),\quad x\in\T^N, t\in(0,T).
\end{equation}
Equation~\eqref{stoSCL} is periodic in the space variable:  $x\in\T^N$ where $\T^N$ is the $N$-dimensional torus. The flux function $A$ in \eqref{stoSCL} is supposed to be of class $C^2$: $A\in C^2(\R;\R^N)$. We assume that $A$ and its derivatives have at most polynomial growth. The right-hand side of \eqref{stoSCL} is a stochastic increment in infinite dimension. It is defined as follows (see \cite{DaPratoZabczyk92} for the general theory): $W$ is a cylindrical Wiener process, $W=\sum_{k\geq 1}\beta_k e_k$, where the coefficients $\beta_k$ are independent Brownian processes and $(e_k)_{k\geq 1}$ is a complete orthonormal system in a Hilbert space $H$. For each $x\in\T^N$, $u\in\R$, $\Phi(x,u)\in L_2(H,\R)$ is defined by $\Phi(x,u)e_k=g_k(x,u)$ where $g_k(\cdot,u)$ is a regular function on $\T^N$. Here, $L_2(H,K)$ denotes the set of Hilbert-Schmidt operator from the Hilbert space $H$ to an other Hilbert space $K$. Since $K=\R$ in our case, this set is isomorphic to $H$, thus we may also define
\begin{equation}\label{defPhigk}
\Phi(x,u)=\sum_{k\geq 1}g_k(x,u) e_k,
\end{equation}
the action of $\Phi(x,u)$ on $e\in H$ being given by $\<\Phi(x,u),e\>_H$. We assume $g_k\in C(\T^N\times\R)$, with the bounds 
\begin{align}
\GG^2(x,u)=\|\Phi(x,u)\|_H^2=\sum_{k\geq 1}|g_k(x,u)|^2\leq D_0(1+|u|^2),\label{D0}\\
\|\Phi(x,u)-\Phi(y,v)\|_H^2=\sum_{k\geq 1}|g_k(x,u)-g_k(y,v)|^2\leq D_1(|x-y|^2+|u-v|h(|u-v|)),\label{D1}
\end{align}
where $x,y\in\T^N$, $u,v\in\R$, and $h$ is a continuous non-decreasing function on $\R_+$ such that $h(0)=0$. We assume also $0\leq h(z)\leq 1$ for all $z\in\R_+$. \bigskip

\textit{Notation:} in what follows, we will use the convention of summation over repeated indices $k$. For example, we write $W=\beta_k e_k$ for the cylindrical Wiener process in \eqref{stoSCL}.
\bigskip

This paper is a preliminary work to the analysis of convergence of the numerical approximation to \eqref{stoSCL} by the Finite Volume method with monotone fluxes, which is done in \cite{DottiVovelle16b}. We give a general notion of family of approximate solutions, see Definition~\ref{def:appsol}, and explain what kind of convergence of such family can be expected. Our main results in this regard are the theorem~\ref{th:martingalesol}, about convergence to martingale solutions, and the theorem~\ref{th:pathcv}, which gives criteria for convergence to pathwise solutions.\bigskip

Problem~\eqref{stoSCL} has already been studied in a series of papers. Like in the deterministic case, the approach to the existence of solutions has been the vanishing viscosity method, see  \cite{EKMS00,Kim03,FengNualart08,ValletWittbold09,DebusscheVovelle10,ChenDingKarlsen12,BauzetValletWittbold12,
BauzetValletWittbold14,KarlsenStorrosten2017} in particular. Approximation by the BGK method has been considered by M.~Hofmanov{\'a} in \cite{HofmanovaBGK15}. Some results of convergence of numerical approximations to \eqref{stoSCL} (by the Finite Volume method in particular) have also be obtained in \cite{KrokerRohde12,BauzetCharrierGallouet14a,BauzetCharrierGallouet14b,Bauzet15,StorrostenKarlsen2016,KoleyMajeeVallet2016}.
\bigskip

The main difference between this present paper and all the works cited above is in the way to answer to the following question: when considering the convergence of approximations to \eqref{stoSCLi}, which mode of convergence regarding the sample variable $\omega$ is used?  Here, we develop an approach based on convergence in law, while in the work referred to\footnote{with the exception of \cite{FengNualart08}, where quite a strong notion of solution is used however}, weak convergence (in Lebesgue spaces, or in the sense of Young measures, \textit{cf.} Section~\ref{sec:ges}) is considered. Convergence in law is the natural mode of convergence for the random variables which manifest in the approximation to \eqref{stoSCL}. Our approach based on convergence in law is successful because we work in the context of c{\`a}dl{\`a}g processes. This is an other difference between this present paper and the references already quoted: our formulation of solution is weak in the space variable, but not weak in the time variable, see \eqref{eq:kineticupre}, \eqref{eq:kineticfpre} for example. This allows to obtain convergence of approximation for each time $t$ (this is the last statement in Theorem~\ref{th:pathcv}), without making any regularity hypothesis on the initial datum at any moment. This paper is also a further development of the approach by kinetic formulation initiated in \cite{DebusscheVovelle10}. We need it crucially in the companion paper \cite{DottiVovelle16b} to obtain the convergence of The Finite Volume method with a standard CFL condition (\textit{cf.} our comment on the Kinetic formulation in the introduction section of \cite{DottiVovelle16b}).
\bigskip

To complete this introduction, let us mention that the approximation of scalar conservation laws with stochastic flux has also been considered in \cite{GessPerthameSouganidis2016} (time-discrete scheme) and
\cite{MohamedSeaidZahri13} (space discrete scheme). For the corresponding Cauchy Problem, see \cite{LionsPerthameSouganidis13a,LionsPerthameSouganidis13,LionsPerthameSouganidis14,
GessSouganidis2015,GessSouganidis2014,Hofmanova2016}.
\bigskip

The plan of the paper is the following one: Section~\ref{sec:kisol} to Section~\ref{sec:martingale} are devoted to the analysis of the Cauchy Problem for \eqref{stoSCL}: we introduce the kinetic formulation of the problem in Section~\ref{sec:kisol}, and prove a uniqueness result in Section~\ref{sec:comparison}. In Section~\ref{sec:martingale}, we develop a general approach to the analysis of convergence of approximate solutions to \eqref{stoSCL} based on martingale methods. Note that Section~\ref{sec:kisol} and Section~\ref{sec:comparison} are for a large part identical to Section~2 and Section~3 in \cite{DebusscheVovelle10}. There are however a lot of modifications, which were needed to prepare Section~\ref{sec:martingale}. In Section~\ref{sec:appli}, we give some applications of our results of convergence of approximation.



\section{Kinetic solution}\label{sec:kisol}


\subsection{Definition}

\subsubsection{Predictable sets and functions}\label{sec:predictable}

For $T>0$, we denote by $\mathcal{B}([0,T])$ the Borel $\sigma$-algebra on $[0,T]$ and we denote by $\mathcal{P}_T\subset \mathcal{B}([0,T])\otimes\mathcal{F}$ the predictable $\sigma$-algebra, \cite[Section 2.2]{ChungWilliams90}. If $E$ is a Banach space, a process $(f(t))$ with values in $E$ is said to be weakly-predictable if the process $(\<f(t),\varphi\>_{E,E'})$ is predictable for every $\varphi$ in the topological dual $E'$. This is equivalent to say that $f$ is weakly $\mathcal{P}_T$-measurable, in the sense of \cite[Definition 1, p.130]{Yosida80}. Similarly, we can define the notion of strong predictability: the process $(f(t))$ is said to be strongly predictable if there exists a sequence of $E$-valued, $\mathcal{P}_T$-measurable simple functions which converges to $f$ at every point $(t,\omega)$ in a set of full measure in $[0,T]\times\Omega$. By Pettis' Theorem, \cite[Theorem p.131]{Yosida80}, the two notions of measurability coincide if $E$ is separable: in this case we say simply "predictable".\medskip 

Let us assume that $E$ is separable to introduce the following notations. Let $p\in[1,+\infty)$. The set $L^p([0,T]\times\Omega;E)$ is the set of $E$-valued, $\mathcal{B}([0,T])\otimes\mathcal{F}$-measurable, Bochner integrable functions $f$ which satisfy
$$
\iint_{[0,T]\times\Omega}\|f(t,\omega)\|_E^p d(\mathcal{L}\times\P)(t,\omega)<+\infty,
$$
where $\mathcal{L}$ is the Lebesgue measure on $[0,T]$. Equivalently, by definition of the product measure $\mathcal{L}\times\P$,
$$
\E\int_0^T\|f(t)\|_E^p dt<+\infty.
$$
We denote by $L^p_{\mathcal{P}}([0,T]\times\Omega;E)$ the set of functions $g$ in $L^p([0,T]\times\Omega;E)$ which are equal $\mathcal{L}\times\P$-almost everywhere to a predictable function $f$. This is the case if, and only if, $\<g,\varphi\>$ is equal $\mathcal{L}\times\P$-almost everywhere to $\<f,\varphi\>$ for all $\varphi\in E'$ (we use the fact that $E'$ is separable since $E$ is separable), so let us briefly consider the case $E=\R$. The class of processes in $L^p_{\mathcal{P}}([0,T]\times\Omega;\R)$ is analysed in \cite[p.~66]{ChungWilliams90} or \cite[p.~172]{RevuzYor99}. In particular, if $X(t)$ is an adapted process with 
$$
\E\int_0^T |X(t)|^pdt<+\infty,
$$
then $X\in L^p_{\mathcal{P}}([0,T]\times\Omega;\R)$. A progressively measurable process $X$ in $L^p([0,T]\times\Omega;\R)$ also is in $L^p_{\mathcal{P}}([0,T]\times\Omega;\R)$.
\medskip

Let $m\in\N^*$. In the case where $E$ is itself a Lebesgue space $E=L^p(D)$, where $D$ is an open subset of $\R^m$, we have $L^p([0,T]\times\Omega;L^p(D))=L^p(D\times [0,T]\times\Omega)$, where $D\times [0,T]\times\Omega$ is endowed with the product measure $\mathcal{L}_{m+1}\times\P$ ($\mathcal{L}_m$ being the $m$-dimensional Lebesgue measure), see \cite[Section~1.8.1]{Droniou01}. Similarly, we have 
$$
L^p_{\mathcal{P}}([0,T]\times\Omega;L^p(D))=L^p_{\mathcal{P}}(D\times[0,T]\times\Omega),
$$
where $L^p_{\mathcal{P}}(D\times [0,T]\times\Omega)$ is the set of functions in 
$L^p(D\times [0,T]\times\Omega)$ which are equal $\mathcal{L}_m\times\mathcal{L}\times\P$-almost everywhere to a $\mathcal{B}(D)\times\mathcal{P}_T$-measurable function (here $\mathcal{B}(D)$ is the Borel $\sigma$-algebra on $D$). We will apply these results with $D=(0,1)^N$, in which case, by periodic extension, we obtain 
\begin{equation}\label{LpLpLp}
L^p([0,T]\times\Omega;L^p(\T^N))=L^p(\T^N\times [0,T]\times\Omega),
\end{equation}
and similarly for spaces $L^p_\mathcal{P}$. 

\subsubsection{Random measure, solution}\label{sec:kim}

Let $\mathcal{M}_b(\T^N\times[0,T]\times\R)$ be the set of bounded Borel signed measures on $\T^N\times[0,T]\times\R$. We denote by $\mathcal{M}^+_b(\T^N\times[0,T]\times\R)$ the subset of non-negative measures. 

\begin{definition}[Random measure] A map $m$ from $\Omega$ to $\mathcal{M}_b(\T^N\times[0,T]\times\R)$ is said to be a random signed measure (on $\T^N\times[0,T]\times\R$) if, for each $\phi\in C_b(\T^N\times[0,T]\times\R)$, $\<m,\phi\>\colon\Omega\to\R$ is a random variable. If almost-surely $m\in\mathcal{M}^+_b(\T^N\times[0,T]\times\R)$, we simply speak of random measure. 
\label{def:randommeasure}\end{definition}

Let $m$ be a random measure with finite first moment
\begin{equation}\label{Firstm}
\E m(\T^N\times[0,T]\times\R)<+\infty.
\end{equation} 
Then $\E m$ is well defined and this is a bounded measure on $\T^N\times[0,T]\times\R$. In particular, it satisfies the following tightness condition
\begin{equation}
\lim_{R\to+\infty}\E m(\T^N\times[0,T]\times B_R^c)=0,
\label{inftym}
\end{equation}
where $B_R^c=\{\xi\in\R,|\xi|\geq R\}$. We note this fact here, since uniform versions of \eqref{inftym} will be required when considering sequences of random measures, see \eqref{inftymn}. 
We will also need the following result.
\begin{lemma}[Atomic points] Let $m$ be a random measure with first moment~\eqref{Firstm}. Let $\pi\colon\T^N\times[0,T]\times\R\to[0,T]$ denote the projection $(x,t,\xi)\mapsto t$. Let $\pi_\# m$ denote the push-forward of $m$ by $\pi$. Let $B_\mathrm{at}$ denote the set of times $t$ such that the event ``$t$ is an atom of $\pi_\# m$" has positive probability:
\begin{equation}\label{defJstar}
B_\mathrm{at}=\left\{t\in[0,T];\P\left(\pi_\# m(\{t\})>0\right)>0\right\}.
\end{equation}
Then $B_\mathrm{at}$ is at most countable.
\label{lem:atomicpoints}\end{lemma}

\textbf{Proof of Lemma~\ref{lem:atomicpoints}.} We have also
\begin{equation*}
B_\mathrm{at}=\left\{t\in[0,T];\E\pi_\# m(\{t\})>0\right\}.
\end{equation*}
The set $B_\mathrm{at}$ is the set of atomic points of the measure $\E\pi_\# m$. It is therefore at most countable.
\qed\medskip


The notion of solution which we introduce below is based on the kinetic formulation of conservation laws introduced in \cite{LionsPerthameTadmor94}. In particular, for a given function $u$ of the variables $(x,t)$, we will need to consider the function
$$
\mathtt{f}(x,t,\xi):=\mathbf{1}_{u(x,t)>\xi},
$$
which is the characteristic function of the subgraph of $u$. We often write $\mathtt{f}:=\mathbf{1}_{u>\xi}$ for short.\medskip

To be flexible enough, we have to impose a c{\`a}dl{\`a}g property on solutions to \eqref{stoSCL} (see Item~\ref{item:1bisdefsol} in the following Definition~\ref{defkineticsol}). We will show however in Corollary~\ref{cor:timecontinuity} that solutions to \eqref{stoSCL} have continuous trajectories. \medskip

\begin{definition}[Solution] Let $u_{0}\in L^\infty(\T^N)$. An $L^1(\T^N)$-valued stochastic process $(u(t))_{t\in[0,T]}$ is said to be a solution to~\eqref{stoSCL} with initial datum $u_0$ if $u$ and $\mathtt{f}:=\mathbf{1}_{u>\xi}$ have the following properties:
\begin{enumerate}
\item\label{item:1defsol} $u\in L^1_{\mathcal{P}}(\T^N\times [0,T]\times\Omega)$,
\item\label{item:1bisdefsol} for all $\varphi\in C^1_c(\T^N\times\R)$, almost-surely, $t\mapsto\<\mathtt{f}(t),\varphi\>$ is c{\`a}dl{\`a}g,
\item\label{item:2defsol} for all $p\in[1,+\infty)$, there exists $C_p\geq 0$ such that
\begin{equation}
\E\left(\sup_{t\in[0,T]}\|u(t)\|_{L^p(\T^N)}^p\right)\leq C_p,
\label{eq:integrabilityu}\end{equation} 
\item\label{item:3defsol} there exists a random measure $m$ with first moment \eqref{Firstm}, such that for all $\varphi\in C^1_c(\T^N\times\R)$, for all $t\in[0,T]$,
\begin{multline}
\<\mathtt{f}(t),\varphi\>=\<\mathtt{f}_0,\varphi\>
+\int_0^t \<\mathtt{f}(s),a(\xi)\cdot\nabla\varphi\>ds\\
+\sum_{k\geq 1}\int_0^t\int_{\T^N}g_k(x,u(x,s))\varphi(x,u(x,s)) dx d\beta_k(s)\\
+\frac{1}{2}\int_0^t\int_{\T^N} \partial_\xi\varphi(x,u(x,s))\GG^2(x,u(x,s)) dx ds-m(\partial_\xi\varphi)([0,t]),
\label{eq:kineticupre}\end{multline}
a.s., where $\mathtt{f}_0(x,\xi)=\mathbf{1}_{u_0(x)>\xi}$, $\GG^2:=\sum_{k=1}^\infty |g_k|^2$ and $a(\xi):=A'(\xi)$.
\end{enumerate}
\label{defkineticsol}\end{definition}


In \eqref{eq:kineticupre}, we have used the brackets $\<\cdot,\cdot\>$ to denote the duality between $C^\infty_c(\T^N\times\R)$ and the space of distributions over $\T^N\times\R$. In what follows, we will denote similarly the integral 
\begin{equation*}
\<F,G\>=\int_{\T^N}\int_\R F(x,\xi)G(x,\xi) dx d\xi,\quad F\in L^p(\T^N\times\R), G\in L^q(\T^N\times\R),
\end{equation*}
where $1\leq p\leq +\infty$ and $q$ is the conjugate exponent of $p$. 
In \eqref{eq:kineticupre} also, we have used (with $\phi=\partial_\xi\varphi$) the shorthand $m(\phi)$ to denote the Borel measure on $[0,T]$ defined by
\begin{equation*}
m(\phi)\colon A\mapsto\int_{\T^N\times A\times\R} \phi(x,\xi)dm(x,t,\xi),\quad \phi\in C_b(\T^N\times\R),
\end{equation*}
for all $A$ Borel subset of $[0,T]$.\smallskip

There is a last point to comment in Definition~\ref{defkineticsol}, which is the measurability of the function $\sup_{t\in[0,T]}\|u(t)\|_{L^p(\T^N)}$ in \eqref{eq:integrabilityu}. Let us denote by $\bar{\mathtt{f}}=1-\mathtt{f}=\mathbf{1}_{u\leq\xi}$ the conjugate function of $\mathtt{f}$. By the identity
\begin{equation}\label{kitouu}
|u|^p=\int_\R\left[\mathtt{f}\mathbf{1}_{\xi>0}+\bar{\mathtt{f}}\mathbf{1}_{\xi<0}\right] p|\xi|^{p-1}d\xi,
\end{equation}
we have, for $p\in[1,+\infty)$,
\begin{equation}\label{upsupf}
\|u(t)\|_{L^p(\T^N)}^p=\sup_{\psi_+\in F_+,\psi_-\in F_-}\<\mathtt{f}(t),\psi_+\>+\<\bar{\mathtt{f}}(t),\psi_-\>,
\end{equation}
where the sup is taken over some countable sets $F_+$ and $F_-$ of functions $\psi$ chosen as follows: $F_\pm=\{\psi_n;n\geq 1\}$, where $(\psi_n)$ is a sequence of non-negative functions in $C^\infty_c(\R)$ which converges point-wise monotonically to $\xi\mapsto p|\xi^\pm|^{p-1}$ if $p>1$ and to $\xi\mapsto\sgn_\pm(\xi)$ if $p=1$. By \eqref{upsupf}, we have
\begin{equation}\label{supsupcadlag}
\sup_{t\in[0,T]}\|u(t)\|_{L^p(\T^N)}^p=\sup_{\psi\pm\in F_\pm}\sup_{t\in[0,T]}\<\mathtt{f}(t),\psi_+\>+\<\bar{\mathtt{f}}(t),\psi_-\>.
\end{equation}
By Item~\eqref{item:1bisdefsol} in Definition~\ref{defkineticsol}, we know that the function 
$$
\sup_{t\in[0,T]}\<\mathtt{f}(t),\psi_+\>+\<\bar{\mathtt{f}}(t),\psi_-\>
$$ 
is $\F$-mea\-su\-ra\-ble for all $\psi_\pm\in F_\pm$. Indeed, the $\sup$ over $[0,T]$ of a c{\`a}dl{\`a}g function is the $\sup$ of the function on any dense countable subset of $[0,T]$ containing the terminal point $T$. By \eqref{supsupcadlag}, the function $\sup_{t\in[0,T]}\|u(t)\|_{L^p(\T^N)}$ is measurable.

\begin{remark}[Initial condition] A limiting argument based on \eqref{eq:kineticupre} leads to the following initial condition for $\mathtt{f}(t)$:
$$
\mathtt{f}(0)=\mathtt{f}_0+\partial_\xi m_0,\quad a.s.,
$$
where $m_0$ is the restriction of $m$ to $\T^N\times\{0\}\times\R$. It is not obvious thus, that \eqref{eq:kineticupre} entails the expected initial condition $\mathtt{f}(0)=\mathtt{f}_0$. This is the case however (and, therefore, $m_0\equiv 0$ a.s.), due to Proposition~\ref{prop:equilibrium} and Corollary~\ref{cor:equilibrium}. See also the discussion on the same topic in Section~5 of \cite{ChenPerthame03}.
\label{rk:IC}\end{remark}

\begin{proposition}[Mass of the random measure] Let $u_{0}\in L^\infty(\T^N)$. Let $(u(t))_{t\in[0,T]}$ be a solution to \eqref{stoSCL} with initial datum $u_0$. Then the total mass of the measure $m$ in \eqref{eq:kineticupre} is
\begin{multline}\label{eq:massm}
m(\T^N\times[0,T]\times\R)=\frac12\|u_0\|_{L^2(\T^N)}^2-\frac12\|u(T)\|_{L^2(\T^N)}^2\\
+\sum_{k\geq 1}\int_0^T\int_{\T^N}g_k(x,u(x,t)) u(x,t) dx d\beta_k(t)
+\frac{1}{2}\int_0^T\int_{\T^N} \GG^2(x,u(x,t)) dx dt,
\end{multline}
almost-surely.
\label{prop:massm}\end{proposition}

\textbf{Proof of Proposition~\ref{prop:massm}.} We start from \eqref{eq:kineticupre}, which we apply with a test-function $\varphi$ independent on $x$. By subtracting $\<\mathbf{1}_{0>\xi},\varphi\>$ to both sides of the equation, we obtain
\begin{multline}
\<\chi(t),\varphi\>=\<\chi_0,\varphi\>
+\sum_{k\geq 1}\int_0^t\int_{\T^N}g_k(x,u(x,s))\varphi(x,u(x,s)) dx d\beta_k(s)\\
+\frac{1}{2}\int_0^t\int_{\T^N} \partial_\xi\varphi(x,u(x,s))\GG^2(x,u(x,s)) dx ds-m(\partial_\xi\varphi)([0,t]),
\label{eq:chiP}\end{multline}
where $\chi(x,t,\xi)=\mathtt{f}(x,t,\xi)-\mathbf{1}_{0>\xi}$, $\chi_0(x,\xi)=\mathtt{f}_0(x,\xi)-\mathbf{1}_{0>\xi}$ are the traditional kinetic functions used in \cite{PerthameBook} for example. We use then an approximation argument to apply \eqref{eq:chiP} with $\varphi(x,\xi)=\xi$. This gives \eqref{eq:massm}.\qed

\subsection{Generalized solutions}\label{sec:ges}

With the purpose to prepare the proof of existence of solution, we introduce the following definitions.

\begin{definition}[Young measure] Let $(X,\mathcal{A},\lambda)$ be a finite measure space. Let $\mathcal{P}_1(\R)$ denote the set of probability measures on $\R$. We say that a map $\nu\colon X\to\mathcal{P}_1(\R)$ is a Young measure on $X$ if, for all $\phi\in C_b(\R)$, the map $z\mapsto \<\nu_z,\phi\>$ from $X$ to $\R$ is measurable. We say that a Young measure $\nu$ vanishes at infinity if, for every $p\geq 1$, 
\begin{equation}
\int_X\int_\R |\xi|^p d\nu_z(\xi)d\lambda(z)<+\infty.
\label{nuvanish}\end{equation}
\label{defYoung}\end{definition}

\begin{proposition}[An alternative definition of Young measures] Let $(X,\mathcal{A},\lambda)$ be a measure space with $\lambda(X)=1$. Let $\mathcal{L}$ be the Lebesgue measure on $\R$ and let $\mathcal{Y}^1$ be the set of probability measures $\nu$ on $(X\times\R,\mathcal{A}\times\mathcal{B}(\R))$ such that $\pi_\#\nu=\lambda$, where $\pi_\#\nu$ is the push forward of $\nu$ by the projection $\pi\colon X\times\R\to X$. Then $\mathcal{Y}^1$ is the set of Young measures as defined in Definition~\ref{defYoung}.
\label{prop:altdefY}\end{proposition}

For the proof of this result, which uses the Disintegration Theorem, we refer to the discussion in \cite[p.19-20]{CastaingRaynauddeFitteValadier04} on the spaces $\mathcal{Y}^1$ and $\mathcal{Y}^1_\mathrm{dis}$ (``dis" for ``disintegration": this corresponds to the Definition~\ref{defYoung}). Note that there is no loss in generality in assuming $\lambda(X)=1$.

\begin{definition}[Kinetic function] Let $(X,\mathcal{A},\lambda)$ be a finite measure space. A measurable function $f\colon X\times\R\to[0,1]$ is said to be a kinetic function if there exists a Young measure $\nu$ on $X$ that vanishes at infinity such that, for $\lambda$-a.e. $z\in X$, for all $\xi\in\R$,
\begin{equation*}
f(z,\xi)=\nu_{z}(\xi,+\infty).
\end{equation*}
We say that $f$ is an {\rm equilibrium} if there exists a measurable function $u\colon X\to\R$ such that $f(z,\xi)=\mathtt{f}(z,\xi)=\mathbf{1}_{u(z)>\xi}$ a.e., or, equivalently, $\nu_z=\delta_{\xi=u(z)}$ for a.e. $z\in X$.
\label{def:kifunction}\end{definition}

\begin{definition}[Conjugate function] If $f\colon X\times\R\to[0,1]$ is a kinetic function, we denote by $\bar f$ the  conjugate function $\bar f:=1-f$.
\label{def:conjugate}\end{definition}

We also denote by $\chi_f$ the function defined by $\chi_f(z,\xi)=f(z,\xi)-\mathbf{1}_{0>\xi}$. This correction to $f$ is integrable on $\R$. Actually, it is decreasing faster than any power of $|\xi|$ at infinity. Indeed, we have $\chi_f(z,\xi)=-\nu_z(-\infty,\xi)$ when $\xi<0$ and $\chi_f(z,\xi)=\nu_z(\xi,+\infty)$ when $\xi>0$. Therefore
\begin{equation}
\label{e10}
|\xi|^p \int_{X} |\chi_f(z,\xi)|d\lambda(z) \le \int_{X}\int_\R |\zeta|^p d\nu_z(\zeta)
d\lambda(z) <\infty,
\end{equation}
for all $\xi\in\R$, $1\leq p<+\infty$.
\medskip

The so-called kinetic functions appear naturally when one examines the stability of a sequence of solutions to \eqref{stoSCL}. We discuss this topic in details in Section~\ref{sec:martingale}, but let us already mention the 
following compactness results.

\begin{theorem}[Compactness of Young measures] Let $(X,\mathcal{A},\lambda)$ be a finite measure 
space such that $\mathcal{A}$ is countably generated. Let $(\nu^n)$ be a sequence of Young measures on $X$ satisfying  \eqref{nuvanish} uniformly for some $p\ge 1$:
\begin{equation}
\sup_n\int_X\int_\R |\xi|^p d\nu^n_z(\xi)d\lambda(z)<+\infty.
\label{nuvanishn}\end{equation}
Then there exists a Young measure $\nu$ on $X$ and a subsequence still denoted $(\nu^n)$ such that, for all $h\in L^1(X)$, for all $\phi\in C_b(\R)$,
\begin{equation}
\lim_{n\to+\infty}\int_X h(z) \int_\R \phi(\xi) d\nu^n_z(\xi) d\lambda(z)= \int_X h(z) \int_\R \phi(\xi) d\nu_z(\xi) d\lambda(z).
\label{cvYoungMeasure}\end{equation}
\label{th:youngmeasure}\end{theorem}

The convergence \eqref{cvYoungMeasure} is the convergence for the $\tau^W_{\mathcal{Y}^1}$ topology defined in \cite[p.21]{CastaingRaynauddeFitteValadier04}. By \cite[Corollary~4.3.7]{CastaingRaynauddeFitteValadier04}, \eqref{nuvanishn} implies that the set $\{\nu_n;n\in\N\}$ is $\tau^W_{\mathcal{Y}^1}$-relatively compact, and for this result, it is not necessary to assume that $\mathcal{A}$ is countably generated. This latter hypothesis is used as a criteria of metrizability of  $\tau^W_{\mathcal{Y}^1}$, \cite[Proposition~2.3.1]{CastaingRaynauddeFitteValadier04}. A consequence of Theorem~\ref{th:youngmeasure} is the following proposition.
 
\begin{corollary}[Compactness of kinetic functions] Let $(X,\mathcal{A},\lambda)$ be a  finite  measure 
space such that $\mathcal{A}$ is countably generated. Let $(f_n)$ be a sequence of kinetic functions on $X\times\R$: $f_n(z,\xi)=\nu^n_z(\xi,+\infty)$ where $\nu^n$ are Young measures on $X$ satisfying \eqref{nuvanishn}. Then there exists a kinetic function $f$ on $X\times\R$ (related to the Young measure $\nu$ in Theorem~\ref{th:youngmeasure} by the formula $f(z,\xi)=\nu_z(\xi,+\infty)$) such that, up to a subsequence, $f_n\rightharpoonup f$ in $L^\infty(X\times\R)$ weak-*.
\label{cor:kineticfunctions}\end{corollary}

We will also need the following result.

\begin{lemma}[Convergence to an equilibrium] Let $(X,\mathcal{A},\lambda)$ be a finite measure space. Let $p> 1$. Let $(f_n)$ be a sequence of kinetic functions on $X\times\R$: $f_n(z,\xi)=\nu^n_z(\xi,+\infty)$ where $\nu^n$ are Young measures on $X$ satisfying \eqref{nuvanishn}. Let $f$ be a kinetic function on $X\times\R$ such that $f_n\rightharpoonup f$ in $L^\infty(X\times\R)$ weak-*. Assume that $f$ is an equilibrium: $f(z,\xi)=\mathtt{f}(z,\xi)=\mathbf{1}_{u(z)>\xi}$ and let
$$
u_n(z)=\int_\R\xi d\nu^n_z(\xi). 
$$
Then, for all $1\leq q<p$, $u_n\to u$ in $L^q(X)$ strong.
\label{lem:weakstrongeq}\end{lemma}

\textbf{Proof of Corollary~\ref{cor:kineticfunctions}.} We apply the theorem~\ref{th:youngmeasure}. The convergence of $(\nu^n)$, which means that
\begin{equation}\label{eq:cvtensorYoung}
\big(z\mapsto\<\nu^n_z,\phi\>\big)\to\big(z\mapsto \<\nu_z,\phi\>\big)\mbox{ in }L^\infty(X)\mbox{ weak}-*,
\end{equation}
for all $\phi\in C_b(\R)$, has the consequence that
\begin{equation}\label{eq:cvCaratheodory}
\int_X\int_\R \alpha(z,\xi)d\nu^n_z(\xi) d\lambda(z)\to\int_X\int_\R \alpha(z,\xi)d\nu_z(\xi) d\lambda(z),
\end{equation}
for every bounded Carath{\'e}odory integrand $\alpha$. This is a consequence of the identity $\tau^M_{\mathcal{Y}^1}=\tau^W_{\mathcal{Y}^1}$ in the Portmanteau Theorem \cite[Theorem~2.1.3]{CastaingRaynauddeFitteValadier04} (see also \cite[Lemma~1.2.3]{CastaingRaynauddeFitteValadier04} about Carath{\'e}odory integrands). We apply \eqref{eq:cvCaratheodory} to 
$$
\alpha(z,\xi)=\int_{-\infty}^\xi\varphi(z,\zeta)d\zeta,
$$
where $\varphi\in L^1\cap L^\infty(X\times\R)$, and apply also the Fubini theorem to obtain
\begin{equation}\label{cvweakstarfnfcor}
\iint_{X\times\R}f_n(z,\xi)\varphi(z,\xi)d\lambda(z)d\xi\to\iint_{X\times\R}f(z,\xi)\varphi(z,\xi)d\lambda(z)d\xi.
\end{equation}
Using the bound by $1$ on the $L^\infty$ norm of $f_n$ and $f$, we deduce by an argument of density that \eqref{cvweakstarfnfcor} holds true when $\varphi\in L^1(X\times\R)$. \qed\medskip

\textbf{Proof of Lemma~\ref{lem:weakstrongeq}.} Let $r\in[1,+\infty]$. By choosing $\theta=\phi$ and $\gamma$ as a test function in $z$ in \eqref{eq:cvtensorYoung}, and by use of a standard approximation procedure, we have
\begin{equation}\label{thetagamma}
\int_X\int_\R \theta(\xi)d\nu^n_z(\xi)\gamma(z) d\lambda(z)\to\int_X\theta(u(z)) \gamma(z) d\lambda(z)
\end{equation}
for all $\theta\in C(\R)$ and $\gamma\in L^r(X)$ such that 
$$
\sup_n\left\|\int_\R \theta(\xi)d\nu^n_z(\xi)\right\|_{L^{r'}(X)}<+\infty,
$$ 
where $r'$ is the conjugate exponent to $r$. Let us assume first that $p>2$ and let us prove the convergence of $(u_n)$ to $u$ in $L^2(X)$. By \eqref{thetagamma}, taking $r=2$, $\theta(\xi)=\xi$ and $\gamma\in L^2(X)$, we obtain the weak convergence of $(u_n)$ to $u$ in $L^2(X)$. By developing the scalar product
$$
\|u-u_n\|_{L^2(X)}^2=\|u\|_{L^2(X)}^2+\|u_n\|_{L^2(X)}^2-2\<u,u_n\>_{L^2(X)},
$$
we see that it is sufficient, in order to establish the strong convergence, to prove that
\begin{equation}\label{strongcvlikegrad}
\ds\limsup_{n\to+\infty}\|u_n\|_{L^2(X)}^2\leq\|u\|_{L^2(X)}^2.
\end{equation}
We obtain \eqref{strongcvlikegrad} by the Jensen inequality, which gives
\begin{equation}\label{JensenJensen}
\|u_n\|_{L^2(X)}^2=\int_X\left|\int_\R\xi d\nu^n_z(\xi)\right|^2d\lambda(z)\leq \int_X\int_\R|\xi|^2d\nu^n_x(\xi)d\lambda(z).
\end{equation}
Indeed the right-hand side of \eqref{JensenJensen} is converging to $\|u\|_{L^2(X)}^2$ (here, we apply \eqref{thetagamma} with $\theta(\xi)=\xi^2$ and $\gamma(z)=1$). Still assuming $p> 2$, the remaining cases $1\leq q<p$ are deduced from the result for $p=2$ by interpolation and by the uniform bound on $\|u_n\|_{L^p(X)}$. Let us consider the case $p\leq2$ now. Let us introduce the truncate functions and truncate sequence $(u^R_n)$ as follows:
$$
T_R(\xi):=\min(R,\max(-R,\xi)),\quad u^R_n(z)=\int_\R T_R(\xi)d\nu^n_z(\xi).
$$
One checks that the study done for $p>2$ can be applied to established the convergence $u_n^R\to T_R(u)$ in $L^r(X)$ strong for every $r<+\infty$.  Then we use the estimate
\begin{align*}
|u^R_n(z)-u_n(z)|\leq\int_{|\xi|>R}|R-\xi|d\nu^n_z(\xi)\leq 2\int_{|\xi|>R}|\xi|d\nu^n_z(\xi),
\end{align*}
from which follows, for $1\leq q<p$, by the Jensen inequality,
\begin{align*}
\|u^R_n-u_n\|_{L^q(X)}^q\leq 2^q\int_X\int_{|\xi|>R}|\xi|^qd\nu^n_z(\xi)d\lambda(z)\leq \frac{2^q}{R^{p-q}}\int_X\int_\R|\xi|^pd\nu^n_z d\lambda(z),
\end{align*}
and, thanks to \eqref{nuvanishn}, the uniform bound
$$
\|u^R_n-u_n\|_{L^q(X)}\leq \frac{2}{R^{p/q-1}}\sup_n\int_X\int_\R|\xi|^pd\nu^n_z d\lambda(z).
$$
We have also $T_R(u)\to u$ in $L^q(X)$ when $R\to+\infty$. Gathering the different results of convergence, we obtain $u_n\to u$ in $L^q(X)$. \qed\medskip

In the \textit{deterministic setting}, if $(u_n(t))$ is a sequence of solutions to \eqref{stoSCL}, then, due to natural bounds and to Theorem~\ref{th:youngmeasure}, the sequence of Young measures $\delta_{u_n(x,t)}$ on $X:=\T^N$ (consider that $t$ is fixed here) has, up to a subsequence, a limit $\nu_t$. Then every non-linear expression $\phi(u_n(t))$ for $\phi\in C_b(\R)$ will converge to $\<\nu_t,\phi\>$ in the sense of \eqref{cvYoungMeasure}. This is why it is natural (\textit{cf.} \cite{Diperna83a}), for such non-linear problems as \eqref{stoSCL}, to introduce the following generalization to Definition~\ref{defkineticsol}. 

\begin{definition}[Generalized solution]
\label{d4} Let $f_0\colon\T^N\times\R\to[0,1]$ be a kinetic function. An $L^\infty(\T^N\times\R;[0,1])$-valued process $(f(t))_{t\in[0,T]}$ is said to be a generalized solution to~\eqref{stoSCL} with initial datum $f_0$ if $f(t)$ and $\nu_t:=-\partial_\xi f(t)$ have the following properties:
\begin{enumerate}
\item\label{item:1d4} for all $t\in[0,T]$, almost-surely, $f(t)$ is a kinetic function, and, for all $R>0$, $f\in L^1_\mathcal{P}(\T^N\times(0,T)\times(-R,R)\times\Omega)$,
\item\label{item:1terd4} for all $\varphi\in C^1_c(\T^N\times\R)$, almost-surely, the map $t\mapsto\<f(t),\varphi\>$ is c{\`a}dl{\`a}g,
\item\label{item:2d4} for all $p\in[1,+\infty)$, there exists $C_p\geq 0$ such that 
\begin{equation}
\E\left(\sup_{t\in[0,T]}\int_{\T^N}\int_\R|\xi|^p d\nu_{x,t}(\xi) dx\right) \leq C_p,
\label{eq:integrabilityf}\end{equation}
\item\label{item:4d4} there exists a random measure $m$  with first moment \eqref{Firstm}, such that for all $\varphi\in C^1_c(\T^N\times\R)$, for all $t\in[0,T]$, almost-surely,
\begin{align}
\<f(t),\varphi\>=&\<f_0,\varphi\>+\int_0^t \<f(s),a(\xi)\cdot\nabla_x\varphi\>ds
\nonumber\\
&+\int_0^t\int_{\T^N}\int_\R g_k(x,\xi)\varphi(x,\xi)d\nu_{x,s}(\xi) dxd\beta_k(s) \nonumber\\
&+\frac{1}{2}\int_0^t\int_{\T^N}\int_\R \GG^2(x,\xi)\partial_\xi\varphi(x,\xi)d\nu_{x,s}(\xi) dx ds
-m(\partial_\xi\varphi)([0,t]).
\label{eq:kineticfpre}
\end{align}
\end{enumerate}
\end{definition}


Let us do a comment about notations: for each $t\in[0,T]$, we have a Young measure $\nu_t$ on $\T^N$. This gives us a set of probability measures $(\nu_{x,t})_{x\in\T^N}$, as they appear in the second line of \eqref{eq:kineticfpre}. There is something misleading in the use of the notation $\nu_{x,t}$, which conveys the idea that we are considering a Young measure $\nu$ with index space $\T^N\times(0,T)$. Such a modification of the point of view is admissible however, and we will use it fully in Section~\ref{sec:compactnun} to obtain the convergence of sequences of Young measures. Indeed, due to item~\ref{item:1d4} and to the fact that, for all $t\in(0,T)$, for a.e. $x\in\T^N$, a.s., 
$$
\int_\R \phi(\xi)d\nu_{x,t}(\xi)=\int_\R f(x,t,\xi)\phi'(\xi)d\xi
$$
if $\phi\in C^1_c(\R)$, the map $(\omega,x,t)\mapsto\<\nu_{x,t},\phi\>$ is measurable (and in $L^1_\mathcal{P}(\T^N\times(0,T)\times\Omega)$ actually). By the Fubini theorem, we deduce that, almost-surely, $(x,t)\mapsto\<\nu_{x,t},\phi\>$ is measurable when $\phi\in C^1_c(\R)$. By an argument of density, this holds true when $\phi\in C_b(\R)$.\medskip

This point about the status of $\nu_{x,t}$ being clear, we have now also to justify that the stochastic integral in \eqref{eq:kineticfpre} is well-defined: the bound \eqref{eq:integrabilityf} implies
\begin{equation}\label{eq:integrabilityftime}
\E\left(\int_0^T\int_{\T^N}\int_\R|\xi|^p d\nu_{x,t}(\xi) dx\right) \leq C_p T.
\end{equation} 
Using successively Jensen's Inequality, the growth hypothesis \eqref{D0}, and \eqref{eq:integrabilityftime} with $p=2$, we obtain, for $\varphi\in C^1_c(\T^N\times\R)$,
\begin{align}
\E\int_0^T\sum_{k\geq 1}\left|\int_{\T^N}\int_\R g_k(x,\xi)\right. &\left.\varphi(x,\xi)d\nu_{x,t}(\xi)dx\right|^2\nonumber\\
&\leq \E\int_0^T\sum_{k\geq 1}\int_{\T^N}\int_\R \left| g_k(x,\xi)\varphi(x,\xi)\right|^2 d\nu_{x,t}(\xi)dx\nonumber\\
&=\E\int_0^T\int_{\T^N}\int_\R \GG^2(x,\xi)\left| \varphi(x,\xi)\right|^2 d\nu_{x,t}(\xi)dx\nonumber\\
&\leq\|\varphi\|_{L^\infty_{x,\xi}}^2D_0 (1+C_2T).\label{sumhk}
\end{align}
The fact that 
$$
t\mapsto\int_{\T^N}\int_\R g_k(x,\xi)\varphi(x,\xi)d\nu_{x,t}(\xi)dx
$$
is predictable is a consequence of item~\ref{item:1d4}. To sum up, we have proved the following result.

\begin{lemma}[Admissible integrand] Let $f_0\colon\T^N\times\R\to[0,1]$ be a kinetic function. Let $(f(t))_{t\in[0,T]}$ be a generalized solution to~\eqref{stoSCL} with initial datum $f_0$. Then, for all $\varphi\in C^1_c(\T^N\times\R)$ the $l^2(\N^*)$-valued process
$$
t\mapsto \left(\int_{\T^N}\int_\R g_k(x,\xi)\varphi(x,\xi)d\nu_{x,t}(\xi)dx\right)_{k\geq 1}
$$ 
is in $L^2_\mathcal{P}([0,T]\times\Omega;l^2(\N^*))$.
\label{lem:admint}\end{lemma}
Let us now state a simple result of reduction from generalized solution to mere solution.

\begin{proposition} Let $u_{0}\in L^\infty(\T^N)$. Let $(f(t))_{t\in[0,T]}$ be a generalized solution to~\eqref{stoSCL} with initial datum $\mathbf{1}_{u_0>\xi}$. If for all $t\in[0,T]$, $f(t)$ is an \emph{equilibrium}:
\begin{equation}\label{hyp:fequilibrium}
f(x,t,\xi,\omega)=\mathtt{f}(x,t,\xi,\omega)=\mathbf{1}_{u(x,t,\omega)>\xi},
\end{equation}
for a.e. $(x,\xi,\omega)\in\T^N\times\R\times\Omega$, then $(u(t))_{t\in[0,T]}$ is a solution to \eqref{stoSCL} with initial datum $u_0$.
\label{prop:gestoes}\end{proposition}

\textbf{Proof of Proposition~\ref{prop:gestoes}.} Under \eqref{hyp:fequilibrium}, we have $\nu_t=\delta_{u(t)}$ a.s. From \eqref{eq:integrabilityftime} with $p=1$, we deduce that $u\in L^1(\T^N\times(0,T)\times\Omega)$. Since
$$
u(x,t)=\int_\R\xi d\nu_{x,t}(\xi),
$$
we obtain $u\in L^1_{\mathcal{P}}(\T^N\times(0,T)\times\Omega)$ as a consequence of Item~\ref{item:1d4} in Definition~\ref{d4}. We have also 
$$
\<\mathtt{f}(t),\varphi\>=\int_{\T^N}\int_\R\psi(x,\xi) d\nu_{x,t}(\xi)dx,\quad\psi(x,\xi):=\int_{-\infty}^\xi\varphi(x,\zeta)d\zeta,
$$
for all $\varphi\in C^\infty_c(\T^N\times\R)$. Therefore Item~\ref{item:1bisdefsol} in Definition~\ref{defkineticsol} follows from Item~\ref{item:1terd4} in Definition~\ref{d4}. Using the identity
$$
\int_{\T^N} |u(t,x)|^p dx = \int_{\T^N}\int_\R |\xi|^p d\nu_{x,t}(\xi)dx,
$$
we obtain Item~\ref{item:2defsol} in Definition~\ref{defkineticsol}. Item~\ref{item:3defsol} in Definition~\ref{defkineticsol} follows from Item~\ref{item:4d4} in Definition~\ref{d4}. \qed
\medskip

We will show in Theorem~\ref{th:Uadd} that \eqref{hyp:fequilibrium}, which we give as an hypothesis in Proposition~\ref{prop:gestoes}, is automatically satisfied by any generalized solution starting from an equilibrium $f_0=\mathtt{f}_0=\mathbf{1}_{u_0>\xi}$.\medskip

We conclude this paragraph with the following result, used in the proof of Corollary~\ref{cor:timecontinuity}.
\begin{lemma}[Distance to equilibrium] Let $(X,\lambda)$ be a finite measure space. Let $f\colon X\times\R\to[0,1]$ be a kinetic function. Then
$$
m(\xi):=\int_{-\infty}^\xi (\mathbf{1}_{u>\zeta}-f(\zeta))d\zeta,\quad\mbox{where } u:=\int_\R\chi_{f}(\zeta)d\zeta,
$$
is well defined and non-negative.
\label{lemdisteq}\end{lemma}

Note in particular that the difference $f(\xi)-\mathbf{1}_{u>\xi}$ writes $\partial_\xi m$ where $m\geq 0$.
\medskip

\textbf{Proof of Lemma~\ref{lemdisteq}.} Let $\nu_z=-\partial_\xi f(z,\cdot)$, $z\in X$. By Jensen's Inequality, we have
\begin{equation}\label{JensenH}
H\left(\int_\R \zeta d\nu_z(\zeta)\right)\leq \int_\R H(\zeta)d\nu_z(\zeta)
\end{equation}
for all convex sub-linear function $H\colon\R\to\R$. Note that
$$
u(z)=\int_\R f(z,\zeta)-\mathbf{1}_{0>\zeta} d\zeta=\int_\R \zeta d\nu_z(\zeta)
$$
by integration by parts. By integration by parts, we also have, for all sub-linear function $H\in C^1(\R)$,
$$
\int_\R H(\zeta)d\nu_z(\zeta)=H(0)+\int_\R H'(\zeta)(f(z,\zeta)-\mathbf{1}_{0>\zeta} )d\zeta
$$
and
$$
H(u(z))=\int_\R H(\zeta)d\delta_{u(z)}(\zeta)=H(0)+\int_\R H'(\zeta)(\mathbf{1}_{u(z)>\zeta}-\mathbf{1}_{0>\zeta} )d\zeta.
$$
By \eqref{JensenH}, it follows that
$$
\int_\R H'(\zeta)(f(z,\zeta)-\mathbf{1}_{u(z)>\zeta} )d\zeta\geq 0
$$
for all convex and sub-linear $H\in C^1(\R)$. Approximating $\zeta\mapsto (\zeta-\xi)^-$ by such functions $H$, we obtain $m(\xi)\geq 0$. \qed

\subsection{Left limits of generalized solutions}

If $(f(t))_{t\in[0,T]}$ is a generalized solution to~\eqref{stoSCL} and $\varphi\in C^\infty_c(\T^N\times\R)$, then, a.s., $t\mapsto\<f(t),\varphi\>$ is c{\`a}dl{\`a}g. In the next proposition, we show that the a.s.-property to be c{\`a}dl{\`a}g is independent on $\varphi$ and that the limit from the left at any point $t_*\in(0,T]$ is represented by a kinetic function.

\begin{proposition} Let $f_0$ be a kinetic initial datum. Let $(f(t))_{t\in[0,T]}$ be a generalized solution to~\eqref{stoSCL} with initial datum $f_0$. Then 
\begin{enumerate}
\item there exists a measurable subset $\hat\Omega\subset\Omega$ of probability $1$ such that, for all $\omega\in\hat{\Omega}$, for all $\varphi\in C_c(\T^N\times\R)$, $t\mapsto\<f(\omega,t),\varphi\>$ is c{\`a}dl{\`a}g,
\item there exists an $L^\infty(\T^N\times\R;[0,1])$-valued process $(f^-(t))_{t\in(0,T]}$ such that: for all $t\in(0,T]$, for all $\omega\in\hat{\Omega}$, for all $\varphi\in C_c(\T^N\times\R)$, $f^-(t)$ is a kinetic function on $\T^N$ which represents the left limit of $s\mapsto\<f(s),\varphi\>$ at $t$:
\begin{equation}\label{fmoinsmoins}
\<f^-(t),\varphi\>=\lim_{s\to t-}\<f(s),\varphi\>.
\end{equation}
\end{enumerate} 
\label{prop:LRlimits}\end{proposition}

{\bf Proof of Proposition~\ref{prop:LRlimits}.} The set of test functions $C^1_c(\T^N\times\R)$ (endowed with the topology of the uniform 
convergence of the functions and their first derivatives) is separable and we fix 
a dense countable subset $\mathcal D_1$ (see the argument about $\Gamma$ in Section~\ref{sec:proofThmartingalestatespace} for a proof of the existence of $\mathcal D_1$). For all $\varphi\in C^1_c(\T^N\times\R)$, a.s., the 
map
\begin{multline}
J_\varphi\colon t\mapsto\int_0^t \<f(s),a(\xi)\cdot\nabla\varphi\>ds+\sum_{k\geq 1}\int_0^t\int_{\T^N}\int_\R g_k(x,\xi)\varphi(x,\xi)d\nu_{x,s}(\xi)dxd\beta_k(s)\\
+\frac{1}{2}\int_0^t\int_{\T^N}\int_\R \partial_\xi\varphi(x,\xi)\GG^2(x,\xi)d\nu_{x,s}(\xi) dx ds
\label{def:Jphi}\end{multline}
is continuous on $[0,T]$. Consequently: a.s., say for $\omega\in\Omega_1$ where $\Omega_1$ is of full measure, for all $\varphi\in \mathcal D_1$, $J_\varphi$ is continuous on $[0,T]$. If $\varphi\in \mathcal D_1$, \eqref{eq:kineticfpre} gives $\<f(t),\varphi\>$ as a sum (up to the constant $\<f_0,\varphi\>$) of $J_\varphi(t)$ with $m(\partial_\xi\varphi)([0,t])$. This latter expression defines a function c{\`a}dl{\`a}g in $t$ for all $\omega\in\Omega_2$, hence $t\mapsto \<f(t),\varphi\>$ is c{\`a}dl{\`a}g if $\omega\in\Omega_1\cap\Omega_2$. Here, $\Omega_2\subset\Omega$ is of full measure. Next, we use the estimate~\eqref{eq:integrabilityf}: there exists a set of full measure $\Omega_3$ in $\Omega$ such that, for every $\omega\in\Omega_3$,
\begin{equation}
\sup_{t\in[0,T]}\int_{\T^N}\int_\R|\xi|^p d\nu_{x,t}(\xi) dx \leq C_p(\omega)<+\infty.
\label{eq:integrabilityfomega}\end{equation}
Let $\omega\in\hat\Omega:=\Omega_1\cap\Omega_2\cap\Omega_3$ be fixed. If $\varphi\in C_c(\T^N\times\R)$, then
\begin{equation}\label{fVSnu}
\<f(t),\varphi\>=\int_{\T^N\times\R}\psi(x,\xi) d\nu_{x,t}(\xi)dx,\quad\psi(x,\xi):=\int_{-\infty}^\xi\varphi(x,\zeta)d\zeta.
\end{equation}
Let $R_\varphi>0$ be such that $\varphi$ is supported in $[-R_\varphi,R_\varphi]$. Since $|\psi(x,\xi)|\leq \|\varphi\|_{L^\infty}(R_\varphi+|\xi|)$, we obtain the bound $\sup_{t\in[0,T]}|\<f(t),\varphi\>|\leq \|\varphi\|_{L^\infty}(R_\varphi+C_1(\omega))$. This gives the continuity of $\<f(t),\varphi\>$ with respect to $\varphi$. Since the space of c{\`a}dl{\`a}g functions is closed under uniform convergence, an argument of density shows that $t\mapsto \<f(t),\varphi\>$ is c{\`a}dl{\`a}g for all $\varphi\in C_c(\T^N\times\R)$. To prove the second assertion of the proposition, let us fix $\omega\in\hat{\Omega}$ and consider an increasing sequence $(t_n)$ in $[0,T]$ converging to a point $t_*\in(0,T]$. Then, by means of \eqref{eq:integrabilityfomega} and since the Borel $\sigma$-algebra of $\T^N$ is countably generated ($\T^N$ being separable),  we can apply  Corollary~\ref{cor:kineticfunctions}: there exist a kinetic function $f^{*,-}$ on $\T^N\times\R$ and a subsequence $(n_k)$ such that $f(t_{n_k})\rightharpoonup f^{*,-}$
weakly-$*$ in $L^\infty(\T^N\times\R)$ as $k\to+\infty$. If an other subsequence $(\tilde{n}_k)$ provides an other weak limit $\tilde{f}^{*,-}$, then we have
$$
\<f^{*,-},\varphi\>=\lim_{t\to t^*-}\<f(t),\varphi\>=\<\tilde{f}^{*,-},\varphi\>
$$
for all $\varphi\in C_c(\T^N\times\R)$. Therefore $f^{*,-}=\tilde{f}^{*,-}$: there is only one possible limit. It follows that the whole sequence $(f(t_n))$ is converging to $f^{*,-}$ in $L^\infty(\T^N\times\R)$ weak-$*$. We establish this fact to ensure that the subsequence $(n_k)$ is independent on $\omega$. Indeed, this shows that, viewed as a function of $(\omega,x,\xi)$, $f^{*,-}$ is measurable. We set $f^-(t_*)=f^{*,-}$ to conclude.\qed\bigskip

\begin{remark}[Left and right limits] Note that we prove a little bit more than what is stated in Proposition~\ref{prop:LRlimits}. Indeed, for $\omega\in\hat{\Omega}$, we have $f(s)\to f^-(t)$ in $L^\infty(\T^N\times\R)$ for the weak-$*$ topology, when $s\uparrow t$, which implies \eqref{fmoinsmoins}. By similar arguments, we can show that $f(s)\to f(t)$ in $L^\infty(\T^N\times\R)$ weak-$*$ when $s\downarrow t$.
\label{rk:LRlimits}\end{remark}

\begin{remark}[Uniform bound] Note that, by construction, $\nu^-=-\partial_\xi f^-$ satisfies the following bounds: for all $\omega\in\hat\Omega$, 
\begin{equation}
\sup_{t\in[0,T]}\int_{\T^N}\int_\R|\xi|^p d\nu^-_{x,t}(\xi) dx\leq C_p(\omega),\quad \E\left(\sup_{t\in[0,T]}\int_{\T^N}\int_\R|\xi|^p d\nu^-_{x,t}(\xi) dx\right)\leq C_p.
\label{eq:integrabilityfpm}\end{equation}
We obtain \eqref{eq:integrabilityfpm} using \eqref{eq:integrabilityf}-\eqref{eq:integrabilityfomega} and Fatou's Lemma. 
\label{rk:integrabilityfpm}\end{remark}

\begin{remark}[Equation for $f^-$] Passing to the limit in \eqref{eq:kineticfpre} for an increasing sequence of times $t$, we obtain the following equation on $f^-$:
\begin{align}
\<f^-(t),\varphi\>=&\<f(0),\varphi\>+\int_0^t \<f(s),a(\xi)\cdot\nabla_x\varphi\>ds
\nonumber\\
&+\int_0^t\int_{\T^N}\int_\R g_k(x,\xi)\varphi(x,\xi)d\nu_{x,s}(\xi) dxd\beta_k(s) \nonumber\\
&+\frac{1}{2}\int_0^t\int_{\T^N}\int_\R \GG^2(x,\xi)\partial_\xi\varphi(x,\xi)d\nu_{x,s}(\xi) dx ds
-m(\partial_\xi\varphi)([0,t)).
\label{eq:kineticfmoinspre0}
\end{align}
In particular, we have
\begin{equation}
\<f(t)-f^-(t),\varphi\>=-m(\partial_\xi\varphi)(\{t\}).
\label{f+VSf-}\end{equation}
Outside the set of atomic points of $A\mapsto m(\partial_\xi\varphi)(A)$, which is at most countable, we have $\<f(t),\varphi\>=\<f^-(t),\varphi\>$. It follows that $f=f^-$ a.e. 
In particular, \eqref{eq:kineticfmoinspre0} gives us the following equation on $f^-$:
\begin{align}
\<f^-(t),\varphi\>=&\<f(0),\varphi\>+\int_0^t \<f^-(s),a(\xi)\cdot\nabla_x\varphi\>ds
\nonumber\\
&+\int_0^t\int_{\T^N}\int_\R g_k(x,\xi)\varphi(x,\xi)d\nu^-_{x,s}(\xi) dxd\beta_k(s) \nonumber\\
&+\frac{1}{2}\int_0^t\int_{\T^N}\int_\R \GG^2(x,\xi)\partial_\xi\varphi(x,\xi)d\nu^-_{x,s}(\xi) dx ds
-m(\partial_\xi\varphi)([0,t)),
\label{eq:kineticfmoinspre}
\end{align}
equation which is also valid for $t=0$ if we set $f^-(0)=f_0$.
\label{rk:fmoinseq}
\end{remark}

In the next proposition, we give a criterion for the continuity of $t\mapsto\<f(t),\varphi\>$ at a given point.

\begin{proposition}[The case of equilibrium] Let $f_0$ be a kinetic initial datum. Let $(f(t))_{t\in[0,T]}$ be a generalized solution to~\eqref{stoSCL} with initial datum $f_0$. Let $t\in(0,T]$. Assume that $f^-(t)$ is at equilibrium: there exists a random variable $v\in L^1(\T^N)$ such that $f^{-}(t,\xi)=\mathbf{1}_{v>\xi}$ a.s. Then $f^-(t)=f(t)$.
\label{prop:equilibrium}\end{proposition}

\textbf{Proof of Proposition~\ref{prop:equilibrium}.} Let $m^*$ denote the restriction of $m$ to $\T^N\times\{t\}\times\R$. Let us also set $f^+=f(t)$. By \eqref{f+VSf-}, we thus have
\begin{equation}\label{kimstarstar}
f^+-\mathbf{1}_{v>\xi}=\partial_\xi m^*.
\end{equation}
There exists a subset $\Omega_4$ of $\Omega$ of probability $1$ such that, for all $\omega\in\Omega_4$, $m$, and thus $m^*$, are finite measures on $\T^N\times[0,T]\times\R$ and $\T^N\times\R$ respectively. Let $\psi$ be a smooth non-negative function such that $0\leq\psi\leq 1$, $\psi\equiv 1$ on $[-1,1]$, $\psi$ being supported in $[-2,2]$. Define the cut-off function $\psi_\eps(\xi)=\psi(\eps\xi)$. Let also $\varphi\in C(\T^N)$. By \eqref{kimstarstar}, we have
\begin{align*}
\iint_{\T^N\times\R}(f^+(x,\xi)-\mathbf{1}_{v(x)>\xi})\varphi(x)\psi_\eps(\xi)dx d\xi&=-\eps \iint_{\T^N\times\R}\varphi(x)\psi'(\eps\xi)d m^*\\
&\leq\eps\|\varphi\|_{L^\infty(\T^N)}\|\psi'\|_{L^\infty(\R)}m^*(\T^N\times\R).
\end{align*}
Taking the limit $\eps\to0$, and taking in consideration the fact that $\varphi$ is arbitrary, we deduce that, for all $\omega\in\hat\Omega\cap\Omega_4$, for a.e. $x\in\T^N$, 
$$
\int_\R( f^+(x,\xi)-\mathbf{1}_{0>\xi}) d\xi=\int_\R(\mathbf{1}_{v(x)>\xi}-\mathbf{1}_{0>\xi}) d\xi=v(x).
$$
Introduce now
$$
p^*\colon\xi\mapsto\int_{-\infty}^\xi (\mathbf{1}_{v>\zeta}-f^+(\zeta))d\zeta.
$$
By Lemma~\ref{lemdisteq}, $p^*$ is non-negative. In addition, $\partial_\xi (m^*+p^*)=0$ due to \eqref{kimstarstar} and the definition of $p^*$. Therefore $m^*+p^*$ is constant, and actually vanishes by the condition at infinity \eqref{inftym} and the obvious fact that $p^*(\T^N\times B_R^c)$ also vanishes when $R\to+\infty$. Since $m^*,p^*\geq 0$, we finally obtain $m^*=0$ and conclude to the identity $f^-(t)=f(t)$. \qed\medskip

Let us consider also the special case $t=0$. By letting $t\to 0+$ in \eqref{eq:kineticfpre}, we have $f(0)-f_0=\partial_\xi m_0$, where $m_0$ is the restriction of $m$ to $\T^N\times\{0\}\times\R$. Consequently, we have the following corollary to Proposition~\ref{prop:equilibrium}.

\begin{corollary} Let $f_0$ be a kinetic initial datum. Let $(f(t))_{t\in[0,T]}$ be a generalized solution to~\eqref{stoSCL} with initial datum $f_0$. Assume that $f_0$ is at equilibrium. Then $f(0)=f_0$ and $m$ does not charge the line $\{t=0\}$: $m(\T^N\times\{0\}\times\R)=0$ a.s.
\label{cor:equilibrium}\end{corollary}

Our final result in this section is about trajectories of solutions to \eqref{stoSCL}. It is an intermediate statement, before the full continuity result given in Corollary~\ref{cor:timecontinuity}.

\begin{proposition} Let $u_{0}\in L^\infty(\T^N)$. Let $(u(t))_{t\in[0,T]}$ be a solution to~\eqref{stoSCL} with initial datum $u_0$. Then, for all $p\in[1,+\infty)$, for all $\omega\in\hat{\Omega}$ (given in Proposition~\ref{prop:LRlimits}), the map $t\mapsto u(t)$ from $[0,T]$ to $L^p(\T^N)$ is continuous from the right.
\label{prop:LRlimitsu}\end{proposition}

{\bf Proof of Proposition~\ref{prop:LRlimitsu}.} We apply Proposition~\ref{prop:LRlimits} to $f(t)=\mathtt{f}(t)=\mathbf{1}_{u(t)>\xi}$. For $\omega\in\hat{\Omega}$, $\varphi\in C_c(\T^N\times\R)$, the map 
$
t\mapsto \<\mathtt{f}(t),\varphi\>
$
is c{\`a}dl{\`a}g. Let $t_*\in [0,T)$ and let $(t_n)$ be a decreasing sequence of $[0,T]$ converging to $t_*$. The sequence $f^n$ of elements $f^n:=\mathtt{f}(t_n)$ takes values in $[0,1]$. For $\omega\in\hat{\Omega}$ fixed, it has a convergent subsequence in $L^\infty(\T^N\times\R)$ weak-*. Since $\<f^n,\varphi\>\to\<\mathtt{f}(t_*),\varphi\>$ for all continuous, compactly supported function $\varphi$ on $\T^N\times\R$, the whole sequence $(f^n)$ is converging to its unique adherence value, $\mathtt{f}(t_*)$. By \eqref{eq:integrabilityu}, the bound \eqref{nuvanishn} is satisfied for all $p\in[1,+\infty)$: we can apply Lemma~\ref{lem:weakstrongeq} to conclude to the convergence $u(t_n)\to u(t_*)$ in $L^p(\T^N)$. \qed


\section{Comparison, uniqueness and reduction of generalized solutions}\label{sec:comparison}

\subsection{Doubling of variables}\label{sec:doubling}
In this paragraph, we prove a technical proposition relating two generalized solutions $f_i$, $i=1,2$ of the equation
\begin{equation}\label{stoSCLi}
du_i(x,t)+\div(A(u_i(x,t)))dt=\Phi_i(x,u_i(x,t))dW(t).
\end{equation}

We use the following convention of notations: if $(f(t))_{t\in[0,T]}$ is a generalized solution to \eqref{stoSCL}, we denote by $f^-$ the left limit defined in Proposition~\ref{prop:LRlimits}, and we denote by $f^+$ the right limit, which is simply $f$: $f^+(t):=f(t)$. This gives more homogeneity to the different statements in this part. Recall also the notation $\bar{f}=1-f$ for the conjugate to $f$, introduced in Definition~\ref{def:conjugate}.

\begin{proposition} Let $f_i$, $i=1,2$, be generalized solution to \eqref{stoSCLi}. Then, for $0\leq t\leq T$, and non-negative test functions $\rho\in C^\infty(\T^N)$, $\psi\in C^\infty_c(\R)$, we have 
\begin{multline}
\E\int_{(\T^N)^2}\int_{\R^2} \rho(x-y)\psi(\xi-\zeta)f_1^\pm(x,t,\xi)\bar f_2^\pm(y,t,\zeta) d\xi d\zeta dx dy \\
\leq\int_{(\T^N)^2}\int_{\R^2} \rho(x-y)\psi(\xi-\zeta) f_{1,0}(x,\xi)\bar f_{2,0}(y,\zeta)d\xi d\zeta dx dy +\mathrm{I}_\rho+\mathrm{I}_\psi,
\label{CR0}\end{multline}
where
$$
\mathrm{I}_\rho=\E\int_0^t\int_{(\T^N)^2}\int_{\R^2} f_1(x,s,\xi)\bar f_2(y,s,\zeta)(a(\xi)-a(\zeta))\psi(\xi-\zeta)d\xi d\zeta
 \cdot\nabla_x\rho(x-y) dx dy ds
$$
and
$$
\mathrm{I}_\psi=\frac{1}{2}\int_{(\T^N)^2}\rho(x-y)\E\int_0^t\int_{\R^2} \psi(\xi-\zeta)
\sum_{k\geq 1}|g_{k,1}(x,\xi)-g_{k,2}(y,\zeta)|^2 d\nu^1_{x,s}\otimes\nu^2_{y,s}(\xi,\zeta) dx dy ds.
$$
\label{prop:CR}\end{proposition}

\begin{remark}\label{intfbarf} Each term in \eqref{CR0} is finite. Let us for instance consider the left-hand side of \eqref{CR0}. Introduce the auxiliary functions 
\begin{equation*}
\psi_1(\xi)=\int_{-\infty}^\xi\psi(s) ds,\quad \psi_2(\zeta)=\int_{-\infty}^\zeta\psi_1(\xi)d\xi.
\end{equation*}
Since $\psi$ is compactly supported, both $\psi_1$ and $\psi_2$ vanish at $-\infty$. When $\xi\to+\infty$, $\psi_1$ remains bounded while $\psi_2$ has linear 
growth. More precisely, if $\psi$ is supported in $[-R,R]$, then 
\begin{equation}\label{psi2sublin}
|\psi_2(\zeta)|\leq (|\zeta|+R)\|\psi\|_{L^1(\R)}.
\end{equation} 
Since
$$
f_1^\pm(x,t,\xi)=\int_{(\xi,+\infty)}d\nu^{1,\pm}_{x,t}(\xi),\quad\bar f_2^\pm(y,t,\zeta)=\int_{(-\infty,\zeta)}d\nu^{2,\pm}_{y,t}(\zeta),
$$
for a.e. $\xi$, $\zeta\in\R$, $x,y\in\T^N$, $t\in[0,T]$, the Fubini Theorem gives us the formula
\begin{equation}\label{psifbarf}
\int_{\R^2}\psi(\xi-\zeta)f_1^\pm(x,t,\xi)\bar f_2^\pm(y,t,\zeta)d\xi d\zeta=\int_{\R^2}\psi_2(u-v) d\nu^{1,\pm}_{x,t}(u) d\nu^{2,\pm}_{y,t}(v).
\end{equation}
By \eqref{psi2sublin}, we deduce that
\begin{multline*}
\left|\int_{\R^2}\psi(\xi-\zeta)f_1^\pm(x,t,\xi)\bar f_2^\pm(y,t,\zeta)d\xi d\zeta\right|\\
\leq \|\psi\|_{L^1(\R)}\left[R+\int_\R |\xi|d\nu^{1,\pm}_{x,t}(\xi)+\int_\R |\xi|d\nu^{2,\pm}_{y,t}(\xi)\right],
\end{multline*}
for a.e. $x,y\in \T^N$, for all $t\in[0,T]$. Using the Young inequality for convolution with indices $1,1,1$, we obtain
\begin{multline}\label{boundCRleft}
\left|\int_{(\T^N)^2}\int_{\R^2} \rho(x-y)\psi(\xi-\zeta)f_1^\pm(x,t,\xi)\bar f_2^\pm(y,t,\zeta)d\xi d\zeta  dx dy\right|\\
\leq \|\psi\|_{L^1(\R)}\|\rho\|_{L^1(\T^N)}(R+C_{1,1}(\omega)+C_{1,2}(\omega)),
\end{multline}
where
$$
C_{1,i}(\omega):=\sup_{t\in[0,T]}\int_{\T^N}\int_\R|\xi| d\nu^{i,\pm}_{x,t}(\xi) dx
$$
is in $L^1(\Omega)$ thanks to \eqref{eq:integrabilityf}-\eqref{eq:integrabilityfpm}.
\end{remark}

{\bf Proof of Proposition~\ref{prop:CR}.} Set 
$$
G_i^2(x,\xi)=\sum_{k=1}^\infty|g_{k,i}(x,\xi)|^2,\quad i\in\{1,2\}.
$$ 
Let $\varphi_1\in C^\infty_c(\T^N_x\times\R_\xi)$ and $\varphi_2\in C^\infty_c(\T^N_y\times\R_\zeta)$ be some given test-functions. Equation~\eqref{eq:kineticfpre} for $f_1=f_1^+$ reads $ \<f_1^{+}(t),\varphi_1\>=\mu_1([0,t])+F_1(t)$,
where $F_1$ is the stochastic integral
$$
F_1(t)=\sum_{k\geq 1}\int_0^t\int_{\T^N}\int_\R g_{k,1}\varphi_1 d\nu^1_{x,s}(\xi)dxd\beta_k(s)
$$
and $t\mapsto \mu_1([0,t])$ is the function of finite variation on $[0,T]$ (\textit{cf.} \cite[p.~5]{RevuzYor99}) defined by
\begin{multline*}
\mu_1([0,t])= \<f_{1,0},\varphi_1\>\delta_0([0,t])+
\int_0^t \<f_1,a\cdot\nabla\varphi_1\>ds\\
+\frac{1}{2}\int_0^t\int_{\T^N}\int_\R \partial_\xi\varphi_1\GG^2_1 d\nu^1_{x,s}(\xi) dx ds-m_1(\partial_\xi\varphi_1)([0,t]).
\end{multline*}
Note that, by Corollary~\ref{cor:equilibrium}, $m_1(\partial_\xi\varphi_1)(\{0\})=0$ and thus the value of $\mu_1(\{0\})$ is $ \<f_{1,0},\varphi_1\>$. Similarly, we write $\<\bar f_2^{+}(t),\varphi_2\>$ as continuous semi-martingale, sum of the stochastic integral
$$
\bar F_2(t)=-\sum_{k\geq 1}\int_0^t\int_{\T^N}\int_\R g_{k,2}\varphi_2 d\nu^2_{y,s}(\zeta)dyd\beta_k(s)
$$
with the function with finite variation given by
\begin{multline*}
\mu_2([0,t])=\<\bar f_{2,0},\varphi_2\>\delta_0([0,t]) +
\int_0^t \<\bar f_2,a\cdot\nabla\varphi_2\>ds\\
-\frac{1}{2}\int_0^t\int_{\T^N}\int_\R \partial_\xi\varphi_2\GG^2_2 d\nu^2_{y,s}(\zeta) dy ds
+m_2(\partial_\zeta\varphi_2)([0,t]).
\end{multline*}
Again, we note that $\mu_2(\{0\})=\<\bar f_{2,0},\varphi_2\>$. Let us define the test-function 
$$
\alpha(x,\xi,y,\zeta)=\varphi_1(x,\xi)\varphi_2(y,\zeta).
$$
We want to compute 
\begin{equation}\label{prodf1f2}
\<\<f_1^+(t)\bar f_2^+(t),\alpha\>\>=\<f_1^+(t),\varphi_1\>\<\bar f_2^+(t),\varphi_2\>,
\end{equation}
where $\<\<\cdot,\cdot\>\>$ denotes the duality product over $\T^N_x\times\R_\xi\times\T^N_y\times\R_\zeta$.
By the It{\={o}} formula for continuous semimartingales, \cite[p.~146]{RevuzYor99}, taking expectation, we obtain
the following identity:
\begin{multline}
\E \<\<f_1^+(t)\bar f_2^+(t),\alpha\>\>
=\<\<f_{1,0}\bar f_{2,0},\alpha\>\>\\
+\E\int_0^t\int_{(\T^N)^2}\int_{\R^2} f_1 \bar f_2 (a(\xi)\cdot\nabla_x +a(\zeta)\cdot\nabla_y)\alpha d\xi d\zeta dx dy ds\\
+\frac12 \E\int_0^t\int_{(\T^N)^2}\int_{\R^2} \partial_\xi\alpha \bar f_2(s)\GG^2_1 d\nu^1_{x,s}(\xi) d\zeta dx dy ds\\
-\frac12 \E\int_0^t\int_{(\T^N)^2}\int_{\R^2} \partial_\zeta\alpha f_1(s) \GG^2_2 d\nu^2_{y,s}(\zeta) d\xi dy dx ds\\
-\E\int_0^t\int_{(\T^N)^2}\int_{\R^2}\GG_{1,2} \alpha d\nu^1_{x,s}(\xi) d\nu^2_{y,s}(\zeta) dx dy ds\\
-\E\int_{(0,t]}\int_{(\T^N)^2}\int_{\R^2} \bar f_2^+(s)\partial_\xi\alpha dm_1(x,s,\xi) d\zeta dy\\
+\E\int_{(0,t]}\int_{(\T^N)^2}\int_{\R^2} f_1^-(s)\partial_\zeta\alpha dm_2(y,s,\zeta) d\xi dx
\label{CR1}\end{multline}
where $\GG_{1,2}(x,y;\xi,\zeta):=\sum_{k\geq 1}g_{k,1}(x,\xi) g_{k,2}(y,\zeta)$. By a 
density argument, \eqref{CR1} remains true for any test-function $\alpha\in C^\infty_c(\T^N_x\times\R_\xi\times\T^N_y\times\R_\zeta)$. Using similar arguments as in Remark \ref{intfbarf}, the assumption that $\alpha$ is compactly supported can be relaxed using the to the condition at infinity \eqref{inftym} on $m_i$ and \eqref{nuvanish} on $\nu^i$, $i=1,2$. Using truncates of $\alpha$, we obtain that \eqref{CR1} remains true if $\alpha\in C^\infty_b(\T^N_x\times\R_\xi\times\T^N_y\times\R_\zeta)$ is supported in a neighbourhood of the diagonal
\begin{equation*}
\{(x,\xi,x,\xi);x\in\T^N,\xi\in\R\}.
\end{equation*} 
We then take $\alpha=\rho\psi$ where $\rho=\rho(x-y)$, $\psi=\psi(\xi-\zeta)$. Note the remarkable identities
\begin{equation}
(\nabla_x+\nabla_y)\alpha=0,\quad (\partial_\xi+\partial_\zeta)\alpha=0.
\label{ddvar}\end{equation}
In particular, the last term in \eqref{CR1} is
\begin{align*}
\E\int_{(0,t]}\int_{(\T^N)^2}\int_{\R^2} f_1^-(s)\partial_\zeta\alpha d\xi dx & dm_2(y,s,\zeta)\\
=&-\E\int_{(0,t]}\int_{(\T^N)^2}\int_{\R^2} f_1^-(s)\partial_\xi\alpha d\xi dx dm_2(y,s,\zeta)\\
=&- \E\int_{(0,t]}\int_{(\T^N)^2}\int_{\R^2}\alpha  d\nu^{1,-}_{x,s}(\xi) dx dm_2(y,s,\zeta)\leq 0 
\end{align*}
since $\alpha\geq 0$. The symmetric term 
\begin{align*}
-\E\int_{(0,t]}\int_{(\T^N)^2}\int_{\R^2} \bar f_2^+(s)\partial_\xi\alpha dm_1(x,s,\xi)& d\zeta dy\\
=&-\E\int_{(0,t]}\int_{(\T^N)^2}\int_{\R^2} \alpha d\nu^{2,+}_{y,s}(\zeta) dy dm_1(x,s,\xi)
\end{align*}
is, similarly, non-positive. Consequently, we have
\begin{equation}
\E \<\<f_1^+(t)\bar f_2^+(t),\alpha\>\>\leq \<\<f_{1,0}\bar f_{2,0},\alpha\>\>+\mathrm{I}_\rho+\mathrm{I}_\psi,
\label{CR2}\end{equation}
where
\begin{equation*}
\mathrm{I}_\rho:=\E\int_0^t\int_{(\T^N)^2}\int_{\R^2} f_1 \bar f_2 (a(\xi)\cdot\nabla_x +a(\zeta)\cdot\nabla_y)\alpha d\xi d\zeta dx  dy ds
\end{equation*}
and
\begin{multline*}
\mathrm{I}_\psi=\frac12\E\int_0^t\int_{(\T^N)^2}\int_{\R^2} \partial_\xi\alpha \bar f_2(s)\GG^2_1 d\nu^1_{x,s}(\xi)  d\zeta dx dy ds\\
-\frac12\E\int_0^t\int_{(\T^N)^2}\int_{\R^2}\partial_\zeta\alpha f_1(s) \GG^2_2 d\nu^2_{y,s}(\zeta)  d\xi dy dx ds\\
-\E\int_0^t\int_{(\T^N)^2}\int_{\R^2}\GG_{1,2} \alpha d\nu^1_{x,s}(\xi) d\nu^2_{y,s}(\zeta) dx dy.
\end{multline*}
Equation~\eqref{CR2} is indeed equation~\eqref{CR0} for $f_i^+$ since, by \eqref{ddvar},
\begin{equation*}
\mathrm{I}_\rho=\E\int_0^t\int_{(\T^N)^2}\int_{\R^2} f_1\bar f_2 (a(\xi)-a(\zeta))\cdot\nabla_x\alpha d\xi  d\zeta dx dy ds
\end{equation*}
and, by \eqref{ddvar} also and integration by parts,
\begin{align*}
\mathrm{I}_\psi&=\frac12\E\int_0^t\int_{(\T^N)^2}\int_{\R^2} \alpha (\GG^2_1+\GG^2_2-2\GG_{1,2}) d\nu^1_{x,s}\otimes\nu^2_{y,s}(\xi,\zeta) dx dy ds\\
&=\frac12\E\int_0^t\int_{(\T^N)^2}\int_{\R^2} \alpha \sum_{k\geq 0}|g_{k,1}(x,\xi)-g_{k,2}(y,\zeta)|^2 d\nu^1_{x,s}\otimes\nu^2_{y,s}(\xi,\zeta) dx dy ds.
\end{align*}
To obtain the result for $f_i^-$, we take $t_n \uparrow t$, write \eqref{CR0} for $f_i^+(t_n)$ and let
$n\to\infty$. \qed

\subsection{Uniqueness, reduction of generalized solution}

In this section we use Proposition~\ref{prop:CR} above to deduce the uniqueness of solutions and the reduction of generalized solutions to solutions.

\begin{theorem}[Uniqueness, Reduction] Let $u_0\in L^\infty(\T^N)$. Assume~\eqref{D0}-\eqref{D1}. 
Then we have the following results:
\begin{enumerate}
\item there is at most one solution with initial datum $u_0$ to \eqref{stoSCL}. 
\item If $f$ is a ge\-ne\-ra\-lized solution to \eqref{stoSCL} with initial 
datum $f_0$ \emph{at equilibrium:} $f_0=\mathbf{1}_{u_0>\xi}$, then there exists a solution $u$ to \eqref{stoSCL} with initial datum 
$u_0$ such that $f(x,t,\xi)=\mathbf{1}_{u(x,t)>\xi}$ a.s., for a.e. $(x,t,\xi)$. 
\item if $u_1$, $u_2$ are two solutions to \eqref{stoSCL} associated to the initial data $u_{1,0}$, $u_{2,0}\in L^\infty(\T^N)$ respectively, then
\begin{equation}
\E\|(u_1(t)-u_2(t))^+\|_{L^1(\T^N)}\leq\|(u_{1,0}-u_{2,0})^+\|_{L^1(\T^N)},
\label{L1compadd}\end{equation}
for all $t\in[0,T]$. This implies the $L^1$-contraction property, and the comparison principle for solutions.
\end{enumerate} 
\label{th:Uadd}\end{theorem}

\begin{corollary}[Continuity in time]Let $u_0\in L^\infty(\T^N)$. Assume~\eqref{D0}-\eqref{D1}. Then, 
for every $p\in[1,+\infty)$, the solution $u$ to \eqref{stoSCL} with initial datum $u_0$ has a representative in $L^p(\Omega;L^\infty(0,T;L^p(\T^N)))$ with almost-sure
continuous trajectories in $L^p(\T^N)$.
\label{cor:timecontinuity}\end{corollary}

\begin{remark}[Uniqueness of the kinetic measure] Let $f$ and $\check{f}$ be two ge\-ne\-ra\-lized solution to \eqref{stoSCL} with initial 
datum $f_0$ at equilibrium, $f_0=\mathbf{1}_{u_0>\xi}$. By Theorem~\ref{th:Uadd}, we have $f=\check{f}$. It follows from \eqref{eq:kineticfpre} that the associated random measures $m$ and $\check{m}$ satisfy: for all $\varphi\in C^1_c(\T^N\times\R)$, for all $t\in[0,T]$, almost-surely, 
\begin{equation}\label{mcheckm}
m(\partial_\xi\varphi)([0,t])=\check{m}(\partial_\xi\varphi)([0,t]).
\end{equation}
At fixed $\varphi$, the two functions of $t$ in \eqref{mcheckm} are c{\`a}dl{\`a}g. Therefore \eqref{mcheckm} is satisfied for all $\varphi\in C^1_c(\T^N\times\R)$, almost-surely, for all $t\in[0,T]$. By an argument of density (as in the proof of Proposition~\ref{prop:LRlimits}), we obtain \eqref{mcheckm} almost-surely, for all $\varphi\in C^1_c(\T^N\times\R)$, for all $t\in[0,T]$.  This implies: almost-surely, $\partial_\xi m=\partial_\xi \check{m}$. By \eqref{eq:massm}, the two measures have the same total mass almost-surely. Consequently, almost-surely, $m=\check{m}$.
\label{rk:uniqmmm}\end{remark}

{\bf Proof of Theorem~\ref{th:Uadd}.} Consider first the additive case: $\Phi(x,u(x))$ independent on $u(x)$. Let $f_i
$, $i=1,2$ be two generalized solutions to \eqref{stoSCL}. Then, we use \eqref{CR0} with $g_k$ 
independent on $\xi$ and $\zeta$. By \eqref{D1}, the last term $\mathrm{I}_\psi$ is bounded by
\begin{equation*}
\frac{t D_1}{2}\|\psi\|_{L^\infty}\int_{(\T^N)^2}|x-y|^2\rho(x-y)dx dy.
\end{equation*}
We then take $\psi:=\psi_{\delta}$ and $\rho=\rho_\eps$ where $(\psi_\delta)$ and $(\rho_\eps)$ are 
approximations to the identity on $\R$ and $\T^N$ respectively, \textit{i.e.}
$$
\psi_{\delta}(\xi)=\frac{1}{\delta}\psi\left(\frac{\xi}{\delta}\right),\quad\rho_\eps(x)=\frac{1}{\eps^N}\rho\left(\frac{x}{\eps}\right),
$$
where $\psi$ and $\rho$ are some given smooth probability densities on $\R$ and $\T^N$ respectively,
to obtain
\begin{equation}
\mathrm{I}_\psi\leq\frac{t D_1}{2}\eps^2\delta^{-1}.
\label{U1}\end{equation}

Denote by $\nu^{i,\pm}_{x,t}$ the Young measures associated to $f_i^\pm$, $i\in\{1,2\}$. By a computation similar to \eqref{psifbarf}, we have, almost-surely, for almost all $x,y\in\T^N$,
\begin{equation}\label{f1f2dd1}
\int_{\R} f_1^\pm(x,t,\xi)\bar f_2^\pm(y,t,\xi) d\xi=\int_{\R^2}(u-v)^+ d\nu^{1,\pm}_{x,t}(u) d\nu^{2,\pm}_{y,t}(v).
\end{equation}
By \eqref{psifbarf}, we have also
\begin{equation}\label{f1f2dd2}
\int_{\R^2}\psi_\delta(\xi-\zeta)  f_1^\pm(x,t,\xi)\bar f_2^\pm(y,t,\zeta)d\xi d\zeta
=\int_{\R^2}\psi_{2,\delta}(u-v) d\nu^{1,\pm}_{x,t}(u) d\nu^{2,\pm}_{y,t}(v),
\end{equation}
where
$$
\psi_{2,\delta}(\xi)=\int_{-\infty}^\xi\psi_{1,\delta}(\zeta)d\zeta,\quad \psi_{1,\delta}(\xi)=\int_{-\infty}^\xi\psi_\delta(\zeta)d\zeta.
$$
Assume that $\psi$ is supported in $(0,1)$. Then $\psi_{2,\delta}(\xi)=0$ if $\xi\leq 0$ and, for $\xi>0$,
\begin{equation}\label{f1f2dd3}
\xi^+-\psi_{2,\delta}(\xi)=\int_0^{\xi^+}\int_{\zeta/\delta}^{+\infty}\psi(u)du d\zeta=\int_0^{+\infty}\xi^+\wedge(\delta u)\psi(u)du.
\end{equation}
Using \eqref{f1f2dd3} in \eqref{f1f2dd1}, \eqref{f1f2dd2} gives
\begin{align*}
0&\leq \int_{\R} f_1^\pm(x,t,\xi)\bar f_2^\pm(y,t,\xi) d\xi-\int_{\R^2}\psi_\delta(\xi-\zeta)  f_1^\pm(x,t,\xi)\bar f_2^\pm(y,t,\zeta)d\xi d\zeta\\
&\leq \int_{\R^2}\int_0^{+\infty}(u-v)^+\wedge(\delta\zeta)\psi(\zeta)d\zeta d\nu^{1,\pm}_{x,t}(u) d\nu^{2,\pm}_{y,t}(v).
\end{align*}
Since $(u-v)^+\wedge(\delta\zeta)\leq |u|\wedge(\delta\zeta)+|v|\wedge(\delta\zeta)$, we have
\begin{align*}
0&\leq \int_{\R} f_1^\pm(x,t,\xi)\bar f_2^\pm(y,t,\xi) d\xi-\int_{\R^2}\psi_\delta(\xi-\zeta)  f_1^\pm(x,t,\xi)\bar f_2^\pm(y,t,\zeta)d\xi d\zeta\\
&\leq \int_0^{+\infty}\left(\int_\R |\xi|d\nu^{1,\pm}_{x,t}(\xi)+\int_\R |\xi|d\nu^{2,\pm}_{y,t}(\xi)\right)\wedge(\delta\zeta)\psi(\zeta)d\zeta.
\end{align*}
It follows that 
\begin{align}
&\left|\int_{(\T^N)^2}\int_{\R}\rho_\eps(x-y) f_1^\pm(x,t,\xi)\bar f_2^\pm(y,t,\xi) d\xi dx dy\right.\nonumber\\
&\hspace*{3cm}\left.-\int_{(\T^N)^2}\int_{\R^2}\rho_\eps(x-y)\psi_\delta(\xi-\zeta)  f_1^\pm(x,t,\xi)\bar f_2^\pm(y,t,\zeta)d\xi d\zeta dx dy\right|\nonumber\\
\leq& \int_0^{+\infty}\left(\int_{\T^N}\int_\R |\xi| (d\nu^{1,\pm}_{x,t}(\xi)+d\nu^{2,\pm}_{x,t}(\xi)) dx \right)\wedge(2\delta\zeta)\psi(\zeta)d\zeta\nonumber\\
\leq &\int_0^{+\infty}\left(C_{1,1}^\pm(\omega)+C_{1,2}^\pm(\omega)  \right)\wedge(2\delta\zeta)\psi(\zeta)d\zeta.\label{f1f2dd4}
\end{align}
We have used  \eqref{eq:integrabilityf}-\eqref{eq:integrabilityfpm} (with a constant $C_{1,i}^\pm$ for $\nu^{i,\pm}$) to obtain \eqref{f1f2dd4}. When $\eps\to0$, we have
\begin{align}
\left|\int_{(\T^N)^2}\int_{\R}\rho_\eps(x-y) f_1^\pm(x,t,\xi)\right. &\left.\bar f_2^\pm(y,t,\xi) d\xi dx dy-\int_{\T^N}\int_{\R} f_1^\pm(x,t,\xi)\bar f_2^\pm(x,t,\xi) d\xi dx\right|\nonumber\\
\leq& \sup_{|z|<\eps}\int_{\T^N}\int_{\R}f_1^\pm(x,t,\xi)\left|\bar f_2^\pm(x-z,t,\xi)-\bar f_2^\pm(x,t,\xi)\right|  d\xi dx\nonumber\\
\leq& \sup_{|z|<\eps}\int_{\T^N}\int_{\R}\left|\chi_{f_2^\pm}(x-z,t,\xi)-\chi_{f_2^\pm}(x,t,\xi)\right|  d\xi dx.\label{f1f2dd5}
\end{align}
Consequently (see \eqref{f1f2dd4}, \eqref{f1f2dd5}),
\begin{multline*}
\lim_{\eps,\delta\to0}\int_{(\T^N)^2}\int_{\R^2}\rho_\eps(x-y)\psi_\delta(\xi-\zeta)  f_1^\pm(x,t,\xi)\bar f_2^\pm(y,t,\zeta)d\xi d\zeta dx dy\\
=\int_{\T^N}\int_{\R} f_1^\pm(x,t,\xi)\bar f_2^\pm(x,t,\xi) d\xi dx,
\end{multline*}
for all $\omega\in\hat\Omega$. 
We apply the estimate \eqref{boundCRleft}. We have the uniform bounds
$$
\|\rho_\eps\|_{L^1(\T^N)}=1,\quad \|\psi_\delta\|_{L^1(\R)}=1,\quad R=\delta\leq 1.
$$
Consequently, we may apply the Lebesgue dominated convergence theorem: we obtain
\begin{multline}
\E\int_{\T^N}\int_{\R} f_1^\pm(x,t,\xi)\bar f_2^\pm(x,t,\xi) dx d\xi\nonumber\\
\leq\E\int_{(\T^N)^2}\int_{\R^2} \rho_\eps(x-y)\psi_\delta(\xi-\zeta)  f_1^\pm(x,t,\xi)\bar f_2^\pm(y,t,\zeta)d\xi d\zeta dx dy +\eta_t(\eps,\delta),\label{U2}
\end{multline}
where $\lim_{\eps,\delta\to 0}\eta_t(\eps,\delta)=0$. 
We need now a bound on the term $\mathrm{I}_\rho$. Since $a'=A''$ has at most polynomial growth, there exists $C\geq 0$, $p\geq 1$, such that
\begin{equation*}
|a(\xi)-a(\zeta)|\leq \Gamma(\xi,\zeta)|\xi-\zeta|,\quad \Gamma(\xi,\zeta)=C(1+|\xi|^{p-1}+|\zeta|^{p-1}).
\end{equation*}
This gives
\begin{equation*}
|\mathrm{I}_\rho|\leq \E\int_0^t\int_{(\T^N)^2}\int_{\R^2}f_1\bar f_2\Gamma(\xi,\zeta)|\xi-\zeta|\psi_\delta(\xi-\zeta)|\nabla_x\rho_\eps(x-y)| d\xi d\zeta dx dy d\sigma.
\end{equation*}
By integration by parts with respect to $(\xi,\zeta)$, we deduce
\begin{equation*}
|\mathrm{I}_\rho|\leq \E\int_0^t\int_{(\T^N)^2}\int_{\R^2}\Upsilon(\xi,\zeta)d\nu^1_{x,\sigma}\otimes\nu^2_{y,\sigma}(\xi,\zeta)|\nabla_x\rho_\eps(x-y)|  dx dy d\sigma,
\end{equation*}
where
\begin{equation*}
\Upsilon(\xi,\zeta)=\int_\zeta^{+\infty}\int_{-\infty}^\xi\Gamma(\xi',\zeta')|\xi'-\zeta'|\psi_\delta(\xi'-\zeta') d\xi'd\zeta'.
\end{equation*}
It is shown below that $\Upsilon$ admits the bound
\begin{equation}
\Upsilon(\xi,\zeta)\leq C(1+|\xi|^{p}+|\zeta|^{p})\delta.
\label{ffp}\end{equation}
Since $\nu^1$ and $\nu^2$ vanish at infinity, \textit{cf.} \eqref{eq:integrabilityf}, we then obtain, for a given constant $C_p$,
\begin{equation*}
|\mathrm{I}_\rho|\leq t C_p\delta \left(\int_{\T^N}|\nabla_x\rho_\eps(x)|dx\right).
\end{equation*}
It follows that, for possibly a different $C_p$,
\begin{equation}
|\mathrm{I}_\rho|\leq t C_p  \delta\eps^{-1}.
\label{U3}\end{equation}
We then gather \eqref{U1}, \eqref{U3} and \eqref{CR0} to deduce for $t\in [0,T]$
\begin{equation}
\E\int_{\T^N}\int_{\R} f_1^\pm(t)\bar f_2^\pm(t) dx d\xi\leq \int_{\T^N}\int_{\R} f_{1,0}\bar f_{2,0}dx d\xi+r(\eps,\delta),
\label{U4}\end{equation}
where the remainder $r(\eps,\delta)$ is
$
\ds r(\eps,\delta)=T C_p\delta\eps^{-1}+\frac{T D_1}{2}\eps^2\delta^{-1}+\eta_t(\eps,\delta)-\eta_0(\eps,\delta).
$
Taking $\delta=\eps^{4/3}$ and letting $\eps\to 0$ gives 
\begin{equation}
\E\int_{\T^N}\int_{\R} f_1^\pm(t)\bar f_2^\pm(t) dx d\xi\leq\int_{\T^N}\int_{\R} f_{1,0}\bar f_{2,0} dx d\xi. 
\label{L1compadd0}\end{equation}
Assume that $f$ is a generalized solution to \eqref{stoSCL} with initial datum $\mathbf{1}_{u_0>\xi}$. Since
 $f_0$ is the (translated) Heavyside function $\mathbf{1}_{u_0>\xi}$, we have the identity $f_0\bar 
 f_0=0$. Taking $f_1=f_2=f$ in \eqref{L1compadd0}, we deduce $f^+(1-f^+)=0$ a.e., {\it i.e.} 
$f^+\in\{0,1\}$ a.e. The fact that $-\partial_\xi f^+$ is a Young measure then gives the conclusion: 
indeed, by Fubini's Theorem, for any $t\in [0,T]$, there is a set $E_t$ of full measure in $\T^N\times
\Omega$ such that, for  $(x,\omega)\in E_t$, $f^+(x,t,\xi,\omega)\in\{0,1\}$ for a.e. $\xi\in\R$. Let 
$$
\tilde E_t=E_t\cap(\T^N\times\hat\Omega).
$$
The set $\tilde E_t$ is of full measure in $\T^N\times\Omega$. For $(x,\omega)\in\tilde E_t$, 
$-\partial_\xi f^+(x,t,\cdot,\omega)$ is a probability measure on $\R$. 
Therefore 
$f^+(t,x,\xi,\omega)=\mathbf{1}_{u(x,t,\omega)>\xi}$ for a.e. $\xi\in\R$, where
$u(x,t,\omega)=\int_\R (f^+(x,t,\xi,\omega)-\mathbf{1}_{0>\xi}) d\xi$. We have a similar result for $f^-$ (this will be used in the proof of Corollary~\ref{cor:timecontinuity}). Proposition~\ref{prop:gestoes} implies that $u$ is a solution in the sense of Definition \ref{defkineticsol}. Since $f=f^+$ (recall the convention of notation introduced at the beginning of Section~\ref{sec:doubling}), this shows the reduction of generalized solutions to solutions.
If now $u_1$ and $u_2$ are two solutions to \eqref{stoSCL}, we deduce from \eqref{L1compadd0} with $f_i=\mathbf{1}_{u_i>\xi}$ and from the identity
\begin{equation*}
\int_\R \mathbf{1}_{u_1>\xi}\overline{\mathbf{1}_{u_2>\xi}}d\xi=(u_1-u_2)^+,
\end{equation*}
the contraction property \eqref{L1compadd}.
\medskip

In the multiplicative case ($\Phi$ depending on $u$), the reasoning is similar, except that there is an additional term in the bound on $\mathrm{I_\psi}$. More precisely, by the hypothesis~\eqref{D1} we obtain in place of \eqref{U1} the estimate 
\begin{equation*}
\ds\mathrm{I_\psi}\leq \frac{T D_1}{2}\eps^2\delta^{-1}
+\frac{D_1}{2}\mathrm{I}^h_\psi,
\end{equation*} 
where
\begin{equation*}
\mathrm{I}_\psi^h=\E\int_0^t\int_{(\T^N)^2}\rho_\eps\int_{\R^2} \psi_\delta(\xi-\zeta)|\xi-\zeta|h(|\xi-\zeta|)d\nu^1_{x,\sigma}\otimes\nu^2_{y,\sigma}(\xi,\zeta)dx dy d\sigma.
\end{equation*}
Choosing $\psi_\delta(\xi)=\delta^{-1}\bar\Psi(\delta^{-1}\xi)$ with $\bar\Psi$ compactly supported in $(0,1)$ gives
\begin{equation}
\mathrm{I_\psi}\leq \frac{T D_1}{2}\eps^2\delta^{-1}+\frac{TD_1C_\psi h(\delta)}{2},\quad C_\psi:=\sup_{\xi\in(0,1)}\|\xi\bar\Psi(\xi)\|.
\label{U1bis}\end{equation}
We deduce \eqref{U4} with a remainder term $\ds r'(\eps,\delta):=r(\eps,\delta)+\frac{TD_1C_\psi h(\delta)}{2}$ and conclude the proof as in the additive case. There remains to prove \eqref{ffp}: setting $\xi''=\xi'-\zeta'$, we have
\begin{align*}
\Upsilon(\xi,\zeta)=&\int_\zeta^{+\infty}\int_{|\xi''|<\delta,\xi''<\xi-\zeta'}\Gamma(\xi''+\zeta',\zeta')|\xi''|\psi_\delta(\xi'') d\xi''d\zeta'\\
\leq&C\int_{\zeta}^{\xi+\delta}\max_{|\xi''|<\delta,\xi''<\xi-\zeta'}\Gamma(\xi''+\zeta',\zeta') d\zeta'\; \delta\\
\leq &C \int_{\zeta}^{\xi+\delta}(1+|\xi|^{p-1}+|\zeta'|^{p-1}) d\zeta'\; \delta,
\end{align*}
which gives \eqref{ffp}. \qed
\bigskip

{\bf Proof of Corollary~\ref{cor:timecontinuity}.} 
We use the notations and the results of Proposition~\ref{prop:LRlimits}. We fix $p\in[1,+\infty)$. Both $(f(t))_{t\in[0,T]}$ and $(f^-(t))_{t\in[0,T]}$ are generalized solutions to \eqref{stoSCL} associated to the initial datum $\mathbf{1}_{u_0>\xi}$ (we use \eqref{eq:kineticfmoinspre} here). By Theorem~\ref{th:Uadd}, they are at equilibrium: $f(t)=\mathbf{1}_{u(t)>\xi}$, $f^-(t)=\mathbf{1}_{u^-(t)>\xi}$. By Proposition~\ref{prop:LRlimitsu}, for all $\omega\in\hat{\Omega}$, the map $t\mapsto u(t)$ from $[0,T]$ to $L^p(\T^N)$ is continuous from the right. Similarly, $t\mapsto u^-(t)$ is continuous from the left. By Proposition~\ref{prop:equilibrium}, the fact that $f^-$ is at equilibrium has the following consequence: at every $t\in(0,T]$, $f(t)=f^-(t)$. In particular, we have $u=u^-$ and thus, almost-surely, $u$ is continuous from $[0,T]$ to $L^p(\T^N)$.
\qed\medskip

We apply \eqref{L1compadd} to infer an $L^\infty$ bound on solutions to \eqref{stoSCL} in the particular case of a multiplicative noise with compact support.

\begin{theorem}[$L^\infty$ bounds] Assume~\eqref{D0}-\eqref{D1} and
\begin{equation}\label{multiplicativeNoise}
g_k(x,u)=0,\quad\forall |u|\geq 1,
\end{equation}
for all $x\in\T^N$, $k\geq 1$. Let $u_0\in L^\infty(\T^N)$ satisfy $-1\leq u_0\leq 1$ almost everywhere. 
Then, for all $t\geq 0$, the solution $u$ to \eqref{stoSCL} with initial datum $u_0$ satisfies: almost-surely,
\begin{equation}\label{Linftyu}
-1\leq u(x,t)\leq 1,
\end{equation}
a.e. in $\T^N$.
\label{th:Linftyu}\end{theorem}

\textbf{Proof of Theorem~\ref{th:Linftyu}.} We use \eqref{L1compadd} to compare $u$ to the two particular constant solutions $(x,t)\mapsto -1$ and $(x,t)\mapsto 1$. \qed

\section{Convergence of approximate solutions}\label{sec:martingale}

In this section, we develop the tools required for the proof of convergence of a certain type of approximate solutions to \eqref{stoSCL}. The basic principle is to generalize the notion of solution introduced in Definition~\ref{defkineticsol}. Indeed, this facilitates the proof of e\-xis\-ten\-ce/convergence. In a second step a result of reduction (or ``rigidity result"), which asserts that a generalized solution is a solution is used. This principle is of much use in the deterministic theory of conservation laws (\textit{cf.} \cite{Diperna83a} with the use of ``measure-valued entropy solutions", \cite{EGH00} with the use of ``entropy process solutions", \cite{PerthameBook} with the use of kinetic solutions as defined here). We have already introduced a generalization of the notion of solution in Definition~\ref{d4}, and have proved a result of reduction in Theorem~\ref{th:Uadd}. Here we will work mainly on the probabilistic aspects of the questions. We will have to consider ``solutions in law", or "martingale solutions" (see the comment after Theorem~\ref{th:martingalesol} for more explanations about the terminology). The plan of this section is the following one: in Section~\ref{sec:approxges}, we define the notion of approximate generalized solution. In Section~\ref{sec:Mintsto}, we give a martingale characterization of the stochastic integral. In Section~\ref{sec:tight}, we give some tightness results on sequences of approximate generalized solutions. The main result, Theorem~\ref{th:martingalesol}, which shows the convergence of a sequence of approximate generalized solutions to a martingale generalized solution, is proved in Section~\ref{sec:proofThmartingale}. Eventually, we obtain a result of pathwise convergence in Section~\ref{sec:pathCV}.

\subsection{Approximate generalized solutions}\label{sec:approxges}
 
Let $\order$ be an integer fixed once and for all.

\begin{definition}[Approximate generalized solutions] Let $f_0^n\colon\T^N\times\R\to[0,1]$ be some kinetic functions. Let $(f^n(t))_{t\in[0,T]}$ be a sequence of $L^\infty(\T^N\times\R;[0,1])$-valued processes. Assume that the functions $f^n(t)$, and the associated Young measures $\nu^n_t=-\partial_\xi\varphi f^n(t)$ are satisfying 
item \ref{item:1d4}, \ref{item:1terd4}, \ref{item:2d4}, in Definition~\ref{d4} and Equation~\eqref{eq:kineticfpre} up to an error term, \textit{i.e.}: for all $\varphi\in C^\order_c(\T^N\times\R)$, there exists an adapted process $\eps^n(t,\varphi)$, with $t\mapsto\eps^n(t,\varphi)$ almost-surely continuous such that
\begin{equation}\label{epsto0}
\lim\limits_{n\to+\infty}\sup_{t\in[0,T]}\left|\eps^n(t,\varphi)\right|=0\mbox{ in probability,}
\end{equation} 
and there exists some random measures $m^n$  with first moment \eqref{Firstm}, such that, for all $n$, for all $\varphi\in C^\order_c(\T^N\times\R)$, for all $t\in[0,T]$, almost-surely,
\begin{align}
\<f^n(t),\varphi\>=&\<f^n_0,\varphi\>+\int_0^t \<f^n(s),a(\xi)\cdot\nabla_x\varphi\>ds
\nonumber\\
&+\int_0^t\int_{\T^N}\int_\R g_k(x,\xi)\varphi(x,\xi)d\nu^n_{x,s}(\xi) dxd\beta_k(s)+\eps^n(t,\varphi) \nonumber\\
&+\frac{1}{2}\int_0^t\int_{\T^N}\int_\R \GG^2(x,\xi)\partial_\xi\varphi(x,\xi)d\nu^n_{x,s}(\xi) dx ds
-m^n(\partial_\xi\varphi)([0,t]).
\label{eq:kineticfpreappt}
\end{align}
Assume also $f^n(0)=f^n_0$. Then we say that $(f^n)$ is a sequence of approximate generalized solutions to~\eqref{stoSCL} with initial datum $f_0^n$.
\label{def:appsol}\end{definition}

\subsection{Martingale characterization of the stochastic integral}\label{sec:Mintsto}
In order to pass to the limit in an equation such as \eqref{eq:kineticfpreappt}, we will first characterize \eqref{eq:kineticfpreappt} in terms of a martingale problem, and then we will use martingale methods to pass to the limit. In the present section, we give the characterization of \eqref{eq:kineticfpreappt} in terms of a martingale problem, see Proposition~\ref{prop:3martingales} and Proposition~\ref{prop:3martingalesCV} below. We refer to \cite[Example~1.4, p.143]{JacodShiryaev03} for characterization of the standard Wiener Process in terms of a martingale problem. In the context of SDEs and SPDEs, such kind of characterizations have been applied in \cite{Ondrejat10,BrzezniakOndrejat11,HofmanovaSeidler12,Hofmanova13b,DebusscheHofmanovaVovelle15} in particular.
\medskip

Let us define the stochastic integrands
\begin{equation}\label{hphikn}
h^n_{\varphi,k}(t)=\int_{\T^N}\int_\R g_k(x,\xi)\varphi(x,\xi)d\nu^n_{x,t}(\xi)dx,\quad h^n_\varphi(t)=\left(h^n_{\varphi,k}(t)\right)_{k\geq 1},
\end{equation} 
and the stochastic integrals
\begin{equation}\label{def:Mnphi}
M_\varphi^n(t)=\sum_{k\geq 1}\int_0^t\int_{\T^N}\int_\R g_k(x,\xi)\varphi(x,\xi)d\nu^n_{x,s}(\xi)dxd\beta_k(s).
\end{equation}
By Lemma~\ref{lem:admint}, we have $h^n_\varphi\in L^2_\mathcal{P}([0,T]\times\Omega;l^2(\N^*))$ for all $n$, $\varphi$. Using It\={o}'s Formula, we deduce from \eqref{def:Mnphi} the following statement.

\begin{proposition}\label{prop:3martingales} Let $(f^n)$ be a sequence of approximate generalized solutions to~\eqref{stoSCL} with initial datum $f_0^n$. Let $\varphi\in C^\order_c(\T^N\times\R)$. Let $M^n_\varphi(t)$ be defined by \eqref{def:Mnphi} and $h^n_{\varphi,k}(t)$ by \eqref{hphikn}. Then the processes
\begin{equation}\label{3martingales}
M^n_\varphi(t),\quad M^n_\varphi(t)\beta_k(t)-\int_0^t h^n_{\varphi,k}(s)ds,\quad |M^n_\varphi(t)|^2-\int_0^t\|h^n_\varphi(s)\|_{l^2(\N^*)}ds,
\end{equation}
are $(\mathcal{F}_t)$-martingales.
\end{proposition}

What will interest us is the reciprocal statement.

\begin{proposition}\label{prop:3martingalesCV} Let $h\in L^2_\mathcal{P}([0,T]\times\Omega;l^2(\N^*))$. Let $X(t)$ be a stochastic process starting from $0$ such that the processes
\begin{equation}\label{3martingalesCV}
X(t),\quad X(t)\beta_k(t)-\int_0^t h_k(s) ds,\quad |X(t)|^2-\int_0^t \|h(s)\|_{l^2(\N^*)}^2ds
\end{equation}
are $(\mathcal{F}_t)$-martingales. Then
\begin{equation}\label{Xintsto}
X(t)=\sum_{k\geq 1}\int_0^t h_k(s) d\beta_k(s),
\end{equation}
for all $t\in[0,T]$.
\end{proposition}

\textbf{Proof of Proposition~\ref{prop:3martingalesCV}.} The proof can be found in \cite[Proposition~A.1]{Hofmanova13b}. Let us give some details about it. We first claim that the following identity is satisfied:
\begin{equation}\label{eq:crossmartingale}
\E\left[(X(t)-X(s))\int_s^t \theta(\sigma)d\beta_k(\sigma)d\sigma-\int_s^t h_k(\sigma)\theta(\sigma)d\sigma\Big|\mathcal{F}_s\right]=0
\end{equation}
for all $0\leq s\leq t\leq T$, all $k\geq 1$ and all $\theta\in L^2_\mathcal{P}([0,T]\times\Omega)$. The proof consists in approximating $\theta$ on the interval $[s,t]$ by predictable simple functions. It is similar to a computation of quadratic variation. Note that \eqref{eq:crossmartingale} uses only the fact that
\begin{equation*}
X(t),\quad X(t)\beta_k(t)-\int_0^t h_k(s) ds
\end{equation*}
are $(\mathcal{F}_t)$-martingales. We apply \eqref{eq:crossmartingale} with $s=0$ and $\theta=h_k$ and sum over $k$ to obtain
\begin{equation}\label{crosscross}
\E[X(t)\bar X(t)]=\E\int_0^t \|h(s)\|_{l^2(\N^*)}^2 ds,\quad \bar X(t):=\sum_{k\geq 1}\int_0^t h_k(s) d\beta_k(s).
\end{equation}
This gives the expression of the cross-product when we expand the term $\E |X(t)-\bar X(t)|^2$. Using the fact that
$$
|X(t)|^2-\int_0^t \|h(s)\|_{l^2(\N^*)}^2 ds
$$
is a $(\mathcal{F}_t)$-martingale and applying It\={o}'s Isometry to $\E|\bar X(t)|^2$ shows that the square terms are also given by
$$
\E|X(t)|^2=\E|\bar X(t)|^2=\int_0^t \|h(s)\|_{l^2(\N^*)}^2ds.
$$
It follows that $X(t)=\bar X(t)$. \qed


\subsection{Tightness}\label{sec:tight}

Let $(f^n)$ be a sequence of approximate generalized solutions, in the sense of Definition~\ref{def:appsol}. Recall that $\mathcal{Y}^1$ is the notation for the set of Young measures on $\T^N\times[0,T]\times\R$ (\textit{cf.} Proposition~\ref{prop:altdefY}) and that $\mathcal{M}_b(\T^N\times[0,T]\times\R)$ is the notation for the set of bounded Borel measures on $\T^N\times[0,T]\times\R$ while $\mathcal{M}^+_b(\T^N\times[0,T]\times\R)$ is the subset of non-negative measures. Let $\nu^n$ be the Young measure associated to $f^n$ ($\nu^n=-\partial_\xi f^n$). The law of $\nu^n$ is a probability measure on the space $\mathcal{Y}^1$. We will see in Section~\ref{sec:compactnun} that, under a natural a priori bound, see \eqref{eq:integrabilityfn}, the sequence $(\mathrm{Law}(\nu^n))$ is tight in $\mathcal{Y}^1$. In Section~\ref{sec:compactmn}, this is the sequence $(\mathrm{Law}(m^n))$ that we will analyse. We show under \eqref{Boundmn} and \eqref{inftymn} that it is tight in $\mathcal{M}^+_b(\T^N\times[0,T]\times\R)$ (see, more specifically, Proposition~\ref{prop:stabextractm}).\medskip

We also need to analyse the tightness of $(\<f_n(t),\varphi\>)$ in the Skorokhod space $D([0,T])$: this is done in Section~\ref{sec:cadlag}.

\subsubsection{Compactness of the Young measures}\label{sec:compactnun}

In this section, we will use the following notions: we say that a sequence $(\nu^n)$ of $\mathcal{Y}^1$ converges to $\nu$ in $\mathcal{Y}^1$ if \eqref{cvYoungMeasure} is satisfied. A \textbf{random Young measure} is by definition a $\mathcal{Y}^1$-valued random variable.

\begin{proposition}  Let $(f^n)$ be a sequence of approximate generalized solutions to~\eqref{stoSCL} with initial datum $f_0^n$. Assume that the following bound is satisfied: for all $p\in[1,+\infty)$, there exists $C_p\geq 0$ independent on $n$ such that $\nu^n:=-\partial_\xi f^n$ satisfies 
\begin{equation}
\E\left[\sup_{t\in[0,T]}\int_{\T^N}\int_\R|\xi|^p d\nu^n_{x,t}(\xi) dx\right]\leq C_p,
\label{eq:integrabilityfn}\end{equation}
Then, there exists a probability space $(\tilde\Omega,\tilde{\mathcal{F}},\tilde\P)$ and some random Young measures $\tilde\nu^n$, $\tilde\nu$, such that
\begin{enumerate}
\item $\tilde\nu^n$ has the same law as $\nu^n$,
\item $\tilde\nu$ satisfies
\begin{equation}
\tilde\E\left(\sup_{J\subset[0,T]}\frac{1}{|J|}\int_J\int_{\T^N}\int_\R|\xi|^p d\tilde\nu_{x,t}(\xi) dx dt \right) \leq C_p,
\label{eq:integrabilityftilde}\end{equation}
where the supremum in \eqref{eq:integrabilityftilde} is a countable supremum over all open intervals $J\subset[0,T]$ with rational extremities,
\item up to a subsequence still denoted $(\tilde\nu^n)$, there is $\tilde\P$-almost-sure convergence of $(\tilde\nu^n)$ to $\tilde\nu$ in $\mathcal{Y}^1$.
\end{enumerate}
Furthermore, if $\tilde f^n,\tilde f\colon\T^N\times[0,T]\times\R\times\tilde\Omega\to[0,1]$ are defined by 
$$
\tilde f^n(x,t,\xi)=\tilde\nu^n_{x,t}(\xi,+\infty),\quad \tilde f(x,t,\xi)=\tilde\nu_{x,t}(\xi,+\infty),
$$
then $\tilde f^n\to\tilde f$ in $L^\infty(\T^N\times[0,T]\times\R)$-weak-* $\tilde\P$-almost-surely, $\tilde f$ being a kinetic function.
\label{prop:stabextractnu}\end{proposition}

\textbf{Proof of Proposition~\ref{prop:stabextractnu}.} Note first that \eqref{eq:integrabilityfn} yields
\begin{equation}
\E\left(\int_0^T\int_{\T^N}\int_\R|\xi|^p d\nu^n_{x,t}(\xi)  dx dt \right) \leq C_p T.
\label{eq:integrabilityfnint}
\end{equation}
For $R>0$, $p\geq 1$, let us denote by $K_{R,p}$ the set of Young measures $\nu\in\mathcal{Y}^1$ such that 
\begin{equation*}
\int_0^T\int_{\T^N}\int_\R |\xi|^p d\nu_{x,t}(\xi)dx dt\leq R.
\end{equation*}
By \cite[Theorem~4.3.2, Theorem~4.3.8,Theorem~2.1.3]{CastaingRaynauddeFitteValadier04}, the set  $K_{R,p}$ is compact in $\mathcal{Y}^1$ for the $\tau^W_{\mathcal{Y}^1}$-topology, which is metrizable, \cite[Theorem~2.3.1]{CastaingRaynauddeFitteValadier04} and corresponds to the convergence \eqref{cvYoungMeasure}. By \eqref{eq:integrabilityfnint}, we have
$$
\P(\nu^n\notin K_{R,p})\leq\frac{C_p T}{R},
$$
which shows that the sequence $(\nu^n)$ of $\mathcal{Y}^1$-valued random variables is tight. The set $\mathcal{Y}^1$ endowed with the $\tau^W_{\mathcal{Y}^1}$-topology is Polish, \cite[Theorem~2.3.3]{CastaingRaynauddeFitteValadier04}: we can use the Prokhorov's metric, \cite[p.~72]{BillingsleyBook}. By Prokhorov's Theorem, \cite[Theorem~5.1]{BillingsleyBook}, there exists a $\mathcal{Y}^1$-valued random variable $\nu$ and a subsequence still denoted $(\nu^n)$ such that $(\nu^n)$ converges in probability to $\nu$. Since the map
$$
\psi_p\colon\mathcal{Y}^1\to[0,+\infty],\quad\nu\mapsto\sup_{J\subset[0,T]}\frac{1}{|J|}\int_J\int_{\T^N}\int_\R|\xi|^p d\nu_{x,t}(\xi) dx dt,
$$
are lower semi-continuous, we have 
$$
\E\psi_p(\nu)\leq\liminf_{n\to+\infty}\E\psi_p(\nu^n)\leq C_p
$$
by \eqref{eq:integrabilityfn} and, consequently, $\nu$ satisfies the condition 
\begin{equation}
\E\left(\sup_{J\subset[0,T]}\frac{1}{|J|}\int_J\int_{\T^N}\int_\R|\xi|^p d\nu_{x,t}(\xi) dx dt \right) \leq C_p.
\label{eq:integrabilityfff}\end{equation}
Let us now apply the Skorokhod Theorem \cite[p.~70]{BillingsleyBook}: there exists a probability space $(\tilde\Omega,\tilde{\mathcal{F}},\tilde\P)$ and some random variables $\tilde\nu^n$, $\tilde\nu$, such that
\begin{enumerate}
\item $\tilde\nu^n$ and $\tilde\nu$ have the same laws as $\nu^n$ and $\nu$ respectively,
\item up to a subsequence still denoted $(\tilde\nu^n)$, there is $\tilde\P$-almost-sure convergence of $(\tilde\nu^n)$ to $\tilde\nu$ in $\mathcal{Y}^1$.
\end{enumerate}
Since $\tilde\nu$ and $\nu$ have same laws, $\tilde\nu$ satisfies the bound \eqref{eq:integrabilityftilde}. If we apply Corollary~\ref{cor:kineticfunctions}, we obtain that $\tilde f^n\to\tilde f$ in $L^\infty(\T^N\times[0,T]\times\R)$-weak-* $\tilde\P$-almost-surely, $\tilde f$ being a kinetic function. \qed

\subsubsection{Compactness of the random measures}\label{sec:compactmn}

\begin{proposition} Let $(f^n)$ be a sequence of approximate generalized solutions to~\eqref{stoSCL} with initial datum $f_0^n$. Assume that
\begin{equation}\label{Boundmn}
\E m^n(\T^N\times[0,T]\times\R)\mbox{ is uniformly bounded},
\end{equation}
and that $m^n$ vanishes for large $\xi$ uniformly in $n$: if $B_R^c=\{\xi\in\R,|\xi|\geq R\}$, then
\begin{equation}
\lim_{R\to+\infty}\E m^n(\T^N\times[0,T]\times B_R^c)=0,
\label{inftymn}\end{equation}
uniformly in $n$. Then, there exists a probability space $(\tilde\Omega,\tilde{\mathcal{F}},\tilde\P)$ and some random measures $\tilde m^n$, $\tilde m\colon\tilde\Omega\to\mathcal{M}_b(\T^N\times[0,T]\times\R)$ such that
\begin{enumerate}
\item $\tilde m^n$ has the same law as $m^n$,
\item up to a subsequence still denoted $(\tilde m^n)$, there is $\tilde\P$-almost-sure convergence of $(\tilde m^n)$ to $\tilde m$ in $\mathcal{M}_b(\T^N\times[0,T]\times\R)$-weak-*.
\end{enumerate}
\label{prop:stabextractm}\end{proposition}

\textbf{Proof of Proposition~\ref{prop:stabextractm}.} Let $\eta\colon\R_+\to\R_+$ be defined by
$$
\eta(R)=\sup_{n\in\N}\E m^n(\T^N\times[0,T]\times B_R^c).
$$
Let $h$ be a fixed function on $\T^N\times[0,T]\times\R$, $h$ continuous, positive, integrable. Proving the statement for  the sequence of measures
$$
B\mapsto m^n(B)+\int_B h(x,t,\xi)dxdtd\xi
$$
is equivalent to prove the statement for the original sequence $(m^n)$. We will assume therefore that $\eta(R)>0$ for all $R\geq 0$ and that 
$$
\|m^n\|:=m^n(\T^N\times[0,T]\times\R)\geq\delta>0,
$$
where $\delta$ is independent on $n$. Let $\mu^n:=\frac{m^n}{\|m^n\|}$. We consider the random variables $X^n=(\mu^n,\|m^n\|)$, taking values in $\mathcal{P}^1(\T^N\times[0,T]\times\R)\times\R_+$, where $\mathcal{P}^1(\T^N\times[0,T]\times\R)$ is the set of probability measures on $\T^N\times[0,T]\times\R$. For $A>0$, let $K_A$ be the set of probability measures $\mu\in\mathcal{P}^1(\T^N\times[0,T]\times\R)$ such that
$$
\sup_{R>1}\frac{\mu(\T^N\times[0,T]\times B_R^c)}{\eta(R)}\leq A.
$$
Then $K_A$ is compact in $\mathcal{P}^1(\T^N\times[0,T]\times\R)$-weak-* by Prokhorov's Theorem and \eqref{inftymn}. Using the Markov Inequality, and the definition of $\eta(R)$, we obtain
$$
\P(\mu^n\notin K_A)\leq\frac{C}{A},
$$
where $C$ is independent on $n$: this shows that $(\mu^n)$ is tight in $\mathcal{P}^1(\T^N\times[0,T]\times\R)$ endowed with the topology of the weak convergence of probability measures. Similarly, using \eqref{Boundmn} and the Markov Inequality, we have
$$
\P(\|m^n\|>A)\leq\frac{C}{A},
$$
where $C$ is independent on $n$:, therefore $(\|m^n\|)$ is tight in $\R$. It follows that $(X^n)$ is tight in $\mathcal{P}^1(\T^N\times[0,T]\times\R)\times\R_+$ endowed with the product topology. This topology is separable, metrizable and there exists a compatible metric which turns the space into a complete space (we can take the sum of the Prokhorov's metric and of the usual metric on $\R_+$). Therefore we can apply the Skorokhod Theorem: there exists a probability space $(\tilde\Omega,\tilde{\mathcal{F}},\tilde\P)$ and some random variables $\tilde X^n=(\tilde\mu^n,\tilde\alpha^n)$, $\tilde X=(\tilde\mu,\tilde\alpha)$ such that $\tilde X^n$ has same law as $X^n$ and,  $\tilde\P$-almost-surely, $\tilde X^n\to\tilde X$ in $\mathcal{P}^1(\T^N\times[0,T]\times\R)\times\R_+$. Set $\tilde m^n=\tilde\alpha^n\tilde\mu^n$ and $\tilde m=\tilde\alpha\tilde\mu$. Then $\tilde m^n$ has the same law as $m^n$ and there is $\tilde\P$-almost-sure convergence of $(\tilde m^n)$ to $\tilde m$ in $\mathcal{M}_b(\T^N\times[0,T]\times\R)$-weak-*. \qed

\subsubsection{Tightness in the Skorokhod space}\label{sec:cadlag}

Let $D([0,T])$ denote the space of c\`adl\`ag functions on $[0,T]$. See \cite[VI.1]{JacodShiryaev03} and \cite[Chapter~3]{BillingsleyBook} for the definition of $D([0,T])$. Let $(f^n)$ be a sequence of approximate generalized solutions to~\eqref{stoSCL} with initial datum $f_0^n$. In Section~\ref{sec:MnM} below, where we analyse the convergence of $(f_n)$, it would be desirable to have a result of tightness of the processes $t\mapsto\<f^n(t),\varphi\>$ since they are random variables in $D([0,T])$ (here, $\varphi$ is a given test-function). It seems difficult to obtain such a result. The only fact which we can infer naturally from \eqref{eq:integrabilityfn}, \eqref{Boundmn}, \eqref{inftymn}, is that the sequence of processes
$$
t\mapsto \<f^n(t),\varphi\>+A^n_\varphi(t),\quad A^n_\varphi(t):=\<m^n,\partial_\xi\varphi\>([0,t]),
$$
is tight in $D([0,T])$, see Proposition~\ref{prop:Ctight} below. Showing additionally that $(A^n_\varphi)$ is tight in $D([0,T])$ seems impossible, however, if no additional properties of $(m^n)$ are known. Indeed, the weak convergence of
$
\mu^n:=\<m^n,\partial_\xi\varphi\>
$
to a measure $\mu$ on $[0,T]$ is not a sufficient condition for the convergence of $A^n_\varphi$ to $A(t)=\mu([0,t])$ in $D([0,T])$. Consider for example the case
$$
\mu^n=\delta_{t_*-s_n}+\delta_{t_*-\sigma_n},
$$
where $t_*\in(0,T)$ and $(s_n)\downarrow 0$, $(\sigma_n)\downarrow 0$ with $s_n<\sigma_n$ for all $n$. Then $(\mu_n)$ converges weakly to $\mu=2\delta_{t_*}$, we have 
$$
\alpha_n(t):=\mu_n([0,t])\to\alpha(t):=\mu([0,t])
$$
for every $t\in[0,T]$, but $(\alpha_n)$, or any subsequence of $(\alpha_n)$, does not converge to $\alpha$ in $D([0,T])$. This example should be compared to \cite[Example~1.20, p.329]{JacodShiryaev03}. See also Theorem~2.15, p.342 in \cite{JacodShiryaev03}.\medskip

As asserted above, we will show that the sequence of processes
$$
t\mapsto \<f^n(t),\varphi\>+A^n_\varphi(t),\quad A^n_\varphi(t):=\<m^n,\partial_\xi\varphi\>([0,t]),
$$
where
$$
\<m^n,\partial_\xi\varphi\>([0,t]):=\iiint_{\T^N\times[0,t]\times\R}\partial_\xi\varphi(x,s,\xi)dm^n(x,s,\xi),
$$
is tight in $D([0,T])$. It is sufficient to show that
\begin{equation}\label{TheJn}
t\mapsto \<f^n(t),\varphi\>+B^n_\varphi(t),\quad B^n_\varphi(t):=\<m^n,\partial_\xi\varphi\>([0,t])-\eps^n(t,\varphi)
\end{equation}
is tight in $D([0,T])$ since each function $t\mapsto\eps^n(t,\varphi)$ converges in probability to $0$ in $C([0,T])$ by \eqref{epsto0}. Since $f^n(0)=f^n_0$, we have 
\begin{equation}\label{equationfn}
\<f^n(t),\varphi\>+B^n_\varphi(t)=\<f^n_0,\varphi\>+J^n_\varphi(t),
\end{equation}
$\P$-almost-surely, where
\begin{multline}
J^n_\varphi\colon t\mapsto\int_0^t \<f^n(s),a(\xi)\cdot\nabla\varphi\>ds+\sum_{k\geq 1}\int_0^t\int_{\T^N}\int_\R g_k(x,\xi)\varphi(x,\xi)d\nu^n_{x,s}(\xi)dxd\beta_k(s)\\
+\frac{1}{2}\int_0^t\int_{\T^N}\int_\R \partial_\xi\varphi(x,\xi)\GG^2(x,\xi)d\nu^n_{x,s}(\xi) dx ds.
\label{def:Jphin}\end{multline}
We will show that $(J^n_\varphi(t))$ is tight in $C([0,T])$.

\begin{proposition} Let $(f^n)$ be a sequence of approximate generalized solutions to~\eqref{stoSCL} with initial datum $f_0^n$. For $\varphi\in C^\order_c(\T^N\times\R)$, set
\begin{align*}
D^n_\varphi(t)&=\int_0^t \<f^n(s),a(\xi)\cdot\nabla\varphi\>ds,\\
M^n_\varphi(t)&=\sum_{k\geq 1}\int_0^t\int_{\T^N}\int_\R g_k(x,\xi)\varphi(x,\xi)d\nu^n_{x,s}(\xi)dxd\beta_k(s),\\
I^n_\varphi(t)&=\frac{1}{2}\int_0^t\int_{\T^N}\int_\R \partial_\xi\varphi(x,\xi)\GG^2(x,\xi)d\nu^n_{x,s}(\xi) dx ds.
\end{align*}
Assume that \eqref{eq:integrabilityfn} is satisfied. Then each sequence $(D^n_\varphi)$, $(M^n_\varphi)$, $(I^n_\varphi)$ is tight in $C([0,T])$. In particular, the sequence $(J^n_\varphi)$ defined by \eqref{def:Jphin} is tight in $C([0,T])$.
\label{prop:Ctight}\end{proposition}

\textbf{Proof of Proposition~\ref{prop:Ctight}.} Note first the trivial uniform bounds
$$
\E|D^n_\varphi(t)|,\ \E|M^n_\varphi(t)|,\ \E|I^n_\varphi(t)| =\mathcal{O}(1),
$$
obtained for $t=0$ since all three terms vanish. We then use the Kolmogorov's criterion to obtain some bounds in some H\"older space $C^\alpha([0,T])$. We have the following estimate on the square of the increments of $D^n_\varphi$:
\begin{equation}\label{KolmogorovD}
\E|D^n_\varphi(t)-D^n_\varphi(\sigma)|^2\leq \|a\cdot\nabla\varphi\|_{L^1(\T^N\times\R)}^2|t-\sigma|^2,
\end{equation}
since $|f_n|\leq 1$ almost-surely. Similarly, using \eqref{D0} and \eqref{eq:integrabilityfn}, we have
\begin{equation}\label{KolmogorovI}
\E|I^n_\varphi(t)-I^n_\varphi(\sigma)|^2\leq D_0^2 T (1+C_4) \|\partial_\xi\varphi\|_{L^\infty(\T^N\times\R)}^2|t-\sigma|^2.
\end{equation}
The estimates \eqref{KolmogorovD} and \eqref{KolmogorovI} give some bounds on $\E\|D^n_\varphi\|_{C^\alpha([0,T])}$ and $\E\|I^n_\varphi\|_{C^\alpha([0,T])}$ respectively, for $\alpha<\frac12$. Furthermore, the Burkholder - Davis - Gundy  Inequality gives, for $p>2$,
\begin{align*}
\E|M^n_\varphi(t)-M^n_\varphi(\sigma)|^p&\leq \E\left[\sup_{\sigma\leq r\leq t}|M^n_\varphi(r)-M^n_\varphi(\sigma)|\right]^p\\
&\leq C_\mathrm{BDG}(p) \E\left[\sum_{k\geq 1}\int_\sigma^t\left|\int_{\T^N}\int_\R g_k(x,\xi)\varphi(x,\xi)d\nu^n_{x,s}(\xi)dx\right|^2ds\right]^{p/2}.
\end{align*}
By Jensen's Inequality, and a bound analogous to \eqref{sumhk}, we obtain
\begin{equation}\label{IncrementsMn}
\E|M^n_\varphi(t)-M^n_\varphi(\sigma)|^p\leq C_\mathrm{BDG}(p) \left[D_0(1+C_2)\right]^{p/2} \|\varphi\|_{L^\infty(\T^N\times\R)}^p |t-\sigma|^{p/2},
\end{equation}
and \eqref{IncrementsMn} gives a bound on $\E\|M^n_\varphi\|_{C^\alpha([0,T])}$ for $\alpha<\frac12-\frac1p$. We obtain in this way some tightness conditions on the laws of $D^n$, $M^n$, $I^n$ respectively on $C([0,T])$. \qed

\subsection{Convergence of approximate generalized solutions}\label{sec:MnM}

We conclude here this section about the stability of generalized solutions by the following statement.

\begin{theorem}[Convergence to martingale solutions] Let $(f^n)$ be a sequence of approximate generalized solutions to~\eqref{stoSCL} with initial datum $f_0^n$, satisfying \eqref{eq:integrabilityfn}, \eqref{Boundmn} and \eqref{inftymn}. We suppose that there exists a kinetic function $f_0$ on $\T^N\times\R$ such that $f^n_0\to f_0$ in $L^\infty(\T^N\times\R)$-weak-*. Then there exists a probability space $(\tilde\Omega,\tilde{\mathcal{F}},\tilde\P)$, a filtration $\tilde{\mathcal{F}}_t$, some $\tilde{\mathcal{F}}_t$-adapted independent Brownian motions $(\tilde{\beta}_k)_{k\geq 1}$, some random Young measures $\tilde\nu^n$, $\tilde\nu$, some random measures $\tilde m^n$, $\tilde m$ on $\T^N\times[0,T]\times\R$ such that
\begin{enumerate} 
\item\label{item:1thmartingale} $\tilde\nu^n$ has the same law as $\nu^n$,
\item\label{item:3thmartingale} up to a subsequence still denoted $(\tilde\nu^n)$, there is $\tilde\P$-almost-sure convergence of $(\tilde\nu^n)$ to $\tilde\nu$ in $\mathcal{Y}^1$.
\item\label{item:Pnu} for all $\psi\in C_b(\R)$, the random map $(x,t)\mapsto\<\psi,\tilde\nu_{x,t}\>$ belongs to $L^2_{\tilde{\mathcal{P}}}(\T^N\times [0,T]\times\tilde\Omega)$,
\item\label{item:5thmartingale} $\tilde m^n$ has the same law as $m^n$,
\item\label{item:7thmartingale} up to a subsequence still denoted $(\tilde m^n)$, there is $\tilde\P$-almost-sure convergence of $(\tilde m^n)$ to $\tilde m$ in $\mathcal{M}_b(\T^N\times[0,T]\times\R)$-weak-*.
\end{enumerate}
Let $\tilde f$ be defined by $\tilde f(x,t,\xi)=\tilde\nu_{x,t}(\xi,+\infty)$, then, $\tilde\P$-almost-surely, $\tilde f$ is a kinetic function and
\begin{enumerate}[resume]
\item\label{item:8thmartingale} up to a subsequence, and $\tilde\P$-almost-surely, $\tilde f^n$  converges in $L^\infty(\T^N\times[0,T]\times\R)$-weak-* to $\tilde f$
\item\label{item:tildefcadlag} $\tilde\P$-almost-surely, for all $\varphi$ in $C_c(\T^N\times\R)$, $t\mapsto\<\tilde{f}(t),\varphi\>$ is c{\`a}dl{\`a}g,
\item\label{item:intnu} $\tilde\nu$ satisfies
\begin{equation}
\tilde\E\left(\sup_{t\in[0,T]}\int_{\T^N}\int_\R|\xi|^p d\tilde{\nu}_{x,t}(\xi) dx\right) \leq C_p,
\label{eq:integrabilityftildeOK}\end{equation}
for all $1\leq p<+\infty$, where $C_p$ is a finite constant,
\item\label{item:soltildef} for all $\varphi\in C^1_c(\T^N\times\R)$, for all $t\in[0,T]$, $\tilde{\P}$-almost-surely, $\tilde{f}$ satisfies
\begin{align}
\<\tilde{f}(t),\varphi\>=&\<f_0,\varphi\>+\int_0^t \<\tilde{f}(s),a(\xi)\cdot\nabla_x\varphi\>ds
\nonumber\\
&+\int_0^t\int_{\T^N}\int_\R g_k(x,\xi)\varphi(x,\xi)d\tilde{\nu}_{x,s}(\xi) dxd\tilde{\beta}_k(s) \nonumber\\
&+\frac{1}{2}\int_0^t\int_{\T^N}\int_\R \GG^2(x,\xi)\partial_\xi\varphi(x,\xi)d\tilde{\nu}_{x,s}(\xi) dx ds
-\tilde{m}(\partial_\xi\varphi)([0,t]).
\label{eq:kineticfpretilde}
\end{align}
\end{enumerate}
\label{th:martingalesol} \end{theorem}

After one does the substitution
$$
(\Omega,\mathcal{F},\P,\mathcal{F}_t,\beta_k(t))\leftarrow (\tilde{\Omega},\tilde{\mathcal{F}},\tilde{\P},\tilde{\mathcal{F}}_t,\tilde{\beta}_k(t)),
$$
which is a substitution of the probabilistic data in the Cauchy Problem for Equation~\eqref{stoSCL}, the points \ref{item:Pnu}, \ref{item:tildefcadlag}, \ref{item:intnu}, \ref{item:soltildef} in Theorem~\ref{th:martingalesol} show that $\tilde f$ is a generalized solution associated to the initial datum $f_0$. Such a function $\tilde f$, which turns out to be a generalized solution to \eqref{stoSCL} after a substitution of the probabilistic data, is called a martingale generalized solution. The term martingale refers to the martingale characterization of \eqref{eq:kineticfpretilde}, \textit{cf.} Proposition~\ref{prop:3martingales} and Proposition~\ref{prop:3martingalesCV}, which we will use to prove Theorem~\ref{th:martingalesol}.\medskip

\subsection{Proof of Theorem~\ref{th:martingalesol}}\label{sec:proofThmartingale}

In this section, we will give the proof of Theorem~\ref{th:martingalesol}. We will use the results (and the proofs) of Proposition~\ref{prop:stabextractnu}, Proposition~\ref{prop:stabextractm}, see Section~\ref{sec:compactnun} and Section~\ref{sec:compactmn} respectively. 

\subsubsection{State space and Skorokhod's Theorem}\label{sec:proofThmartingalestatespace}
Recall that
$$
W(t)=\sum_{k\geq 1}\beta_k(t) e_k,
$$
where $(e_k)_{k\geq 1}$ is the orthonormal basis of the Hilbert space $H$. Let $\mathfrak{U}$ be an other separable Hilbert space such that $H\hookrightarrow\mathfrak{U}$ with Hilbert-Schmidt injection. Then the trajectories of $W$ are $\P$-a.s. in the path-space $\mathcal{X}_W=C([0,T];\mathfrak{U})$ (see \cite[Theorem~4.3]{DaPratoZabczyk92}). We consider the $C^\order$-norm
$$
\|\varphi\|_{C^\order}=\sup\{\| D^m \varphi\|_{L^\infty(\T^N\times\R)}; m\in\{0,\ldots,\order\}^{N+1}\}
$$
on $C^\order_c(\T^N\times\R)$. Let 
$$
\Gamma=\{\varphi_1,\varphi_2,\ldots\}
$$ 
be a dense countable subset of $C^\order_c(\T^N\times\R)$ for this norm. We can construct $\Gamma$ as follows: let 
$$
\rho_\eps(x,\xi):=\frac{1}{\eps^{N+1}}\rho(\eps^{-N}x,\eps^{-1}\xi)
$$ 
be a compactly supported approximation of the unit on $\T^N\times\R$. Let $\{\theta_p;p\in\N\}$ be a dense subset of $L^1(\T^N\times\R)$. We can assume that all the functions $\theta_p$ are compactly supported (otherwise, we use a process of truncation). Then any function in $C^\order_c(\T^N\times\R)$ can be approximated by functions in 
$$
\Gamma:=\{\rho_{k^{-1}}*\theta_p;p\in\N, k\in\N^*\}\subset C^\order_c(\T^N\times\R)
$$
for the convergence measured by the $C^\order$-norm. Indeed, given $\varphi\in C^\order_c(\T^N\times\R)$, $a>0$, and $m\in\{0,\ldots,\order\}^{N+1}$, we have, by the triangular inequality,
\begin{align}
\|D^m\varphi-D^m\rho_{\eps}*\theta_p\|_{L^\infty}
&\leq \|D^m\varphi-D^m\rho_{\eps}*\varphi\|_{L^\infty}+\|D^m\rho_{\eps}*(\varphi-\theta_p)\|_{L^\infty}\nonumber\\
&\leq\omega_{D^m\varphi}(\eps)+\frac{\|\rho\|_{L^\infty}}{\eps^{N+1+|m|}}\|\varphi-\theta_p\|_{L^1},\label{GammaDense1}
\end{align}
since the norm of $D^m\rho_\eps$ in $L^\infty$ is bounded by $\frac{\|\rho\|_{L^\infty}}{\eps^{N+1+|m|}}$. In \eqref{GammaDense1}, $\omega_{D^m\varphi}$ denotes the modulus of continuity of $D^m\varphi$. We choose $\eps=k^{-1}$ with $k$ large enough to ensure $\omega_{D^m\varphi}(\eps)<a$ for all  $m\in\{0,\ldots,\order\}^{N+1}$. Taking then $p\in\N$ such that $\|\varphi-\theta_p\|_{L^1}<a\eps^{(\order+1)(N+1)}$, we obtain $\|\varphi-\rho_{k^{-1}}*\theta_p\|_{C^\order}<2a$.\medskip

Let also $\R^\infty$ denote the product space $\prod_{\varphi\in\Gamma}\R$ endowed with the topology of point-wise convergence. As such, $\R^\infty$ is separable, complete and admits a compatible metric. Define the Polish space
$$
\mathcal{E}:=C([0,T];\R^\infty)\times C([0,T];\R^\infty)\times C([0,T];\R^\infty)\times\R^\infty,
$$
and
\begin{equation}\label{deffleftarrow}
\eps^n_\varphi(t)=\eps^n(t,\varphi),\quad J^n_\varphi(t)=\<f^n(t),\varphi\>-\<f^n_0,\varphi\>+\<m^n,\partial_\xi\varphi\>([0,t])-\eps^n_\varphi(t),
\end{equation}
for all $\varphi\in C^\order_c(\T^N\times\R)$.
Note that, as a consequence of Equation~\eqref{equationfn} and Proposition~\ref{prop:Ctight}, we know that, for all $\varphi\in C^\order_c(\T^N\times\R)$,
\begin{equation}\label{fntightD}
(J^n_\varphi)\mbox{ is tight in }C([0,T]).
\end{equation}
By \eqref{epsto0}, we also have $\eps^n_\varphi\to 0$ in $C([0,T])$ in probability, for all $\varphi\in  C^\order_c(\T^N\times\R)$. 
We introduce the four following sequences:
$$
\{J^n(t)\}:=(J^n_{\varphi}(t))_{\varphi\in\Gamma}, \quad \{M^n(t)\}:=(M^n_{\varphi}(t))_{\varphi\in\Gamma},\quad \{\eps^n(t)\}:=(\eps^n_{\varphi}(t))_{\varphi\in\Gamma},
$$
and $\{f^n_\mathrm{in}\}:=(\<f^n_0,\varphi\>)_{\varphi\in\Gamma}$, where $M^n_\varphi$ is defined by \eqref{def:Mnphi}. We will consider the multiplet
$$
Z^n=(\nu^n,\{J^n\},\{M^n\},\{\eps^n\},\{f^n_\mathrm{in}\},\mu^n,\|m^n\|,W)\in\mathcal{X},
$$
where the state space $\mathcal{X}$ is 
$$
\mathcal{X}:=\mathcal{Y}^1\times\mathcal{E}\times\mathcal{P}^1(\T^N\times[0,T]\times\R)\times\R_+\times\mathcal{X}_W.
$$

Let $\eps>0$. By \eqref{fntightD}, there exists for each $j\in\N$ a compact $K_j$ in $C([0,T])$ such that
$$
\inf_{n\in\N}\P\left(J^n_{\varphi_j}\in K_j\right)\geq 1-\frac{\eps}{2^j}.
$$
Let $K=\prod_{j\in\N}K_j$. Then $K$ is compact\footnote{since $C([0,T];\R^\infty)$ is homeomorphic to the countable product, over $\Gamma$, of copies of $C([0,T];\R)$} in $C([0,T];\R^\infty)$ and 
$$
\P(\{J^n\}\in K^c)\leq \sum_{j\in\N} \P\left(J^n_{\varphi_j}\in K_j^c\right)\leq \sum_{j\in\N}\frac{\eps}{2^j}=2\eps,
$$
for all $n\in\N$. This shows that $(\{J^n\})$ is tight in $C([0,T];\R^\infty)$. We have similar results about $(\{M^n\})$ and $(\{\eps^n\})$ by Proposition~\ref{prop:Ctight}. On 
$\mathcal{X}_W$ we consider the topology induced by the norm
$$
\|v\|=\sup_{t\in[0,T]}\|v(t)\|_{\mathfrak{U}}
$$
Then $\mathcal{X}_W$ is separable and complete. A first consequence of this is the fact that the law of the single random variable $W$ is tight in $\mathcal{X}_W$. A second consequence is the fact that $\mathcal{X}$ is a separable completely-metrizable space. By Section~\ref{sec:compactnun} and Section~\ref{sec:compactmn}, we conclude that $(Z^n)$ is tight in the Polish space $\mathcal{X}$. We may thus apply Skorokhod's Theorem to $(Z^n)$: there exists a probability space $ (\tilde{\Omega},\tilde{\mathcal{F}},\tilde{\P})$ and some random variable $\tilde Z^n$, $\tilde Z$ such that $\tilde Z^n$ has the same law as $Z^n$ and, up to a subsequence, $\tilde\P$-almost-surely, $\tilde Z^n$ converges to $\tilde Z$ in $\mathcal{X}$. 


\subsubsection{Identification of the limit: c{\`a}dl{\`a}g version}\label{sec:proofThmartingalecadlag}

Let us denote the component of $\tilde{Z}$ as follows:
$$
\tilde Z=(\tilde\nu,\{\tilde{J}\},\{\tilde{M}\},\{\tilde{\eps}\},\{\tilde{f}_\mathrm{in}\},\tilde\mu,\tilde\alpha,\tilde{W}).
$$
Note first that $\tilde{\eps}=0$ by \eqref{epsto0}. We have also
$$
\{\tilde{f}_\mathrm{in}\}=(\<f_0,\varphi\>)_{\varphi\in\Gamma}
$$
since $f^n_0\to f_0$ in $L^\infty(\T^N\times\R)$-weak-* by hypothesis. Recall (see Proposition~\ref{prop:stabextractnu} and Proposition~\ref{prop:stabextractm}) that $\tilde{f}(x,t,\xi)=\tilde{\nu}_{x,t}(\xi,+\infty)$ and $\tilde{m}=\tilde{\alpha}\tilde{\mu}$.
It was shown in the proof of Proposition~\ref{prop:stabextractnu} that item \ref{item:1thmartingale}, \ref{item:3thmartingale}, \ref{item:8thmartingale} of Theorem~\ref{th:martingalesol} are satisfied and that the moments of $\tilde{\nu}$ are bounded as in \eqref{eq:integrabilityftilde}. By the proof of Proposition~\ref{prop:stabextractm}, we have also \ref{item:5thmartingale}, \ref{item:7thmartingale} of Theorem~\ref{th:martingalesol}. We will first establish the following result.

\begin{lemma}\label{lem:negN0} We have the following identities: $\tilde{\P}$-almost-surely, 
\begin{equation}\label{negN0n}
\mbox{for all }t\in[0,T],\mbox{ for all }\varphi\in\Gamma,\; \<\tilde{f}^n(t),\varphi\>=\<f_0,\varphi\>+\tilde{J}^n_\varphi(t)-\<\tilde{m}^n,\partial_\xi\varphi\>([0,t])-\tilde{\eps}^n_\varphi(t),
\end{equation}
and: $\tilde{\P}$-almost-surely, there exists a negligible set $N_0\subset[0,T]$ such that,
\begin{equation}\label{negN0}
\mbox{for all }t\in[0,T]\setminus N_0,\mbox{ for all }\varphi\in\Gamma,\; \<\tilde{f}(t),\varphi\>=\<f_0,\varphi\>+\tilde{J}_\varphi(t)-\<\tilde{m},\partial_\xi\varphi\>([0,t]).
\end{equation}
\end{lemma}

We will use \eqref{negN0} to prove Proposition~\ref{prop:LRtilde}, where we obtain a c{\`a}dl{\`a}g version of $\tilde{f}$ (c{\`a}dl{\`a}g in the sense that $\tilde{\P}$-almost-surely, for all $\varphi\in C_c(\T^N\times\R)$, $t\mapsto\<\tilde{f}(t),\varphi\>$ is c{\`a}dl{\`a}g).\medskip

\textbf{Proof of Lemma~\ref{lem:negN0}.} Let $\theta\in C([0,T])$. Let us integrate the identity \eqref{deffleftarrow} against $\theta$. Using the Fubini theorem, we obtain: $\P$-almost-surely,
\begin{equation}\label{negN01}
\int_0^T (J^n_\varphi(t)+\<f^n_0,\varphi\>+\eps^n_\varphi(t))\theta(t)dt-\<\nu^n,\Phi\>-\<m^n,\Psi\>=0,
\end{equation}
where
$$
\Psi(x,t,\xi)=\partial_\xi\varphi(x,\xi)\int_t^T\theta(s)ds,\quad \Phi(x,t,\xi)=\int_{-\infty}^\xi\varphi(x,\zeta)d\zeta\theta(t).
$$
Note that $\Psi$ and $\Phi$ are continuous and bounded functions. Taking the square, then expectancy in \eqref{negN01} gives $\E F(Z^n)=0$, where $F\colon\mathcal{X}\to\R$ defined by
$$
F(Z^n)=\left|\int_0^T (J^n_\varphi(t)+\<f^n_0,\varphi\>+\eps^n_\varphi(t))\theta(t)dt-\<\nu^n,\Phi\>-\|m^n\|\<\mu^n,\Psi\>\right|^2
$$
is a continuous function. By identity of the laws of $Z^n$ and $\tilde{Z}^n$, we have $\tilde{\E}F(\tilde{Z}^n)=0$ for all $n$. Since $F$ is non-negative, this means $F(\tilde{Z}^n)=0$,  $\tilde{\P}$-almost-surely. Since $\theta$ is arbitrary and $\Gamma$ is countable, we deduce \eqref{negN0n}, \textit{a priori} for $t\in [0,T]\setminus N_n$, where $N_n$ is a measurable negligible set. We can take $N_n=\emptyset$ because  both sides of \eqref{negN0n} are c{\`a}dl{\`a}g functions. By almost-sure convergence, that $F(\tilde{Z}^n)=0$,  $\tilde{\P}$-almost-surely implies $F(\tilde{Z})=0$,  $\tilde{\P}$-almost-surely. Similarly, hence, we obtain \eqref{negN0}. \qed\medskip
\begin{proposition}\label{prop:LRtilde} There exists a measurable subset $\tilde\Omega^+$ of $\tilde\Omega$ of probability one, a random Young measure $\tilde{\nu}^+$ on $\T^N\times(0,T)$ such that
\begin{enumerate}
\item\label{item:nunuplus} for all $\tilde\omega\in\tilde\Omega^+$, for a.e. $(x,t)\in\T^N\times(0,T)$, the probability measures $\tilde{\nu}^+_{x,t}$ and $\tilde{\nu}_{x,t}$ coincide,
\item\label{item:cadlagplusOK} the kinetic function $\tilde{f}^+(x,t,\xi):=\tilde{\nu}^+_{x,t}(\xi,+\infty)$ satisfies: for all $\tilde\omega\in\tilde\Omega^+$, for all $\varphi\in C_c(\T^N\times(0,T))$, $t\mapsto\<\tilde{f}^+(t),\varphi\>$ is c{\`a}dl{\`a}g,
\item\label{item:nuplusOKbound} the random Young measure $\tilde{\nu}^+$ satisfies \eqref{eq:integrabilityftildeOK}.
\end{enumerate} 
\end{proposition}

\textbf{Proof of Proposition~\ref{prop:LRtilde}.} The proof is quite similar to the proof of Proposition~\ref{prop:LRlimits}. For $\varphi\in\Gamma$, let $F_\varphi(t)$ denote the right-hand side of \eqref{negN0}:
\begin{equation}\label{def:Fvarphi}
F_\varphi(t)=\<f_0,\varphi\>+\tilde{J}_\varphi(t)-\<\tilde{m},\partial_\xi\varphi\>([0,t]).
\end{equation}
We define $\tilde\Omega^+$ as the intersection of the three following events: first \eqref{negN0}, second: ``for all $\varphi\in\Gamma$, $F_\varphi$ is c{\`a}dl{\`a}g", third the event
$$
\sup_{J\subset[0,T]}\frac{1}{|J|}\int_J\int_{\T^N}\int_\R|\xi|^p d\tilde\nu_{x,t}(\xi) dx dt<+\infty,
$$ 
where the supremum over intervals $J$ is as in \eqref{eq:integrabilityftilde} (a countable supremum over all open intervals $J\subset[0,T]$ with rational extremities). Assume that $\tilde\Omega^+$ is realized (say we draw a particular $\tilde{\omega}\in\tilde{\Omega}^+$). Assume in particular that
\begin{equation}\label{Cppp}
\frac{1}{|J|}\int_J\int_{\T^N}\int_\R|\xi|^p d\tilde\nu_{x,t}(\xi) dx dt\leq C_p(\tilde\omega),
\end{equation}
for all open intervals $J\subset[0,T]$ with rational extremities. Then the map
$$
t\mapsto \int_{\T^N}\int_\R|\xi|^p d\tilde\nu_{x,t}(\xi) dx
$$
is integrable on $(0,T)$. A simple approximation procedure shows then that \eqref{Cppp} holds true when $J$ is any interval in $[0,T]$.\smallskip

Let $t_*\in[0,T)$. Let $(\eps_l)$ be a sequence of positive numbers decreasing to $0$ such that $t_*+\eps_1<T$. Let $J_l=(t_*,t_*+\eps_l)$. Consider the sequence of Young measures, and corresponding kinetic functions
$$
\tilde\nu^{(l)}_{x}=\frac{1}{|J_l|}\int_{J_l}\tilde\nu_{x,t} dt,\quad \tilde{f}^{(l)}(x,\xi)=\tilde\nu^{(l)}_{x}(\xi,+\infty)=\frac{1}{|J_l|}\int_{J_l}\tilde{f}(x,t,\xi )dt.
$$
Since the Borel $\sigma$-algebra of $\T^N$ is countably generated ($\T^N$ being separable),  we can apply Theorem~\ref{th:youngmeasure} and Corollary~\ref{cor:kineticfunctions}. There exists a subsequence $(l_m)$ and a Young measure $\tilde\nu^*$ such that $\tilde\nu^{(l_m)}\to \tilde\nu^*$ in the sense of \eqref{cvYoungMeasure} and $\tilde{f}^{(l_m)}\to\tilde{f}^*$ in $L^\infty(\T^N\times\R)$ weak-*, where $\tilde{f}^*(x,\xi)=\tilde\nu^*_x(\xi,+\infty)$. The limit $\tilde{f}^*$ is unique. Indeed, if $\varphi\in\Gamma$, then, due to \eqref{negN0} and to the Fubini theorem, and due to the right-continuity of $F_\varphi$, we have
 $$
\<\tilde{f}^{(l)},\varphi\>=\frac{1}{|J_l|}\int_{J_l} F_\varphi(t) dt\to F_\varphi(t_*).
$$
This implies 
\begin{equation}\label{starN0ornot}
\<\tilde{f}^*,\varphi\>=F_\varphi(t_*).
\end{equation}
Since $\Gamma$ is dense, $\tilde{f}^*$ and $\tilde\nu^*=-\partial_\xi\tilde{f}^*$ are unique. We deduce that the convergence holds along the whole sequence $l=1,2.\ldots$, independently on $\tilde{\omega}\in\tilde{\Omega}^+$ and on $t_*\in[0,T)$. Consequently, setting
$$
\tilde{\nu}^+_{x,t_*}=\tilde\nu^*_x,\quad \tilde{f}^+(x,t_*,\xi)=\tilde{\nu}^+_{x,t_*}(\xi,+\infty),
$$
we have: for all $\tilde{\omega}\in\tilde{\Omega}^+$, for all $t\in[0,T)$, for all $\phi\in C_b(\T^N\times\R)$,
\begin{equation}\label{cvtildenuplus}
\frac{1}{\eps}\int_t^{t+\eps}\iint_{\T^N\times\R}\phi(x,\xi)d\tilde{\nu}_{x,s}(\xi) dx ds\to\iint_{\T^N\times\R}\phi(x,\xi)d\tilde{\nu}^+_{x,t}(\xi) dx.
\end{equation}
Since $T$ is arbitrary, we can as well work on $[0,T+1]$, instead of $[0,T]$. In that way, we can give a meaning to $\tilde{\nu}^+_{x,t}$ for $t=T$ also. By \eqref{negN0} and \eqref{starN0ornot}, we have $\tilde{f}^+(x,t,\xi)=\tilde{f}(x,t,\xi)$ and $\tilde{\nu}^+_{x,t}=\tilde{\nu}_{x,t}$ for all 
$\tilde{\omega}\in\tilde{\Omega}^+$, for all $t\in (0,T)\setminus N_0$, for a.e. $(x,\xi)\in\T^N\times\R$. If $\phi\in C_b(\R)$ and $\tilde{\omega}\in\tilde{\Omega}^+$, then $(x,t)\mapsto\<\tilde\nu_{x,t},\phi\>$ is measurable and $(x,t)\mapsto\<\tilde\nu^+_{x,t},\phi\>$ differs from the latter function on a negligible subset of $\T^N\times(0,T)$. Therefore $(x,t)\mapsto\<\tilde\nu^+_{x,t},\phi\>$ itself is measurable. We deduce that $\tilde{\nu}^+$ and $\tilde{f}^+$ satisfy the measurability properties of  a random Young measure and a random kinetic function respectively, and point~\ref{item:nunuplus} of the proposition is proved. The point~\ref{item:cadlagplusOK} of the proposition follows from \eqref{starN0ornot}, which gives $\<\tilde{f}(t),\varphi\>=F_\varphi(t)$ for all $t$.
To obtain the last point~\ref{item:nuplusOKbound} of the proposition. We note first that $\tilde{\nu}^+$, like $\tilde\nu$, satisfies \eqref{eq:integrabilityftilde}. If 
$$
\sup_{J\subset[0,T]}\frac{1}{|J|}\int_J\int_{\T^N}\int_\R|\xi|^p d\tilde\nu^+_{x,t}(\xi) dx dt<+\infty,
$$
which happen $\tilde{\P}$-almost-surely, then 
$$
\sup_{J\subset[0,T]}\frac{1}{|J|}\int_J\int_{\T^N}\int_\R|\xi|^p d\tilde\nu^+_{x,t}(\xi) dx dt=\sup_{t\in[0,T]}\int_{\T^N}\int_\R|\xi|^p d\tilde\nu^+_{x,t}(\xi) dx
$$
by right-continuity of $t\mapsto\tilde{\nu}^+_t$. This gives the desired result. \qed\medskip

We will now consider only the c{\`a}dl{\`a}g versions: we replace $\tilde{\nu}$ by $\tilde{\nu}^+$ and $\tilde{f}$ by $\tilde{f}^+$. This amounts to a modification on a negligible set. Therefore, this does not affect the results \ref{item:1thmartingale}, \ref{item:3thmartingale}, \ref{item:5thmartingale}, \ref{item:7thmartingale}, \ref{item:8thmartingale} in Theorem~\ref{th:martingalesol}. We have now also items~\ref{item:tildefcadlag} and \ref{item:intnu} of the theorem. There remains to define the filtration $(\tilde{\F}_t)$, the Wiener processes $\tilde{\beta}_k$ and to prove the points~\ref{item:Pnu} and \ref{item:soltildef} of the theorem. We define $(\tilde{\F}_t)$, $\tilde{\beta}_k$ and show item~\ref{item:Pnu} in the proof of convergence of the stochastic integral in the next section~\ref{sec:proofThmartingalecvsto}. The equation~\eqref{eq:kineticfpretilde} is established in Section~\ref{sec:proofThmartingaleEquation}. To finish the current section, let us first record the fact that \eqref{negN0} is now true for all $t$, due to our re-definition of $\tilde{f}$ and to \eqref{starN0ornot}: $\tilde{\P}$-almost-surely,
\begin{equation}\label{negN0tilde}
\mbox{for all }t\in[0,T),\mbox{ for all }\varphi\in\Gamma,\; \<\tilde{f}(t),\varphi\>=\<f_0,\varphi\>+\tilde{J}_\varphi(t)-\<\tilde{m},\partial_\xi\varphi\>([0,t]).
\end{equation}
We deduce from \eqref{negN0tilde} the following lemma.

%
%

\begin{proposition}\label{prop:CVfftilde} There exists a countable subset $\tilde{B}\subset[0,T]$ such that, $\tilde{\P}$-almost-surely, for all $t\in[0,T]\setminus\tilde{B}$, for all $\varphi\in C_c(\T^N\times\R)$, $\<\tilde{f}^n(t),\varphi\>\to\<\tilde{f}(t),\varphi\>$.
\end{proposition}

\textbf{Proof of Proposition~\ref{prop:CVfftilde}.} It is sufficient to obtain the convergence for $\varphi\in\Gamma$. We apply Lemma~\ref{lem:atomicpoints}. Let  
\begin{equation}\label{defJstartilde}
\tilde{B}=\left\{t\in[0,T];\tilde{\P}\left(\pi_\#\tilde{m}(\{t\})>0\right)>0\right\}.
\end{equation}
Then $\tilde{B}$ is countable. Since $|\<\tilde{m},\partial_\xi\varphi\>(\{t\})|\leq\|\partial_\xi\varphi\|_{L^\infty}\pi_\#\tilde{m}(\{t\})$,
$\tilde\P$-almost-surely, we have $\<\tilde{m},\partial_\xi\varphi\>(\{t\})=0$ for all $t\in[0,T]\setminus\tilde{B}$. For $t\in[0,T]\setminus\tilde{B}$ then, the right-hand side of \eqref{negN0n} is converging to the right-hand side of \eqref{negN0tilde}. We deduce the convergence of the left-hand sides, \textit{i.e.} $\<\tilde{f}^n(t),\varphi\>\to\<\tilde{f}(t),\varphi\>$.
%
\qed\medskip

%
%
%

\subsubsection{Identification of the limit: convergence of the stochastic integral}\label{sec:proofThmartingalecvsto}

Let us set
\begin{align}
M^*_\varphi(t)&=\sum_{k\geq 1}\int_0^t\int_{\T^N}\int_\R g_k(x,\xi)\varphi(x,\xi)d\tilde{\nu}_{x,s}(\xi)dxd\tilde{\beta}_k(s),\label{def:tildeMphi}
\end{align}
($\tilde{\beta}_k$ is defined in Lemma~\ref{lem:tildeW} below). Our aim is to prove the identification $\{\tilde{M}\}=\{M^*\}$. To obtain this result, we will use the martingale characterization developed in Section~\ref{sec:martingale}. The proof is decomposed in several steps.\medskip

\paragraph{Step 1. Filtration}
The approximation procedures to \eqref{stoSCL} (vanishing viscosity me\-thod, Finite Volume method as here) construct approximate solutions on arbitrary time intervals $[0,T]$. We will therefore consider the functions as defined on the whole time interval $\R_+$. This is simply to avoid the special case of the final time in the definition of the Skorokhod space $D([0,T])$, \textit{cf.} \cite{BillingsleyBook}, \cite[Remark~1.10, p.~326]{JacodShiryaev03}. Let $E$ be a Polish space. Let us introduce the following notations (see \cite[Definition~1.1 p.~325]{JacodShiryaev03} in the case $E=\R^m$): on the space $D(\R_+;E)$, $\mathscr{D}^0_t(E)$ is the $\sigma$-algebra generated by the maps $\alpha\mapsto\alpha(s)$, $s\leq t$; 
$$
\mathscr{D}_t(E)=\bigcap_{t<s}\mathscr{D}^0_s(E),\quad\mathscr{D}_{t-}(E)=\bigvee_{s<t}\mathscr{D}_s(E).
$$
Note that $\mathscr{D}_t(E)\not=\mathscr{D}^0_t(E)$: the time of entrance in an open subset $U$ of $E$,
$$
\tau_U(\alpha)=\inf\left\{t\geq 0;\alpha(t)\in U\right\}
$$
is a stopping time with respect to $(\mathscr{D}_t(E))$, but not with respect to $(\mathscr{D}^0_t(E))$ \cite[Proposition~I.4.6]{RevuzYor99}.

\begin{proposition} Let $t>0$. Given a continuous bounded function $\theta\colon E\to\R$, $s\in[0,t)$ and $\eps>0$, let $\theta_{\# s}$ denote the evaluation map $\alpha\mapsto \theta(\alpha(s))$ on $D(\R_+;E)$, and let $\theta^\eps_{\# s}$ denote the regularization
\begin{equation}\label{fseps}
\theta_{\# s}^\eps\colon\alpha\to\frac{1}{\eps}\int_s^{t\wedge(s+\eps)}\theta(\alpha(\sigma))d\sigma
\end{equation}
of $\theta_{\# s}$. Then $\theta_{\# s}^\eps$ is a $\mathscr{D}_{t-}(E)$-measurable bounded function, continuous for the Skorokhod topology. Let $\mathcal{H}$ denote the set of functions
\begin{equation}\label{defHHHH}
H=\theta_{\# s_1}^{1,\eps_1}\cdots \theta_{\# s_k}^{k,\eps_k},
\end{equation}
where $k\geq 1$, $0\leq s_1<\cdots<s_k<t$, $0<\eps_1,\ldots,\eps_k$, $\theta^1,\ldots,\theta^k\in C_b(E)$. Then every characteristic function $\mathbf{1}_A$ of a cylindrical set $A\in\mathscr{D}_{t-}(E)$ of the form
\begin{equation}\label{defAAAA}
A=\left\{\alpha\in D(\R_+;E);\alpha(\tau_1)\in B_1,\ldots,\alpha(\tau_k)\in B_k\right\},
\end{equation}
for $B_1,\ldots,B_k$ closed subsets of $E$ and $0\leq \tau_1<\cdots<\tau_k<t$, is the bounded pointwise limit of a sequence of functions in $\mathcal{H}$.
\label{prop:Dtminus}\end{proposition}

\textbf{Proof of Proposition~\ref{prop:Dtminus}.} This is essentially the proof of \cite[Lemma~1.45 p.335]{JacodShiryaev03}. Let $\alpha\in D(\R_+;E)$ and let $(\alpha_n)$ be a sequence in $D(\R_+;E)$ such that $\alpha_n\to\alpha$ a.e. on $[0,t]$: this is the case if $\alpha_n\to\alpha$ in $D(\R_+;E)$ since $\alpha_n(\sigma)\to\alpha(\sigma)$ for every $\sigma$ not in the (countable) jump set of $\alpha$. Then, by the dominated convergence theorem, $\theta_{\# s}^\eps(\alpha_n)\to \theta_{\# s}^\eps(\alpha)$. Therefore $\theta_{\# s}^\eps$ is a bounded function, continuous for the Skorokhod topology. It is $\mathscr{D}_{t-}(E)$-measurable since it is the bounded pointwise limit when $\eta\to 0$ of the sequence of $\mathscr{D}_{t-}(E)$-measurable functions
$$
\alpha\to\frac{1}{\eps}\int_s^{(t-\eta)\wedge(s+\eps)}\theta(\alpha(\sigma))d\sigma.
$$
Let us prove the last point. We can choose some sequences of continuous bounded functions $\theta_1^n,\ldots,\theta_k^n\colon E\to\R$ converging simply to the characteristic functions $\mathbf{1}_{B_1},\ldots,\mathbf{1}_{B_k}$ (by considering, for example, the function distance to $B_j$, which is continuous). Since $\theta^\eps_{\# s}$ is approaching $\theta_{\# s}$ for the bounded pointwise convergence, the result follows. \qed\bigskip

\begin{remark}\label{sta.e.} Note that the function $H$ defined by \eqref{defHHHH} is more than merely continuous for the Skorokhod topology. Indeed, what we have seen in the proof of Proposition~\ref{prop:Dtminus} is that, for any $\alpha\in D(\R_+;E)$ and any sequence $(\alpha_n)$ in $D(\R_+;E)$ such that $\alpha_n\to\alpha$ a.e. on $[0,t]$, we have $H(\alpha_n)\to H(\alpha)$.
\end{remark}

Let us set 
$$
\{\tilde{f}\}=(\<\tilde{f},\varphi_j\>)_{j\in\N},\quad E=\R^\infty\times\R^\infty\times\mathfrak{U}.
$$
Recall that $\R^\infty$ is the product space $\prod_{\varphi\in\Gamma}\R$ endowed with the topology of point-wise convergence. Since $E$ is a product of Polish spaces, it is a Polish space. Since the product of $D(\R_+;\R^\infty)$ with $C(\R_+;\R^\infty\times\mathfrak{U})$ is, topologically, a subset of $D(\R_+;E)$, the triplet $(\{\tilde{f}\},\{\tilde{M}\},\tilde{W})$ is an element of $D(\R_+;E)$.

\begin{definition} The filtration $(\tilde{\mathcal{F}}_t)$ is the completion of the filtration generated by the triplet $(\{\tilde{f}\},\{\tilde{M}\},\tilde{W})$:
\begin{equation}\label{deftildeFt}
\tilde{\mathcal{F}}_t=\sigma\big((\{\tilde{f}\},\{\tilde{M}\},\tilde{W})^{-1}\left(\mathscr{D}_t(E)\right)\cup\big\{N\in\tilde{\mathcal{F}};\;\tilde{\P}(N)=0\big\}\big),\quad t\in[0,T].
\end{equation}
\end{definition}

Note that $(\tilde{\mathcal{F}}_t)$ is right-continuous since $(\mathscr{D}_t(E))$ is, and complete by definition.

\paragraph{Step 2. Wiener process}
Let $j\colon H\to\mathfrak{U}$ denote the injection of $H$ into $\mathfrak{U}$. Note that $j\circ j^*$ is a Trace-class operator on $\mathfrak{U}$. The Brownian motions $\tilde{\beta}_k(t)$ are the components of $\tilde{W}(t)$ on the orthonormal basis $(e_k)$:
\begin{lemma}\label{lem:tildeW}
The process $\tilde{W}$ has a modification which is a $(\tilde{\mathcal{F}}_t)$-adapted  $j\circ j^*$-Wiener process, and there exists a collection of mutually independent real-valued $(\tilde{\mathcal{F}}_t)$-Brownian motions $\{\tilde{\beta}_k\}_{k\geq1}$ such that 
\begin{equation}\label{Wienertildeeps}
\tilde{W}=\sum_{k\geq1}\tilde{\beta}_k e_k
\end{equation}
in $C([0,T];\mathfrak{U})$.
\end{lemma}

Note: see \cite[Paragraph~4.1]{DaPratoZabczyk92} for the definition of a $Q$-Wiener process.\medskip

\textbf{Proof of Lemma~\ref{lem:tildeW}.} It is clear that $\tilde{W}$ is a $j\circ j^*$-cylindrical Wiener process (this notion is stable by convergence in law; actually it can be characterized in terms of the law of $\tilde{W}$ uniquely if we drop the usual hypothesis of a.s. continuity of the trajectories. This latter property of continuity can be recovered, after a possible modification of the process, by using Kolmogorov's Theorem). Also $\tilde{W}$ is
$(\tilde{\mathcal{F}}_t)$-adapted by definition of the filtration $(\tilde{\mathcal{F}}_t)$. By \cite[Proposition~4.1]{DaPratoZabczyk92}, we obtain the decomposition~\eqref{Wienertildeeps}. The $\tilde{\P}$-a.s. convergence of the sum in \eqref{Wienertildeeps} in the space $C([0,T];\mathfrak{U})$ is proved as in \cite[Theorem~4.3]{DaPratoZabczyk92}. \qed\medskip

Note that the last component $\tilde{W}^n$ of $\tilde{Z}^n$ depends a priori on $n$. Without loss of generality, we will replace $\tilde{W}^n$ by $\tilde{W}$. Of course, this does not affect the almost-sure convergence of $\tilde{Z}^n$ to $\tilde{Z}$, and Lemma~\ref{lem:tildeW} asserts that this does not modify the law of $\tilde{Z}^n$. This operation is not mandatory for the validity of what follows, and quite natural since the original sequence $(Z^n)$ is stationary (as a sequence) with respect to its last argument.

\paragraph{Step 3. Martingales}

\begin{proposition}\label{prop:3martingalestilde} Let $\varphi_j\in\Gamma$. Let $\tilde{h}_{j,k}(t)$ be defined by 
\begin{equation*}
\tilde{h}_{j,k}(t)=\int_{\T^N}\int_\R g_k(x,\xi)\varphi_j(x,\xi)d\tilde{\nu}_{x,t}(\xi)dx.
\end{equation*} 
Then, for $j\in\N$, $k\geq 1$, the processes
\begin{equation}\label{3martingalestilde}
\tilde{M}_j(t),\quad \tilde{M}_j(t)\tilde{\beta}_k(t)-\int_0^t \tilde{h}_{j,k}(s)ds,\quad |\tilde{M}_j(t)|^2-\int_0^t\|\tilde{h}_j(s)\|_{l^2(\N^*)}^2 ds,
\end{equation}
and $(\tilde{W}(t))$ are $(\tilde{\mathcal{F}}_t)$-martingales.
\end{proposition}

\textbf{Proof of Proposition~\ref{prop:3martingalestilde}.} The proof is similar to the proof of \cite[Proposition~1.1 p.522]{JacodShiryaev03}, except that we do not use any hypothesis of boundedness here since we use the $\tilde{\P}$-almost-sure convergence  and the Vitali Theorem to pass to the limit in the expectation of the quantities of interest (an other minor difference with the proof of \cite[Proposition~1.1 p.522]{JacodShiryaev03} is that $\tilde{M}$ is known to be continuous $\tilde{\P}$-a.s., not only c\`adl\`ag). \smallskip

Let $t_1,t_2\in\R_+$, $t_1<t_2$ and let $H$ be a $\mathscr{D}_{t_1-}(E)$-measurable bounded function as in \eqref{defHHHH}.
By identities of the laws of $M^n_\varphi$ and $\tilde{M}^n_\varphi$, we have
$$
\tilde{\E}|\tilde{M}^n_{\varphi_j}(t_2)-\tilde{M}^n_{\varphi_j}(t_1)|^2=\E|M^n_{\varphi_j}(t_2)-M^n_{\varphi_j}(t_1)|^2.
$$
Using \eqref{IncrementsMn}, it follows that
$$
\sup_n\tilde{\E}\left|H\left(\{\tilde{f}^n\},\{\tilde{M}^n\},\tilde{W}\right)\left[\tilde{M}^n_{\varphi_j}(t_2)-\tilde{M}^n_{\varphi_j}(t_1)\right]\right|^2<+\infty,
$$
since $H$ is bounded. We have in addition
\begin{equation}\label{tripleae}
\left(\{\tilde{f}^n\},\{\tilde{M}^n\},\tilde{W}\right)\to\left(\{\tilde{f}\},\{\tilde{M}\},\tilde{W}\right)
\end{equation}
a.e., $\tilde{\P}$-almost-surely by Proposition~\ref{prop:CVfftilde} and thus,
$$
H\left(\{\tilde{f}^n,\},\{\tilde{M}^n\},\tilde{W}\right)\to H\left(\{\tilde{f}\},\{\tilde{M}\},\tilde{W}\right),
$$
$\tilde{\P}$-almost-surely. There is also convergence
$$
\tilde{M}^n_{\varphi_j}(t_2)-\tilde{M}^n_{\varphi_j}(t_1)\to\tilde{M}_j(t_2)-\tilde{M}_j(t_1)
$$
$\tilde{\P}$-almost-surely. By Vitali's Theorem, we obtain
\begin{multline}\label{cvMnj}
\tilde{\E}\left[H\left(\{\tilde{f}^n,\},\{\tilde{M}^n\},\tilde{W}\right)\left(\tilde{M}^n_{\varphi_j}(t_2)-\tilde{M}^n_{\varphi_j}(t_1)\right)\right]\\
\to 
\tilde{\E}\left[H\left(\{\tilde{f}\},\{\tilde{M}\},\tilde{W}\right)\left(\tilde{M}_j(t_2)-\tilde{M}_j(t_1)\right)\right].
\end{multline}
By identities of the laws, the left-hand side of \eqref{cvMnj} is
$$
\E\left[H\left(\{f^n,\},\{M^n\},W\right)\left(M^n_{\varphi_j}(t_2)-M^n_{\varphi_j}(t_1)\right)\right]=0,
$$
since $M^n_\varphi$ is a $(\mathcal{F}_t)$-martingale. We deduce from \eqref{cvMnj} thus that
\begin{equation}\label{barMjMartingale1}
\tilde{\E}\left[H\left(\{\tilde{f}\},\{\tilde{M}\},\tilde{W}\right)\left(\tilde{M}_j(t_2)-\tilde{M}_j(t_1)\right)\right]=0.
\end{equation}
Due to Proposition~\ref{prop:Dtminus}, we deduce from \eqref{barMjMartingale1} that
\begin{equation}\label{barMjMartingale1A}
\tilde{\E}\left[\mathbf{1}_A\left(\{\tilde{f}\},\{\tilde{M}\},\tilde{W}\right)\left(\tilde{M}_j(t_2)-\tilde{M}_j(t_1)\right)\right]=0,
\end{equation}
for all cylindrical sets $A$ as in \eqref{defAAAA}. The left-hand side of \eqref{barMjMartingale1A} defines a finite measure (due to \eqref{IncrementsMn}) which coincides with the trivial measure $A\mapsto 0$ for sets $A$ as in \eqref{defAAAA}. Since such sets form a $\pi$-system which generates $\mathcal{D}_{t_1-}(E)$, hence a separating class, we deduce that \eqref{barMjMartingale1A} holds true for all $A\in\mathcal{D}_{t_1-}(E)$. It follows then also that \eqref{barMjMartingale1} is satisfied for all $\mathcal{D}_{t_1-}(E)$-measurable bounded function $H$. Let now $s,t\in[0,T)$ with $s<t$. Let $(s_n)$ and $(t_n)$ be some decreasing sequences in $\R_+$, converging to $s$ and $t$ respectively. Let $H$ be a $\mathscr{D}_{s}(E)$-measurable bounded function. Then $H$ is a $\mathscr{D}_{s_n-}(\R^{2+m})$-measurable bounded function since $s<s_n$. By passing to the limit in \eqref{barMjMartingale1} written with $t_1=s_n$, $t_2=t_n$ (we use the right-continuity of the processes here), we obtain
\begin{equation}\label{barMjMartingale2}
\tilde{\E}\left[H\left(\{\tilde{f}\},\{\tilde{M}\},\tilde{W}\right)\left(\tilde{M}_j(t)-\tilde{M}_j(s)\right)\right]=0.
\end{equation}
This shows that $(\tilde{M}_j(t))$ is a $\tilde{\mathcal{F}}_t$-martingale. The proof that $(\tilde{W}(t))$ is a $\tilde{\mathcal{F}}_t$-martingale is similar, we do not give the details of that point. To go on, let us define now the processes
$$
\tilde{H}^n_{j,k}(t)=\int_0^t\tilde{h}^n_{j,k}(s)ds,\quad\tilde{H}_{j,k}(t)=\int_0^t\tilde{h}_{j,k}(s)ds,
$$ 
and
$$
\tilde{\mathcal{H}}^n_{j}(t)=\int_0^t \|\tilde{h}^n_j(s)\|_{l^2(\N^*)}^2 ds,\quad\tilde{\mathcal{H}}_{j}(t)=\int_0^t\|\tilde{h}_j(s)\|_{l^2(\N^*)}^2 ds,
$$
and the processes
\begin{equation*}
\begin{cases}
\displaystyle\tilde Y^n_{j,k}(t)=\tilde{M}^n_j(t)\tilde{\beta}_k(t)-\tilde{H}^n_{j,k}(t), 
&\displaystyle \tilde{Y}_{j,k}(t)=\tilde{M}_j(t)\tilde{\beta}_k(t)-\tilde{H}_{j,k}(t),\\
\\
\displaystyle\tilde{V}^n_j(t)= |\tilde{M}^n_j(t)|^2-\tilde{\mathcal{H}}^n_{j}(t),
&\displaystyle\tilde{V}_j(t)= |\tilde{M}_j(t)|^2-\tilde{\mathcal{H}}_{j}(t).
\end{cases}
\end{equation*}

To complete the proof of Proposition~\ref{prop:3martingalestilde}, we have to show that $(\tilde{Y}_k(t))$ and $(\tilde{V}_j(t))$ are $\tilde{\mathcal{F}}_t$-martingale. We will use the following result.

\begin{lemma}\label{lem:cvaehnjk} Let $T>0$. Then, up to a subsequence, for all $j\in\N$, $k\in\N^*$, $\tilde{\P}$-almost-surely, $\tilde{h}^{n}_{j,k}\to \tilde{h}_{j,k}$ and $\|\tilde{h}^{n}_j(\cdot)\|_{l^2(\N^*)}^2\to\|\tilde{h}_j(\cdot)\|_{l^2(\N^*)}^2$ in $L^1(0,T)$,  when $n\to+\infty$. 
\end{lemma}

Lemma~\ref{lem:cvaehnjk} implies that, $\tilde{\P}$-almost-surely, for every $t\in[0,T]$, $\tilde{H}^{n}_{j,k}(t)$ and $\tilde{\mathcal{H}}^{n}_{j}(t)$ are converging to $\tilde{H}_{j,k}(t)$ and $\tilde{\mathcal{H}}_{j}(t)$ respectively. We have also $\tilde{M}^n_j\to\tilde{M}_j$ in $C(\R_+)$, from which follows the convergences
$\tilde{M}^n_j\tilde{\beta}_k\to\tilde{M}_j\tilde{\beta}_k$ and 
$|\tilde{M}^n_j|^2\to|\tilde{M}_j|^2$ in $C(\R_+)$, $\tilde{\P}$-almost-surely. We deduce that, $\tilde{\P}$-almost-surely, 
\begin{equation}\label{cvYZ}
\tilde Y^n_{j,k}(t)\to\tilde{Y}_{j,k}(t),\quad\tilde V^n_j(t)\to\tilde{V}_j(t), 
\end{equation}
for all $t\geq 0$. With the estimate \eqref{IncrementsMn}, it is easy to obtain the bounds 
\begin{equation}\label{EquiIntYZ}
\tilde{\E}|\tilde Y^n_{j,k}(t)-\tilde{Y}^n_{j,k}(s)|^2\leq C,\quad \tilde{\E}|\tilde V^n_j(t)-\tilde{V}^n_j(s)|^2\leq C,
\end{equation}
where the constant $C$ depend on $s,t\in[0,T]$, $k$, but not on $n$. By \eqref{cvYZ}	and \eqref{EquiIntYZ} (this last condition shows the equi-integrability of $(\tilde Y^n_{j,k}(t)-\tilde Y^n_{j,k}(s))$ and $(\tilde V^n_j(t)-\tilde V^n_j(s))$ respectively), we can use the arguments applied to the martingale $\tilde{M}^n_\varphi(t)$ in the first part of the proof: it will establish that $\tilde{Y}_{j,k}(t)$ and $\tilde{V}_j(t)$ are $(\tilde{\F}_t)$-martingales. \qed\medskip

Let us now give the\smallskip

\textbf{Proof of Lemma~\ref{lem:cvaehnjk}.} Let us first show that, for all $j,k$, we have the following convergence :
\begin{equation}\label{CVHL2}
\tilde{h}^n_{j,k}\to\tilde{h}_{j,k}\quad\mbox{in}\quad L^2((0,T)\times\tilde\Omega).
\end{equation}
Define, for every $\psi\in C_b(\T^N\times\R)$,
\begin{equation}\label{hpsi}
\tilde{h}^n_{\psi}(t)=\int_{\T^N}\int_\R \psi(x,\xi)d\tilde{\nu}^n_{x,t}(\xi)dx,\quad 
\tilde{h}_{\psi}(t)=\int_{\T^N}\int_\R \psi(x,\xi)d\tilde{\nu}_{x,t}(\xi)dx.
\end{equation}
If  $\psi\in C^1_c(\T^N\times\R)$, then $\tilde{h}^n_{\psi}(t)=\<\tilde{f}^n(t),\partial_\xi\psi\>$. By Proposition~\ref{prop:CVfftilde}, we have then, $\tilde{\P}$-almost-surely,
\begin{equation}\label{eq:cvhpsi}
\mbox{for all }t\in[0,T]\setminus\tilde{B},\;\tilde{h}^n_{\psi}(t)\to\tilde{h}_{\psi}(t).
\end{equation}
Using the Jensen inequality, we have
\begin{equation}\label{Domhl2}
\|\tilde{h}^n_{\psi}-\tilde{h}_\psi\|_{L^2((0,T)\times\tilde\Omega)}^2\leq 4T\|\psi\|^2_{C_b(\T^N\times\R)}.
\end{equation}
By the Vitali Theorem, we obtain the convergence $\tilde{h}^n_{\psi}\to\tilde{h}_\psi$ in $L^1((0,T)\times\tilde\Omega)$. Using \eqref{Domhl2} also, we see that this convergence can be extended to the case of a general integrand $\psi\in C_b(\T^N\times\R)$. Let us then take $\psi=g_k\varphi_j$. We obtain first 
$\tilde{h}^n_{j,k}\to\tilde{h}_{j,k}$ in $L^2((0,T)\times\tilde\Omega)$. It follows that, up to a subsequnce, $\tilde{\P}$-almost-surely, $\tilde{h}^n_{j,k}\to\tilde{h}_{j,k}$ in $L^2(0,T)$, hence in $L^1(0,T)$. The subsequence and the $\tilde{\P}$-almost-sure property can be made independent on $j,k$ since $\Gamma\times\N^*$ is countable. The growth hypothesis \eqref{D0} also shows that
$$
\sum_k\|\tilde{h}^n_{j,k}-\tilde{h}_{j,k}\|_{L^2((0,T)\times\tilde\Omega)}^2\leq 4D_0(1+C_2)T\|\varphi_j\|^2_{C_b(\T^N\times\R)}.
$$
Again, using the dominated convergence theorem, we deduce that 
$$
\|\tilde{h}^{n}_j\|_{l^2(\N^*)}^2\to\|\tilde{h}_j\|_{l^2(\N^*)}^2
$$ 
in $L^1((0,T)\times\tilde{\Omega})$, which allows to conclude the proof of the lemma. \qed

\paragraph{Step 4. Conclusion of the martingale method}\label{sec:step4M}

Let us first prove that $M^*_\varphi(t)$ given in \eqref{def:tildeMphi} is well-defined. 

\begin{lemma}\label{lem:PnuOK} Item~\ref{item:Pnu} in Theorem~\ref{th:martingalesol} is satisfied, \textit{i.e.}: for all $\psi\in C_b(\R)$,
$(x,t)\mapsto\<\psi,\tilde\nu_{x,t}\>$ belongs to $L^2_{\tilde{\mathcal{P}}}(\T^N\times [0,T]\times\tilde\Omega)$. 
\end{lemma}

\textbf{Proof of Lemma~\ref{lem:PnuOK}.} For $\psi\in C_b(\R)$, set 
$
\tilde{X}_\psi(x,t)=\<\psi,\tilde\nu_{x,t}\>.
$
We have $\tilde{X}_\varphi\in L^2(\T^N\times [0,T]\times\tilde\Omega)$, with
\begin{equation}\label{BountildeXpsi}
\tilde{\E}\|\tilde{X}_\psi\|^2_{L^2(\T^N\times[0,T])}\leq \|\psi\|_{C_b(\R)}^2 T.
\end{equation}
If $\theta\in C(\T^N)$, and if $\psi$ is $C^1$, vanishes in the neighbourhood of $-\infty$ and satisfies $\psi'\in C_c(\R)$, then, due to \eqref{fVSnu}, we have
$$
\<\tilde{X}_\psi(t),\theta\>_{L^2(\T^N)}=\<\tilde{f}(t),\varphi\>,\quad\varphi(x,\xi):=\theta(x)\psi'(\xi).
$$ 
By Item~\ref{item:tildefcadlag} of Theorem~\ref{th:martingalesol}, the process $Y_t:=\<\tilde{X}_\psi(t),\theta\>_{L^2(\T^N)}$ is c\`adl\`ag. Since $(Y_t)$ is adapted by definition of $(\tilde{\F}_t)$, it is an optional process \cite[p.~172]{RevuzYor99}. In par\-ti\-cu\-lar, $(Y_t)$ is progressively measurable \cite[Proposition~4.8]{RevuzYor99}, hence $Y\in L^2_{\tilde{\mathcal{P}}}([0,T]\times\tilde\Omega)$. 
A limiting argument (by approximation and truncation of the function $\psi$ in particular), using \eqref{BountildeXpsi} and the fact that $\tilde{\nu}$ vanishes at infinity shows that the result holds true when $\psi$ is merely a function in $C_b(\R)$ and $\theta$ any function in $L^2(\T^N)$. We obtain, therefore, that, for all $\psi\in C_b(\R)$, $X_\psi$ belongs to $L^2_{\tilde{\mathcal{P}}}([0,T]\times\tilde\Omega; L^2(\T^N)-\mbox{weak})$. Since being weakly or strongly $\tilde{\P}$-measurable is the same thing, (\textit{cf.} Section~\ref{sec:predictable}), we have established the result.\qed\medskip

We can apply now Proposition~\ref{prop:3martingalesCV}. Indeed, due to Lemma~\ref{lem:PnuOK}, the processes $\tilde{h}_{j,k}$ in 
Proposition~\ref{prop:3martingalestilde} are in $L^2_{\tilde{\mathcal{P}}}([0,T]\times\tilde\Omega)$. By the martingale property~\eqref{3martingalestilde}, we conclude that $\tilde{M}_\varphi(t)=M^*_\varphi(t)$, with $M^*_\varphi(t)$ defined by \eqref{def:tildeMphi}, for every $\varphi\in\Gamma$.

\subsubsection{Identification of the limit: equation}\label{sec:proofThmartingaleEquation}

We prove now \eqref{eq:kineticfpretilde}. Let $\varphi\in\Gamma$. By item~\ref{item:3thmartingale} and item~\ref{item:8thmartingale} of Theorem~\ref{th:martingalesol}, using also the identity $\tilde{M}_\varphi(t)=M^*_\varphi(t)$, we have the identification
\begin{multline*}
\tilde{J}_\varphi(t)=\int_0^t \<\tilde{f}(s),a(\xi)\cdot\nabla\varphi\>ds+\sum_{k\geq 1}\int_0^t\int_{\T^N}\int_\R g_k(x,\xi)\varphi(x,\xi)d\tilde{\nu}_{x,s}(\xi)dx d\tilde{\beta}_k(s)\\
+\frac{1}{2}\int_0^t\int_{\T^N}\int_\R \partial_\xi\varphi(x,\xi)\GG^2(x,\xi)d\tilde{\nu}_{x,s}(\xi) dx ds.
\end{multline*}
The equation \eqref{eq:kineticfpretilde} follows therefore from the identity \eqref{negN0tilde}.

\subsection{Pathwise solutions and almost-sure convergence}\label{sec:pathCV}

If $f_0$ is at equilibrium in Theorem~\ref{th:martingalesol}, then we have seen in Theorem~\ref{th:Uadd} that \eqref{stoSCL} admits a unique solution for a given initial datum. We can use this uniqueness result to obtain existence of pathwise solution and convergence in $L^p$ of the sequence of approximate solutions in that case.

\begin{theorem}[Pathwise solution] Suppose that there exists a sequence of approximate generalized solutions $(f^n)$ to~\eqref{stoSCL} with initial datum $f_0^n$ satisfying \eqref{eq:integrabilityfn}, \eqref{Boundmn} and the tightness condition \eqref{inftymn} and such that $(f^n_0)$ converges to the equilibrium function $\mathtt{f}_0(\xi)=\mathbf{1}_{u_0>\xi}$ in $L^\infty(\T^N\times\R)$-weak-*, where $u_0\in L^\infty(\T^N)$. We have then
\begin{enumerate}
\item there exists a unique solution $u\in L^1(\T^N\times [0,T]\times\Omega)$ to \eqref{stoSCL} with initial datum $u_0$;
\item let
$$
u^n(x,t)=\int_\R\xi d\nu^n_{x,t}(\xi)=\int_\R\left( f^n(x,t,\xi)-\mathbf{1}_{0>\xi} \right)d\xi.
$$
Then, for all $p\in[1,\infty[$, $(u^n)$ is converging to $u$ with the following two different modes of convergence: $u_n\to u$ in $L^p(\T^N\times(0,T)\times\Omega)$ and, for a subsequence $(n_k)$, almost surely, for all $t\in[0,T]$, $u^{n_k}(t)\to u(t)$ in $L^p(\T^N)$.
\end{enumerate}
\label{th:pathcv} \end{theorem}

\textbf{Proof of Theorem~\ref{th:pathcv}.} We use the Gy\"ongy-Krylov argument, \cite[Lemma~1.1]{GyongyKrylov96} (the basis of the Gy\"ongy-Krylov argument is this simple fact: if a couple $(X_n,Y_n)$ of random variables converges in law to a random variable written $(Z,Z)$, \textit{i.e.} concentrated on the diagonal, then $X_n-Y_n$ converges to $0$ in probability). Let us go back to Section~\ref{sec:proofThmartingalestatespace}. We introduce the random variable
$$
Z^{n,q}=(\nu^n,\{J^n\},\{M^n\},\{\eps^n\},\{f^n_\mathrm{in}\},\mu^n,\|m^n\|,\nu^q,\{J^q\},\{M^q\},\{\eps^q\},\{f^q_\mathrm{in}\},\mu^q,\|m^q\|,W)
$$
in the state space $\mathcal{Z}$ equal to
$$
\mathcal{Y}^1\times\mathcal{E}\times\mathcal{P}^1(\T^N\times[0,T]\times\R)\times\R_+\times\mathcal{Y}^1\times\mathcal{E}\times\mathcal{P}^1(\T^N\times[0,T]\times\R)\times\R_+\times\mathcal{X}_W.
$$
We repeat the arguments used in Section~\ref{sec:proofThmartingale} to show that $Z^{n,q}$ is tight in $\mathcal{Z}$ and that there exists a probability space $(\tilde{\Omega},\tilde{\mathcal{F}},\tilde{\P})$ and a new random variable $\tilde{Z}^{n,q}$ with the same law as $Z^{n,q}$, such that a subsequence $(\tilde{Z}^{n_l,q_l})_l$ is converging $\tilde{\P}$-almost-surely in $\mathcal{Z}$ to a random variable $\tilde{Z}$. Let $\tilde{\nu}$ be the the first component of $\tilde{Z}$ and $\check{\tilde{\nu}}$ be the seventh component of $\tilde{Z}$. Repeating all steps from Section~\ref{sec:proofThmartingalecadlag}, \ref{sec:proofThmartingalecvsto}, \ref{sec:proofThmartingaleEquation}, we obtain two generalized solutions
$$
\tilde{f}(x,t,\xi)=\tilde{\nu}_{(x,t)}(\xi,+\infty),\quad \check{\tilde{f}}(x,t,\xi)=\check{\tilde{\nu}}_{(x,t)}(\xi,+\infty),
$$
to Equation \eqref{stoSCL} with probabilistic data $(\tilde{\Omega},\tilde{\mathcal{F}},\tilde{\P},(\tilde{\mathcal{F}}_t),\tilde{W})$, where $(\tilde{\mathcal{F}}_t)$ is the completion of the filtration generated by the five-uplet $(\{\tilde{f}\},\{\tilde{M}\},\{\check{\tilde{f}}\},\{\check{\tilde{M}}\},\tilde{W})$:
\begin{equation*}
\tilde{\mathcal{F}}_t=\sigma\big((\{\tilde{f}\},\{\tilde{M}\},\{\check{\tilde{f}}\},\{\check{\tilde{M}}\},\tilde{W})^{-1}\left(\mathscr{D}_t(E)\times\mathscr{D}_t(\check{E})\right)\cup\big\{N\in\tilde{\mathcal{F}};\;\tilde{\P}(N)=0\big\}\big),
\end{equation*}
for $t\in[0,T]$, with 
$$
E:=\R^\infty\times\R^\infty,\quad \check{E}:=\R^\infty\times\R^\infty\times\mathfrak{U}.
$$
Note that $\mathscr{D}_t(E)\times\mathscr{D}_t(\check{E})\not=\mathscr{D}_t(E\times\check{E})$ since the natural topologies of $D(\R_+;E)\times D(\R_+;\check{E})$ and $D(\R_+;E\times\check{E})$ are different (the topology of the former is the product topology of the Skorokhod topologies on each space: this authorizes two changes of times, one for each coordinate; for the Skorokhod topology on $D(\R_+;E\times\check{E})$, only one change of time is admissible). The solutions $\tilde{f}$ and $\check{\tilde{f}}$ have the same initial condition $f_0$, which is an equilibrium function $\mathtt{f}_0$. By Theorem~\ref{th:Uadd}, we have 
\begin{equation}\label{tildefcheckf}
\tilde{f}=\check{\tilde{f}}=\mathtt{f},
\end{equation}
where $\mathtt{f}$ is the equilibrium function $\mathbf{1}_{\tilde{u}>\xi}$, where
$$
\tilde u(x,t):=\int_\R\xi d\tilde\nu_{(x,t)}(\xi).
$$
A first consequence of \eqref{tildefcheckf} is that $\tilde{\nu}=\check{\tilde{\nu}}$, $\tilde{\P}$-almost-surely. By Remark~\ref{rk:uniqmmm} on the uniqueness of the kinetic measure, we have also  $\tilde{m}=\check{\tilde{m}}$, $\tilde{\P}$-almost-surely.
We apply the Gy\"ongy-Krylov argument: we obtain that $(\nu^n)$ is converging in probability in $\mathcal{Y}^1$ and $(m^n)$ is converging in probability in $\mathcal{M}_b(\T^N\times[0,T]\times\R)$-weak-*. Extracting an additional subsequence if necessary, we can assume that the convergences are also $\P$-almost-sure. By the arguments of the sections~\ref{sec:proofThmartingalecadlag}, \ref{sec:proofThmartingalecvsto}, \ref{sec:proofThmartingaleEquation}, it follows that $f(t,x,\xi):=\nu_{x,t}(\xi,+\infty)$ is a generalized solution to \eqref{stoSCL}. Note, to give few details, that we do not need to follow Step~1. and Step~2. of Section~\ref{sec:proofThmartingalecvsto} here, since the filtration $(\mathcal{F}_t)$ and the Wiener processes $\beta_k(t)$ are already known here. The convergence of the stochastic integral in $J^n_\varphi(t)$ does not require the martingale method of Step~3. of Section~\ref{sec:proofThmartingalecvsto} either. Using the $L^2$ convergence of the integrand (\textit{cf.} Lemma~\ref{lem:cvaehnjk}) is sufficient by the It\^o isometry.\smallskip

We use the second identity in \eqref{tildefcheckf} now. It states, equivalently, that $\tilde{\P}$-almost-surely, for a.e. $(x,t)$, $\tilde\nu_{(x,t)}=\delta_{\tilde u(x,t)}$. The fact that $\tilde\nu$ is a Dirac mass can be characterized in terms of equality in the Jensen Inequality:
\begin{equation}\label{EJensen}
\tilde{\E}\iint_{\T^N\times(0,T)}\Phi\left(\int_\R\xi d\tilde{\nu}_{(x,t)}(\xi)\right) dx dt=\tilde{\E}\iint_{\T^N\times(0,T)}\int_\R\Phi(\xi)d\tilde{\nu}_{(x,t)}(\xi) dx dt,
\end{equation}
where $\Phi$ is a strictly convex, polynomially bounded function, like $\Phi(\xi)=\xi^2$ for example.
The identity \eqref{EJensen} depends on $\mathrm{Law}(\tilde{\nu})=\mathrm{Law}(\nu)$ uniquely. Therefore $\nu$ also is almost-surely a Dirac mass: $\P$-almost-surely, for a.e. $(x,t)$, $\nu_{(x,t)}=\delta_{u(x,t)}$, where
$$
u(x,t):=\int_\R\xi d\nu_{(x,t)}(\xi).
$$
(Remark that $\nu_{(x,t)}=\delta_{u(x,t)}$ a.s., a.e., is also a consequence of Theorem~\ref{th:Uadd}. However this theorem is difficult to show, and, although we have already used Theorem~\ref{th:Uadd}, the argument based on \eqref{EJensen} is simple and natural). By  Proposition~\ref{prop:gestoes}, $u$ is a solution to \eqref{stoSCL}: it is the unique solution by Theorem~\ref{th:Uadd}. 
Using Lemma~\ref{lem:weakstrongeq}, also, we have $\|u^n-u\|_{L^p(\T^N\times(0,T))}^p\to 0$ in probability. We also have the uniform bound
\begin{equation}\label{appCVproba}
\E\|u^n-u\|_{L^p(\T^N\times(0,T))}^{pr}\leq C,
\end{equation}
where $r>1$ and $C$ is independent on $n$. Taking \eqref{appCVproba} for granted, we deduce, with the convergence in probability, that 
$\E\|u^n-u\|_{L^p(\T^N\times(0,T))}^p\to 0$ and we obtain the first part of the second point of Theorem~\ref{th:pathcv}. The bound \eqref{appCVproba} follows from the estimate 
$$
\E\|u^n-u\|_{L^p(\T^N\times(0,T))}^{pr}\leq \E\|u^n-u\|_{L^{pr}(\T^N\times(0,T))}^{pr} T^{r-1},
$$
and \eqref{eq:integrabilityu}, \eqref{eq:integrabilityfn}. To prove that almost surely, for all $t\in[0,T]$, $u_{n_k}(t)\to u(t)$ in $L^p(\T^N)$, we use Lemma~\ref{lem:weakstrongeq} and Proposition~\ref{prop:CVfftilde}. It gives: almost surely, for all $t\in[0,T]\setminus B_\mathrm{at}$, $u_{n_k}(t)\to u(t)$ in $L^p(\T^N)$, where $B_\mathrm{at}$ is defined in Lemma~\ref{lem:atomicpoints}. Since, almost-surely, $u$ is continuous in time with values in $L^p(\T^N)$ by Corollary~\ref{cor:timecontinuity}, it follows from \eqref{f+VSf-} that $B_\mathrm{at}$ is empty. This gives the desired result. \qed

\section{Some applications}\label{sec:appli}

\subsection{Vanishing viscosity method}\label{sec:appParabolic}
Assume that \eqref{D0} and \eqref{D1} are satisfied. Consider the parabolic approximation to \eqref{stoSCL}:
\begin{equation}\label{eq:parabolic}
du^\eta+\div(A(u^\eta))dt-\eta\Delta u^\eta dt=\Phi^\eta(x,u^\eta) dW(t).
\end{equation}
For $\eta>0$ and $u_0^\eta\in L^\infty(\T^N)$, the existence of solutions to \eqref{eq:parabolic} under the initial condition $u^\eta(0)=u^\eta_0$ has been proved in \cite{GyongyRovira00} provided the noise has a finite number of components. Therefore, we assume (compare to \eqref{defPhigk})
\begin{equation}\label{defPhigketa}
\Phi^\eta(x,u)=\sum_{1\leq k\leq K_\eta}g_k(x,u) e_k,
\end{equation}
where $K_\eta$ is finite, $K_\eta\to+\infty$ when $\eta\to0$. Let $(\eta_n)\downarrow 0$. In \cite{DebusscheVovelle10}, we have shown that the sequence $(u^{\eta_n})$ gives rise to a sequence of approximate generalized solutions $(f^n)$, with random measure $m^n$, given by
\begin{align*}
f^n&=\mathtt{f}^n=\mathbf{1}_{u^{\eta_n}>\xi},\\
\<m^n,\varphi\>&=\iint_{\T^N\times(0,T)}\varphi(x,t){\eta_n} |\nabla_x u^{\eta_n}(x,t)|^2 dx dt,\\
\eps_n(t,\varphi)&=\eta_n\int_0^t\int_{\T^N} \mathtt{f}^n(x,s,\xi)\Delta\varphi(x,\xi) d\xi dx ds.
\end{align*}
Here the order is $\order=2$. Let $p\in[1,+\infty)$. By Theorem~\ref{th:pathcv}, we recover the result given in \cite{DebusscheVovelle10} of convergence $u^\eta\to u$ in $L^p(\T^N\times(0,T)\times\Omega)$, where $u$ is the solution to \eqref{stoSCL} with initial datum $u_0$. We also obtain that, if $(\eta_n)\downarrow 0$, then, for a subsequence $(n_k)$, almost surely, for all $t\in[0,T]$, $u^{\eta_{n_k}}(t)\to u(t)$ in $L^p(\T^N)$.

\subsection{BGK approximation}\label{sec:appBGK}
We consider now the following BGK approximation to \eqref{stoSCL}:
\begin{align}
df^\eta+a(\xi)\cdot\nabla_x f^\eta dt=&\frac{\mathtt{f}^\eta-f^\eta}{\eta}dt-\partial_\xi f^\eta \Phi^\eta dW(t)-\frac12\partial_\xi(\GG^2\partial_\xi(-f^\eta)),\label{BGK1}\\
\mathtt{f}^\eta:=&\mathbf{1}_{u^\eta>\xi},\quad u^\eta=\int_\R (f^\eta (\xi)-\mathbf{1}_{0>\xi})d\xi.\label{BGK2}
\end{align}
Assume \eqref{defPhigketa}, assume that \eqref{D1} is satisfied and that (instead of \eqref{D0}), either $\GG^2(x,\xi)\leq D_0|\xi^2|$ or $\GG^2(x,\xi)\leq D_0$ is satisfied. M.~Hofmanov{\'a} has proved in \cite{HofmanovaBGK15} the existence ot solutions to \eqref{BGK1}-\eqref{BGK2} with given initial datum $f^\eta_0=\mathtt{f}^\eta_0=\mathbf{1}_{u^\eta_0>\xi}$ (the fact that the initial datum is at equilibrium can be relaxed). Let $(\eta_n)\downarrow 0$. Then $f^n:=f^{\eta_n}$ provides a sequence of generalized approximate solutions of order $\order=0$, with 
\begin{align*}
\partial_\xi m^n&=\frac{\mathtt{f}^n-f^n}{\eta_n},\\
\eps_n(t,\varphi)&=0.
\end{align*}
Let $u_0\in L^\infty(\T^N)$. Assume $u_0^{\eta_n}\to u_0$ in $L^p(\T^N)$ for all $p\in[1,+\infty)$, and let $u$ be the solution to \eqref{stoSCL} with initial datum $u_0$. By Theorem~\ref{th:pathcv}, we recover the convergence $u^n\to u$ proved in \cite{HofmanovaBGK15}. We have also: for a subsequence $(n_k)$, almost surely, for all $t\in[0,T]$, $u^{n_k}(t)\to u(t)$ in $L^p(\T^N)$

\subsection{Approximation by the Finite Volume method}\label{sec:appFV}

The approximation of \eqref{stoSCL} by the Finite Volume method is considered in the companion paper \cite{DottiVovelle16b}.


\begin{thebibliography}{10}

\bibitem{Bauzet15}
C.~Bauzet.
\newblock Time-splitting approximation of the {C}auchy problem for a stochastic
  conservation law.
\newblock {\em Math. Comput. Simulation}, 118:73--86, 2015.

\bibitem{BauzetCharrierGallouet14a}
C.~Bauzet, J.~Charrier, and T.~Gallou{\"e}t.
\newblock Convergence of flux-splitting finite volume schemes for hyperbolic
  scalar conservation laws with a multiplicative stochastic perturbation.
\newblock {\em Math. Comp.}, 85(302):2777--2813, 2016.

\bibitem{BauzetCharrierGallouet14b}
C.~Bauzet, J.~Charrier, and T.~Gallou{\"e}t.
\newblock Convergence of monotone finite volume schemes for hyperbolic scalar
  conservation laws with multiplicative noise.
\newblock {\em Stoch. Partial Differ. Equ. Anal. Comput.}, 4(1):150--223, 2016.

\bibitem{BauzetValletWittbold12}
C.~Bauzet, G.~Vallet, and P.~Wittbold.
\newblock The {C}auchy problem for conservation laws with a multiplicative
  stochastic perturbation.
\newblock {\em J. Hyperbolic Differ. Equ.}, 9(4):661--709, 2012.

\bibitem{BauzetValletWittbold14}
C.~Bauzet, G.~Vallet, and P.~Wittbold.
\newblock The {D}irichlet problem for a conservation law with a multiplicative
  stochastic perturbation.
\newblock {\em J. Funct. Anal.}, 266(4):2503--2545, 2014.

\bibitem{BillingsleyBook}
P.~Billingsley.
\newblock {\em Convergence of probability measures}.
\newblock Wiley Series in Probability and Statistics: Probability and
  Statistics. John Wiley \& Sons Inc., New York, second edition, 1999.
\newblock A Wiley-Interscience Publication.

\bibitem{BrzezniakOndrejat11}
Z.~Brze{\'z}niak and M.~Ondrej{\'a}t.
\newblock Weak solutions to stochastic wave equations with values in
  {R}iemannian manifolds.
\newblock {\em Comm. Partial Differential Equations}, 36(9):1624--1653, 2011.

\bibitem{CastaingRaynauddeFitteValadier04}
C.~Castaing, P.~Raynaud~de Fitte, and M.~Valadier.
\newblock {\em Young measures on topological spaces}, volume 571 of {\em
  Mathematics and its Applications}.
\newblock Kluwer Academic Publishers, Dordrecht, 2004.
\newblock With applications in control theory and probability theory.

\bibitem{ChenDingKarlsen12}
G.-Q. Chen, Q.~Ding, and K.~H. Karlsen.
\newblock On nonlinear stochastic balance laws.
\newblock {\em Arch. Ration. Mech. Anal.}, 204(3):707--743, 2012.

\bibitem{ChenPerthame03}
G.-Q. Chen and B.~Perthame.
\newblock Well-posedness for non-isotropic degenerate parabolic-hyperbolic
  equations.
\newblock {\em Ann. Inst. H. Poincar\'e Anal. Non Lin\'eaire}, 20(4):645--668,
  2003.

\bibitem{ChungWilliams90}
K.~L. Chung and R.~J. Williams.
\newblock {\em Introduction to stochastic integration}.
\newblock Probability and its Applications. Birkh\"auser Boston, Inc., Boston,
  MA, second edition, 1990.

\bibitem{DaPratoZabczyk92}
G.~Da~Prato and J.~Zabczyk.
\newblock {\em Stochastic equations in infinite dimensions}, volume~44 of {\em
  Encyclopedia of Mathematics and its Applications}.
\newblock Cambridge University Press, Cambridge, 1992.

\bibitem{DebusscheHofmanovaVovelle15}
A.~Debussche, M.~Hofmanov{\'a}, and J.~Vovelle.
\newblock Degenerate parabolic stochastic partial differential equations:
  Quasilinear case.
\newblock {\em Annals of Probability}, 2015.

\bibitem{DebusscheVovelle10}
A.~Debussche and J.~Vovelle.
\newblock Scalar conservation laws with stochastic forcing.
\newblock {\em J. Funct. Anal.}, 259(4):1014--1042, 2010.

\bibitem{Diperna83a}
R.~J. DiPerna.
\newblock Convergence of approximate solutions to conservation laws.
\newblock {\em Arch. Rational Mech. Anal.}, 82(1):27--70, 1983.

\bibitem{DottiVovelle16b}
S.~Dotti and J.~Vovelle.
\newblock Convergence of the finite volume method for scalar conservation laws
  with multiplicative noise: an approach by kinetic formulation.
\newblock {\em hal-01391073}, 2016.

\bibitem{Droniou01}
J.~Droniou.
\newblock Int\'egration et espaces de sobolev \`a valeurs vectorielles.
\newblock {\em http://www-gm3.univ-mrs.fr/polys/}, 2001.

\bibitem{EKMS00}
W.~E, K.~Khanin, A.~Mazel, and Y.~Sinai.
\newblock Invariant measures for {B}urgers equation with stochastic forcing.
\newblock {\em Ann. of Math. (2)}, 151(3):877--960, 2000.

\bibitem{EGH00}
R.~Eymard, T.~Gallou{\"e}t, and R.~Herbin.
\newblock Finite volume methods.
\newblock In {\em Handbook of numerical analysis, {V}ol. {VII}}, Handb. Numer.
  Anal., VII, pages 713--1020. North-Holland, Amsterdam, 2000.

\bibitem{FengNualart08}
J.~Feng and D.~Nualart.
\newblock Stochastic scalar conservation laws.
\newblock {\em J. Funct. Anal.}, 255(2):313--373, 2008.

\bibitem{GessPerthameSouganidis2016}
B.~Gess, B.~Perthame, and P.~E. Souganidis.
\newblock Semi-discretization for stochastic scalar conservation laws with
  multiple rough fluxes.
\newblock {\em SIAM Journal on Numerical Analysis}, 54(4):2187--2209, 2016.

\bibitem{GessSouganidis2014}
B.~Gess and P.~E. Souganidis.
\newblock Long-time behavior, invariant measures and regularizing effects for
  stochastic scalar conservation laws.
\newblock {\em arXiv:1411.3939 [math]}, Nov. 2014.

\bibitem{GessSouganidis2015}
B.~Gess and P.~E. Souganidis.
\newblock Scalar conservation laws with multiple rough fluxes.
\newblock {\em Communications in Mathematical Sciences}, 13(6):1569--1597,
  2015.

\bibitem{GyongyKrylov96}
I.~Gy{\"o}ngy and N.~Krylov.
\newblock Existence of strong solutions for {I}t\^o's stochastic equations via
  approximations.
\newblock {\em Probab. Theory Related Fields}, 105(2):143--158, 1996.

\bibitem{GyongyRovira00}
I.~Gy{\"o}ngy and C.~Rovira.
\newblock On {$L^p$}-solutions of semilinear stochastic partial differential
  equations.
\newblock {\em Stochastic Process. Appl.}, 90(1):83--108, 2000.

\bibitem{Hofmanova13b}
M.~Hofmanov{\'a}.
\newblock Degenerate parabolic stochastic partial differential equations.
\newblock {\em Stochastic Process. Appl.}, 123(12):4294--4336, 2013.

\bibitem{HofmanovaBGK15}
M.~Hofmanov{\'a}.
\newblock A {B}hatnagar--{G}ross--{K}rook approximation to stochastic scalar
  conservation laws.
\newblock {\em Ann. Inst. Henri Poincar\'e Probab. Stat.}, 51(4):1500--1528,
  2015.

\bibitem{Hofmanova2016}
M.~Hofmanov{\'a}.
\newblock Scalar conservation laws with rough flux and stochastic forcing.
\newblock {\em Stochastic Partial Differential Equations. Analysis and
  Computations}, 4(3):635--690, 2016.

\bibitem{HofmanovaSeidler12}
M.~Hofmanov{\'a} and J.~Seidler.
\newblock On weak solutions of stochastic differential equations.
\newblock {\em Stoch. Anal. Appl.}, 30(1):100--121, 2012.

\bibitem{JacodShiryaev03}
J.~Jacod and A.~N. Shiryaev.
\newblock {\em Limit theorems for stochastic processes}, volume 288 of {\em
  Grundlehren der Mathematischen Wissenschaften [Fundamental Principles of
  Mathematical Sciences]}.
\newblock Springer-Verlag, Berlin, second edition, 2003.

\bibitem{KarlsenStorrosten2017}
K.~H. Karlsen and E.~B. {Storr{\o} sten}.
\newblock On stochastic conservation laws and {{Malliavin}} calculus.
\newblock {\em Journal of Functional Analysis}, 272(2):421--497, 2017.

\bibitem{Kim03}
Y.~Kim.
\newblock Asymptotic behavior of solutions to scalar conservation laws and
  optimal convergence orders to {$N$}-waves.
\newblock {\em J. Differential Equations}, 192(1):202--224, 2003.

\bibitem{KoleyMajeeVallet2016}
U.~Koley, A.~K. Majee, and G.~Vallet.
\newblock A finite difference scheme for conservation laws driven by {{Levy}}
  noise.
\newblock {\em arXiv:1604.07840 [math]}, Apr. 2016.

\bibitem{KrokerRohde12}
I.~Kr{\"o}ker and C.~Rohde.
\newblock Finite volume schemes for hyperbolic balance laws with multiplicative
  noise.
\newblock {\em Appl. Numer. Math.}, 62(4):441--456, 2012.

\bibitem{LionsPerthameSouganidis13}
P.-L. Lions, B.~Perthame, and P.~E. Souganidis.
\newblock Scalar conservation laws with rough (stochastic) fluxes.
\newblock {\em Stoch. Partial Differ. Equ. Anal. Comput.}, 1(4):664--686, 2013.

\bibitem{LionsPerthameSouganidis13a}
P.-L. Lions, B.~Perthame, and P.~E. Souganidis.
\newblock Stochastic averaging lemmas for kinetic equations.
\newblock In {\em S\'eminaire {L}aurent {S}chwartz---\'{E}quations aux
  d\'eriv\'ees partielles et applications. {A}nn\'ee 2011--2012}, S\'emin.
  \'Equ. D\'eriv. Partielles, pages Exp. No. XXVI, 17. \'Ecole Polytech.,
  Palaiseau, 2013.

\bibitem{LionsPerthameSouganidis14}
P.-L. Lions, B.~Perthame, and P.~E. Souganidis.
\newblock Scalar conservation laws with rough (stochastic) fluxes: the
  spatially dependent case.
\newblock {\em Stoch. Partial Differ. Equ. Anal. Comput.}, 2(4):517--538, 2014.

\bibitem{LionsPerthameTadmor94}
P.-L. Lions, B.~Perthame, and E.~Tadmor.
\newblock A kinetic formulation of multidimensional scalar conservation laws
  and related equations.
\newblock {\em J. Amer. Math. Soc.}, 7(1):169--191, 1994.

\bibitem{MohamedSeaidZahri13}
K.~Mohamed, M.~Seaid, and M.~Zahri.
\newblock A finite volume method for scalar conservation laws with stochastic
  time--space dependent flux functions.
\newblock {\em J. Comput. Appl. Math.}, 237(1):614--632, 2013.

\bibitem{Ondrejat10}
M.~Ondrej{\'a}t.
\newblock Stochastic nonlinear wave equations in local {S}obolev spaces.
\newblock {\em Electron. J. Probab.}, 15:no. 33, 1041--1091, 2010.

\bibitem{PerthameBook}
B.~Perthame.
\newblock {\em Kinetic formulation of conservation laws}, volume~21 of {\em
  Oxford Lecture Series in Mathematics and its Applications}.
\newblock Oxford University Press, Oxford, 2002.

\bibitem{RevuzYor99}
D.~Revuz and M.~Yor.
\newblock {\em Continuous martingales and {B}rownian motion}, volume 293 of
  {\em Grundlehren der Mathematischen Wissenschaften [Fundamental Principles of
  Mathematical Sciences]}.
\newblock Springer-Verlag, Berlin, third edition, 1999.

\bibitem{StorrostenKarlsen2016}
E.~B. Storr{\o}sten and K.~H. Karlsen.
\newblock Analysis of a splitting method for stochastic balance laws.
\newblock {\em arXiv:1601.02428 [math]}, Jan. 2016.

\bibitem{ValletWittbold09}
G.~Vallet and P.~Wittbold.
\newblock On a stochastic first-order hyperbolic equation in a bounded domain.
\newblock {\em Infin. Dimens. Anal. Quantum Probab. Relat. Top.},
  12(4):613--651, 2009.

\bibitem{Yosida80}
K.~Yosida.
\newblock {\em Functional analysis}, volume 123 of {\em Grundlehren der
  Mathematischen Wissenschaften [Fundamental Principles of Mathematical
  Sciences]}.
\newblock Springer-Verlag, Berlin-New York, sixth edition, 1980.

\end{thebibliography}

\def\ocirc#1{\ifmmode\setbox0=\hbox{$#1$}\dimen0=\ht0 \advance\dimen0
  by1pt\rlap{\hbox to\wd0{\hss\raise\dimen0
  \hbox{\hskip.2em$\scriptscriptstyle\circ$}\hss}}#1\else {\accent"17 #1}\fi}
  \def\cprime{$'$}

\end{document}